\documentclass[11pt,final]{article}
\usepackage[pagewise]{lineno}

\usepackage{amsfonts,amsthm,amsmath, amssymb, mathtools, bm}
\usepackage{graphicx,color}
\usepackage{hyperref}
\usepackage{a4wide}
\usepackage{multirow}
\usepackage{enumitem}
\usepackage{mathrsfs}
\usepackage{aliascnt}
\usepackage[nameinlink,capitalise,sort]{cleveref}
\newtheorem{proposition}{Proposition}[section]

\newaliascnt{theorem}{proposition}
\newtheorem{theorem}[theorem]{Theorem}
\aliascntresetthe{theorem}

\newaliascnt{lemma}{proposition}
\newtheorem{lemma}[lemma]{Lemma}
\aliascntresetthe{lemma}

\newaliascnt{corollary}{proposition}

\aliascntresetthe{corollary}

\newaliascnt{definition}{proposition}
\newtheorem{definition}[definition]{Definition}
\aliascntresetthe{definition}

\newaliascnt{remark}{proposition}

\aliascntresetthe{remark}

\newaliascnt{example}{proposition}

\aliascntresetthe{example}

\crefname{definition}{Definition}{Definitions}
\Crefname{definition}{Definition}{Definitions}

\crefname{theorem}{Theorem}{Theorems}
\Crefname{theorem}{Theorem}{Theorems}

\crefname{lemma}{Lemma}{Lemmas}
\Crefname{lemma}{Lemma}{Lemmas}

\crefname{corollary}{Corollary}{Corollaries}
\Crefname{corollary}{Corollary}{Corollaries}

\crefname{remark}{Remark}{Remarks}
\Crefname{remark}{Remark}{Remarks}

\crefname{example}{Example}{Examples}
\Crefname{example}{Example}{Examples}

\newenvironment{proofof}[1]{\smallskip\noindent{\textbf{Proof~of~#1.}}%
  \hspace{1pt}}{\hspace{-5pt}{\nobreak\quad\nobreak\hfill\nobreak%
    $\square$\vspace{2pt}\par}\smallskip\goodbreak}

\numberwithin{equation}{section}

\allowdisplaybreaks[4]

\renewcommand{\phi}{\varphi}
\renewcommand{\theta}{\vartheta}
\renewcommand{\epsilon}{\varepsilon}
\renewcommand{\L}[1]{\mathbf{L^#1}}

\renewcommand{\d}[1]{\mathinner{\mathrm{d}{#1}}}

\newcommand{\dd}[1]{\mathinner{\mathrm{d}{#1}}}
\newcommand{\Lloc}[1]{\mathbf{L^{#1}_{loc}}}
\newcommand{\Wloc}[2]{\mathbf{W^{#1,#2}_{loc}}}
\newcommand{\C}[1]{\mathbf{C^{#1}}}

\newcommand{\Cc}[1]{\mathbf{C_c^{#1}}}

\newcommand{\W}[2]{\mathbf{W^{#1,#2}}}
\newcommand{\BV}{\mathbf{BV}}
\newcommand{\cBV}[1]{\mathcal{BV}^{#1}}
\newcommand{\BVloc}{\mathbf{BV_{loc}}}
\newcommand{\modulo}[1]{{\left|#1\right|}}
\newcommand{\norma}[1]{{\left\|#1\right\|}}
\newcommand{\reali}{{\mathbb{R}}}
\newcommand{\naturali}{{\mathbb{N}}}
\newcommand{\interi}{{\mathbb{Z}}}
\newcommand{\Lip}{\mathop\mathbf{Lip}}
\newcommand{\tv}{\mathop\mathrm{TV}}
\newcommand{\essinf}{\mathop\mathrm{ess~inf}}
\newcommand{\esssup}{\mathop\mathrm{ess~sup}}

\newcommand{\sgn}{\mathop\mathrm{sgn}}
\newcommand{\spt}{\mathop\mathrm{spt}}
\newcommand{\Id}{\mathop{\mathbf{Id}}}
\newcommand{\pint}[1]{\mathaccent23{#1}}

\newcommand{\Su}{\mathbb{S}^1}
\newcommand{\const}[1]{\mathcal{C}_{#1}}

\DeclareFontFamily{U}{mathx}{\hyphenchar\font45}
\DeclareFontShape{U}{mathx}{m}{n}{
  <-6> mathx5 <6-7> mathx6 <7-8> mathx7
  <8-9> mathx8 <9-10> mathx9
  <10-12> mathx10 <12-> mathx12
}{}
\DeclareSymbolFont{mathx}{U}{mathx}{m}{n}
\DeclareMathAccent{\widecheck}{0}{mathx}{"71}

\begin{document}

\title{General Periodic Conservation Laws \\ Motivated by Traffic
  Modeling}

\author{Rinaldo M.~Colombo$^1$ \and Mauro Garavello$^2$ \and Claudia
  Nocita$^2$}

\maketitle

\footnotetext[1]{Unità INdAM \& Dipartimento di Ingegneria
  dell'Informazione, Università di Brescia, Italy.\hfill\\
  \texttt{rinaldo.colombo@unibs.it}}

\footnotetext[2]{Dipartimento di Matematica e Applicazioni,
  Università di Milano-Bicocca, Italy.\hfill\\
  \texttt{mauro.garavello@unimib.it} and
  \texttt{c.nocita@campus.unimib.it}}

\begin{abstract}
  This paper deals with the well-posedness, in the spatially periodic
  setting, of a mixed system of nonlocal-local conservation laws
  modeling traffic flow where autonomous vehicles --- informed on
  traffic conditions over long distances --- interact with
  human-driven bulk traffic.  We also provide a comprehensive
  well-posedness theory in the periodic framework separately for
  general local conservation laws and for systems of nonlocal
  conservation laws, together with \emph{a priori} total variation
  bounds and stability estimates.

  \medskip

    \medskip

  \noindent\textbf{MSC~2020:} 35L65, 35M10, 76A30.
\end{abstract}

\tableofcontents

\section{Introduction}
\label{sec:introduction}

Macroscopic traffic modeling is intimately related to conservation
laws. A circular road provides an optimal environment to observe
\emph{stop \& go waves}. A pioneering work in this context
is~\cite{Sugiyamal2008}, where the onset of \emph{stop \& go waves} is
measured, see also~\cite{Tadaki2013}. Then, the advent of autonomous
vehicles (AVs) leads to search for driving strategies that allow AVs
to mitigate the effect of \emph{stop \& go waves}~\cite{STERN2018205},
in particular reducing fuel consumption and
pollution~\cite{STERN2019351}.

The present work provides the necessary analytical framework to deal
with (systems of) conservation laws in a periodic setting where some
of the vehicles (the AVs) are informed on traffic conditions at long
distances, coexist and interact with the bulk traffic. This work
establishes the well-posedness and stability estimates necessary to
formalize optimal management problems, tackle them and analyze the
long time behavior of solutions.

To get these results, it has been necessary to extend known results in
both local and nonlocal setting. In the local case, while
well-posedness directly follows from the classical
work~\cite{Kruzkov1970}, we provide new $\BV$ bounds and stability
estimates of the solutions on the speed law.  On the other hand, we
develop from scratch the nonlocal periodic framework for general
systems, comprehensive of stability estimates on the speed, which is
only partially present in the literature; see for
example~\cite{MR3959349,MR3854328, MR3057143, keimer_pflug}.

Denote by $r$ the density of human-driven vehicles (which we refer to
as \emph{local population}) and $\rho^i$ the density of the $i$-th
class of AVs (the $i$-th \emph{nonlocal population}), with
$i \in \{1, \ldots, k\}$ and $k \in \naturali \setminus
  \{0\}$. Human-driven vehicles tune their speed according to a
pointwise evaluation of the total density, which at time
$t \in \reali_+$ and position $x \in \Su$ is
\begin{equation}
  \label{eq:60}
  (\Sigma\rho + r) (t,x)
  \coloneqq
  \sum_{i=1}^k \rho^i (t,x) + r (t,x) \,,
\end{equation}
while the $i$-th class of AVs is assigned a speed depending on spatial
convolutions of the total density with smooth compactly supported
kernels $\eta^i \colon \reali \to \reali^{k+1}$, that is
\begin{equation}\label{eq:98}
  \left((\Sigma \rho + r)*\eta^i\right) (t,x)
  \coloneqq
  \sum_{j=1}^k \left( \rho^j (t) *\eta^i_j \right) (x)
  +
  \left( r (t) *\eta^i_{k+1} \right) (x)\,.
\end{equation}
From the modeling point of view, for $j\in \{1, \ldots, k\}$, the
component $\eta_j^i$ represents how the $i$-th \emph{nonlocal}
population reacts to the presence of the $j$-th one and $\eta^i_{k+1}$
describes its interaction with the \emph{local} population. In
particular, $\spt \eta_j^i$ (respectively, $\spt \eta_{k+1}^i$)
encodes the interaction horizon between the $i$-th population and the
$j$-th one (resp. the local population). (See~\eqref{eq:14} for the
definition of the convolution in a periodic setting).

Introducing the speed law $v_L$ of the local population and $v_{NL}^i$
of the $i$-th nonlocal one, it is natural to consider the Cauchy
problem
\begin{equation}
  \label{eq:3}
  \begin{cases}
    \partial_t \rho^i
    +
    \partial_x\left(\rho^i \;
    v_{NL}^i\left(\left(\Sigma\rho + r\right)*\eta^i\right) \right)= 0
     & \text{in } \, \reali_+ \times \Su
    \qquad i \in \{1, \ldots, k\}
    \\
    \partial_t r +
    \partial_x \left(r \;
    v_L\left(\Sigma\rho + r\right)\right) =0
     & \text{in } \, \reali_+ \times \Su
    \\
    \rho(0,x)=\rho_o(x)
     & \text{in } \,  \Su
    \\
    r(0,x)=r_o(x)
     & \text{in } \,  \Su \,,
  \end{cases}
\end{equation}
which extends the model considered in~\cite{opera-seconda} framing it
in the periodic setting of $\Su$. Below, we prove the existence and
uniqueness of solutions, their locally Lipschitz continuous dependence
on the initial data and on time, and various \emph{a priori}
estimates.

\medskip

An alternative approach to the modeling of the AVs -- bulk traffic
interaction has been proposed in~\cite{MR5009512}, where the authors
use nonlocal conservation laws with time delay for the evolution of
the human-driven cars, giving an emphasis on the different reaction
times between AVs and human-driven vehicles.  A different modeling is
obtained through mixed PDE--ODE systems,
see~\cite{ciaramaglia:hal-05547499, MR4002732, MR4104950, MR4720869},
where each AV is singularly described by an ODE. In this setting, the
speed of AVs depends on the bulk density, possibly in a nonlocal way,
at the AV location but does not allow for interaction between
different AVs.  In this context, various control problems related to
traffic management have been considered; see for instance~\cite{bayen,
  MR4442432, MR3623691, MR2801323, MR2852203, MR3264413,
  KeimerPflugRodestock2025, MR2765683, MR4247151}.

In the periodic setting, { especially} with nonlocal terms,
only few results are available.  In~\cite{MR4384589} the authors study
the well-posedness of a nonlocal scalar conservation law on a ring
road where, as in our case, nonlocal terms allow backwards-looking
interactions. Moreover, they show that these interaction kernels, if
suitably designed, are able to stabilize traffic flow.  Finally, in
the setting therein, the asymptotic-in-time $\L2$-convergence of
solutions to the uniform constant profile is proved. Similar
asymptotic results in the periodic setting are also obtained
in~\cite{MR4371088}, in the case of a linear speed law with
anisotropic kernels.  In the local framework, the existence and
uniqueness of entropy periodic solutions to general conservation laws
traces back to~\cite{Kruzkov1970} and, for fluxes depending only on
the unknown,~\cite{MR3072229} proves the asymptotic-in-time
convergence of the entropic solution to the mean value of the initial
datum.

\medskip The first aim of this paper is to fill the gap in the
well-posedness theory in the periodic setting both for local and
nonlocal (systems of) conservation laws with time and space dependent
fluxes.  In particular, we also deduce a full set of \emph{a priori}
$\BV$ estimates for the solutions, and we provide stability estimates
on the speed, extending to the periodic setting the results available
in~\cite{opera-prima, stability_loc}.  While the proof for local
conservation laws largely relies on the finite speed of propagation,
the same argument fails when applied to nonlocal equations due to the
infinite speed of propagation of information induced by the
convolution operator.

The results separately proved at the local and nonlocal levels are
then exploited to establish our main result, i.e. the well-posedness
of the mixed system~\eqref{eq:3}, which is stated in \Cref{thm:WP} for
solutions
$(\rho, r) \in \C0([0,T]; \W11(\Su; \reali^k)\times \L1(\Su;
  \reali))$.  Note that, for $t \in [0, T]$,
$\rho(t) \in \W11(\Su; \reali^k)$ while $r(t) \in \L1(\Su;
  \reali)$. This asymmetry is essentially necessary for including and
applying Kružkov theory in our setting; see also~\cite{opera-seconda}
for a detailed discussion.  Note that in case the nonlocal density
$\rho$ is constantly equal to $0$, then system~\eqref{eq:3} results in
a nonlinear local equation where our theory is completely equivalent
to the Kružkov one. Instead, in the reverse case, i.e. when the local
density $r$ is constantly equal to $0$, then our theory for the mixed
system~\eqref{eq:3} is more restrictive than the general one for
nonlocal systems; see \cref{thm:WP} and \cref{thm:4}.  Moreover, the
lack of a general stability theory for systems of (local) conservation
laws { prevents us} from considering vector-valued density $r(t)$.

The well-posedness of~\eqref{eq:3} is a basic step toward
investigation, in the periodic setting, of a broad class of problems
ranging from control theory and traveling-wave analysis,
see~\cite{traveling}, to the long-time behavior of solutions.

\medskip The next Section presents the analytical results, first in
\S~\ref{sec:mixed-nonlocal-local} about the mixed
system~\eqref{eq:3}. Then, \S~\ref{sec:kruv-zkov-periodic-1} is
devoted to an extension to the periodic case of the classical
well-posedness result~\cite{Kruzkov1970} and of the stability
estimates~\cite{stability_loc}. Similarly, nonlocal systems are
considered in \S~\ref{sec:nonl-peri-case-1}. All proofs are deferred
to Section~\ref{sec:proofs}, while the appendices collect relevant
results not easily available.

\section{Main Result}
\label{sec:main-result}

We represent elements of $\Su$ with real numbers in
$\mathopen[0,1\mathclose[$. Given a map $f \colon \Su \to \reali$ we
denote by $\widehat f \colon \reali \to \reali$ its $1$-periodic
extension to $\reali$. Conversely, given a $1$-periodic function
$g \colon \reali \to \reali$ by $\widecheck{g} \colon \Su \to \reali$
we denote its restriction to $\Su$. We defer the rigorous definitions
of function spaces, norms, and of the convolution in the periodic
setting to Section~\ref{sec:proofs}.

The space $\cBV{}(\Su; \reali^m)$, for $m \in \naturali$, denotes the
space of functions with bounded total variation on $\Su$ and
$\cBV{1}(\Su; \reali^m)$ is the set of functions in
$\W11(\Su; \reali^m)$ whose spatial derivative has bounded total
variation, see~\eqref{eq:20}. We understand that for any function in
$\cBV{} (\Su; \reali^k)$ we use its right continuous representative
and introduce the norms
\begin{displaymath}
  \norma{f}_{\cBV{}(\Su; \reali^k)}
  \coloneqq
  \norma{f}_{\L1(\Su; \reali^k)} + \tv(f) \,,
  \quad\qquad
  \norma{f}_{\cBV{1}(\Su; \reali^k)}
  \coloneqq
  \norma{f}_{\W11(\Su; \reali^k)} + \tv(f') \,.
\end{displaymath}
For a precise definition of Total Variation on $\Su$ we defer
to~\eqref{eq:44}.

We denote by $\mathcal{C} ({\cdots})$ a non-negative continuous
function, non-decreasing in each of its arguments, whose actual value
does not need to be quantified. Moreover, given $\ell > 0$ we set
$B_{\ell} = [-\ell, \ell]^h$ for a suitable integer $h$ depending on
the context.
\subsection{The Mixed Nonlocal-Local Periodic Case}
\label{sec:mixed-nonlocal-local}

Consider the mixed system~\eqref{eq:3} of nonlocal and local
conservation laws on $\Su$. Throughout this section we consider the
following assumptions.
\begin{enumerate}[label=\bf(vNL)]
  \item\label{item:17}
        $v_{NL} \in (\C3\cap \W3\infty)(\reali;\reali^k)$.
\end{enumerate}

\begin{enumerate}[label=\bf(vL)]
  \item\label{item:21} {  $v_L \in \C2(\reali; [0, V_L])$ and
          there exist $V_L, R_L > 0$ such that $v_L(r) = V_L$ for
          $r \leq 0$ and $v_L(r) = 0$ for $r \geq R_L$.}
\end{enumerate}

\begin{enumerate}[label=$\bm{(\eta)}$]
  \item\label{item:hyp-eta}
        $\eta \in \C3(\reali;\reali^{k\times (k+1)})$ with
        $\mathop{\rm diam} (\spt \eta) \leq H$ for a suitable
        $H \in \naturali \setminus \{0\}$.
\end{enumerate}

\noindent We here provide a definition of solution particularized to
the periodic setting. The asymmetry in the regularities of $\rho$ and
$r$ is justified by the nature of~\eqref{eq:3}.
\begin{definition}
  \label{def:sol}
  Let $I = \reali_+$ or $I=[0,T]$ for a positive $T$. The couple
  \begin{displaymath}
    (\rho,r) \in
    \C0\left(I; \W11 (\Su; \reali^k) \times \L1 (\Su; \reali)\right)
  \end{displaymath}
  is a \emph{solution} to~\eqref{eq:3} on $I$ if, setting for
  $(t,x,u) \in I \times \reali \times \reali$ and for
  $i \in \{1, \ldots, k\}$
  \begin{eqnarray*}
    v^i (t,x)
    & \coloneqq
    & v_{NL}^i
    \left(\left(\Sigma\widehat{\rho (t)} * \eta^i\right) (x)
    +
    \left(\widehat{r(t)} *\eta^{i}_{k+1}\right) (x)\right)
    \\
    w (t,x,u)
    & \coloneqq
    & v_L\left(\Sigma\widehat{\rho (t)} \, (x) + u \right),
  \end{eqnarray*}
  then for all $\phi \in \Cc1(\pint{I} \times \reali; \reali)$,
  $\psi \in \Cc1(\pint{I} \times \reali; \reali_+)$,
  $\kappa \in \reali$ and $i \in \{1, \ldots, k\}$
  \begin{eqnarray*}
    \int_I \int_{\reali}
    \left(\widehat\rho^i(t,x) \,\partial_t \phi(t,x) + \widehat\rho^i(t,x) \,
    v^i(t,x) \, \partial_x \phi(t,x) \right) \d{x} \d{t}
    & =
    & 0 \,;
    \\
    \int_I \int_{\reali} \modulo{\widehat{r}(t,x) - \kappa} \, \partial_t \psi(t,x)
    \d{x}\d{t}
    \\
    +
    \int_I \int_{\reali} \sgn(\widehat{r}(t,x) - \kappa) \,
    \left[
      \widehat{r}(t,x) \, w\left(t,x,\widehat{r} (t,x)\right)
      - \kappa \, w(t,x,\kappa)
      \right] \partial_x \psi(t,x) \d{x}\d{t}
    \\
    -
    \int_I \int_{\reali} \sgn(\widehat{r}(t,x) - \kappa) \,
    \kappa \, \partial_x w(t,x,\kappa) \, \psi(t,x) \d{x}\d{t}
    & \geq
    & 0 \,;
    \\
    \rho (0)
    & =
    & \rho_o  \,;
    \\
    r (0)
    & =
    & r_o  \,.
  \end{eqnarray*}
\end{definition}

We present our main result, namely the well-posedness of~\eqref{eq:3}.

\begin{theorem}[Well-Posedness of~\eqref{eq:3}]
  \label{thm:WP}
  Let~\ref{item:17}, \ref{item:21} and~\ref{item:hyp-eta} hold. Then,
  the problem~\eqref{eq:3} generates a unique map
  \begin{equation}
    \label{eq:79}
    \mathcal{S} \colon \reali_+ \times
    \cBV{1}(\Su; \reali_+^k)
    \times
    \cBV{}(\Su; [0,R_L])
    \to
    \cBV{1}(\Su; \reali_+^k)
    \times
    \cBV{}(\Su; [0,R_L])
  \end{equation}
  such that
  \begin{enumerate}[label=\textup{\textbf{($\mathcal{S}$\arabic*)}},
      ref=\textup{\textbf{($\mathbf{\mathcal{S}}$\arabic*)}}]
    \item \label{item:29} $\mathcal{S}$ is a global semigroup, i.e.,
          $\mathcal{S}_0 = \Id$ and for all $t,t_* \in \reali_+$,
          $\mathcal{S}_{t} \circ \mathcal{S}_{t_*} = \mathcal{S}_{t+t_*}$.

    \item \label{item:40} The semigroup $\mathcal{S}$ is locally
          Lipschitz continuous: for all $M,T > 0$, there exists a constant
          \begin{displaymath}
            \const{1} \coloneqq
            \mathcal{C} \left(
            \norma{v_{NL}}_{\W3\infty(\reali; \reali^k)},
            R_L,
            \norma{v_L}_{\W2\infty(\reali; \reali)},
            \norma{\eta}_{\W{3}\infty(\reali; \reali^{k \times (k+1)})},
            H,M, T
            \right)
          \end{displaymath}
          such that for $t, t_* \in [0,T]$ and
          $(\rho_o, r_o), \, (\rho_{o,*}, r_{o,*}) \in \cBV{1}(\Su;
            \reali_+^k) \times \cBV{}(\Su; [0,R_L])$ with
            {
              \begin{displaymath}
                \max
                \left\{
                \norma{\rho_o}_{\cBV{1} (\Su; \reali^k)}
                \,, \;
                \norma{\rho_{o,*}}_{\cBV{1} (\Su; \reali^k)}
                \,, \;
                \norma{r_o}_{\L1 (\Su; \reali)}
                \,, \;
                \norma{r_{o,*}}_{\cBV{} (\Su; \reali)}
                \right\}
                \leq
                M \,,
              \end{displaymath}
            } the following estimates hold:
          \begin{eqnarray*}
            \begin{array}{c}
              \begin{array}{lr}
                 & \norma{ \mathcal{S}_t(\rho_o, r_o)
                -
                \mathcal{S}_t(\rho_{o,*},r_{o,*}) }_{\W11(\Su; \reali^k)
                \times
                \L1(\Su; \reali)}
                \\
                \leq
                 & (1 + \const{1} \, t) \,
                \norma{(\rho_o - \rho_{o,*},
                r_o {-} r_{o,*})}_{\W11(\Su; \reali^k)
                \times
                \L1(\Su; \reali)}\,;
              \end{array}
              \\[18pt]
              \norma{ \mathcal{S}_t(\rho_o, r_o)
                -
                \mathcal{S}_{t_*}(\rho_o, r_o) }_{\W11(\Su; \reali^k)
                \times
                \L1(\Su; \reali)}
              \leq
              \const{1} \, \modulo{t - t_*} \,.
            \end{array}
          \end{eqnarray*}

    \item \label{item:23} For any
          $(\rho_o,r_o) \in \cBV{1}(\Su; \reali_+^k) \times \cBV{}(\Su;
            [0,R_L])$ and any $T > 0$, the orbit
          $t \mapsto \mathcal{S}_t (\rho_o, r_o)$ solves~\eqref{eq:3} on
          $[0, T]$ in the sense of \Cref{def:sol}.

    \item \label{item:39} {$\mathcal{S}$ is uniquely
            characterized by~\ref{item:29}, \ref{item:40} and by the
            property that for all
            $(\rho_o, r_o) \in \C2 (\Su; \reali_+^k) \times \cBV{} (\Su;
              [0,R_L])$, the orbit $t \mapsto \mathcal{S}_t (\rho_o,r_o)$ is
            the unique solution to~\eqref{eq:3} in the sense of
            \Cref{def:sol}.}

    \item \label{item:3} For all
          $(\rho_o, r_o) \in \cBV{1}(\Su; \reali^k_+) \times \cBV{} (\Su;
            [0,R_L])$, call
          $\left(\rho(t), r(t)\right) \coloneqq \mathcal{S}_t(\rho_o, r_o)$
          for $t \in \reali_+$. Then, there exists a positive constant
          \begin{displaymath}
            \const{2} \coloneqq
            \mathcal{C} \left(
            \begin{array}{c}
                \norma{v_{NL}}_{\W3\infty(\reali; \reali^k)},
                R_L,
                \norma{v_L}_{\W2\infty(\reali; \reali)},
                \\
                \norma{\eta}_{\W{3}\infty(\reali; \reali^{k \times (k+1)})},
                H,
                \norma{\rho_o}_{\cBV{1}(\Su; \reali^k)},
                \norma{r_o}_{\cBV{}(\Su; \reali)},
                t
              \end{array}
            \right)
          \end{displaymath}
          such that
          \begin{displaymath}
            \begin{array}{l@{\qquad}l}
              \norma{\rho(t)}_{\L1(\Su; \reali^k)} =
              \norma{\rho_o}_{\L1(\Su; \reali^k)},
               & \norma{r(t)}_{\L1(\Su; \reali )} =
              \norma{r_o}_{\L1(\Su; \reali )},
              \\[.2cm]
              \norma{\rho(t)}_{\cBV{1}(\Su; \reali^k)} \leq
              \const{2},
               & \tv\left(r(t)\right) \leq
              \const{2}.
            \end{array}
          \end{displaymath}
  \end{enumerate}
\end{theorem}

\noindent Note that in~\ref{item:40} the constant $M$ bounds different
norms of $r_o$ and $r_{o,*}$.

Problem~\eqref{eq:3} is autonomous and admits a globally defined
solution. The proof is obtained dealing, separately, with the local
equation and the nonlocal one, both non autonomous. Therefore, in the
next section, we specialize the analysis of these two equations to the
time interval $[0,T]$, for a fixed $T >0$.

\subsection{The {Kru\v zkov} Periodic Case}
\label{sec:kruv-zkov-periodic-1}

We consider the Cauchy problem for a scalar local conservation law on
$\Su$
\begin{equation}
  \label{eq:6}
  \left\{
  \begin{array}{lc@{\,\,}l}
    \partial_t r + \partial_x \left( r \, w(t,x,r) \right) = 0 & \mbox{in} & [0,T] \times \Su
    \\
    r(0,x) = r_o (x)                                           & \mbox{in} & \Su.
  \end{array}
  \right.
\end{equation}
The following assumptions on $w$ are used --- selectively ---
throughout this section.

\begin{enumerate}[label=\bf (K\arabic *), series=requisiti]
  \item \label{item:26}
        $w \in (\C1 \cap \L\infty)([0,T] \times \Su \times \reali;
          \reali_+)$.
  \item \label{item:27}
        $\partial_{rt}^2 w, \partial_{rx}^2 w \in \Lloc\infty([0,T] \times
          \Su \times \reali; \reali)$.
  \item \label{item:34} There exists $R \geq 0$ such that for all
        $r \leq -R$,
        $w(t,x,r) = \norma{w}_{\L\infty ([0,T]\times\Su\times\reali;
          \reali)}$ and for all $r \geq R$, $w(t,x,r) = 0$.
  \item \label{item:28}
        $\partial_{xx}^2 w, \partial_{xr}^2 w \in \C0([0,T] \times \Su
          \times \reali; \reali)$.
  \item \label{item:31}
        $ \partial^2_{rt} w, \partial^2_{tx} w \in
          \L\infty([0,T]\times\Su\times\reali; \reali)$.
\end{enumerate}

\noindent To tackle problem~\eqref{eq:6}, we first provide the
definition of solution and then the corresponding well-posedness and
stability result.

\begin{definition}\label{def:solLoc}
  Let~\ref{item:26} hold and $r_o \in \L\infty(\Su; \reali)$. By
  \emph{solution} to~\eqref{eq:6} we mean a map
  $r \in \L\infty([0,T]\times \Su; \reali)$ such that $\widehat r$ is
  a {Kru\v zkov} solution (in the sense
  of~\cite[Definition~1]{Kruzkov1970}) to the Cauchy problem
  \begin{equation}\label{eq:38}
    \left\{
    \begin{array}{lr@{\,\,}l}
      \partial_t \widehat r + \partial_x \left( \widehat r \, \widehat w(t,x,\widehat r) \right) = 0
       &
      \mbox{in}
       &
      [0,T] \times \reali
      \\
      \widehat r(0,x) = \widehat r_o (x)
       &
      \mbox{in}
       &
      \reali  \,,
    \end{array}
    \right.
  \end{equation}
  where $\widehat w$ is the spatial periodic extension of $w$.
\end{definition}

\begin{theorem}[Well-posedness and stability of~\eqref{eq:6}]
  \label{thm:K}
  Assume that $w$ satisfies~\ref{item:26}--\ref{item:27}.
  \begin{enumerate}[label=\bf(\roman*)]
    \item \label{item:16} For all
          $r_o, r_{o,*} \in \L\infty(\Su; \reali)$, if
          $r, r_* \in \L\infty([0,T] \times \Su; \reali)$ are the
          corresponding solutions to~\eqref{eq:6} in the sense of
          \Cref{def:solLoc}, then for a.e. $t \in [0,T]$,
          \begin{displaymath}
            \norma{r(t,x) - r_*(t,x)}_{\L1(\Su; \reali)}
            \leq \norma{r_o - r_{o,*}}_{\L1(\Su; \reali)}.
          \end{displaymath}

    \item \label{item:18} If also~\ref{item:34} holds, for all
          $r_o \in \L\infty(\Su; [0,R])$, if $r$ is the solution
          to~\eqref{eq:6} in the sense of \Cref{def:solLoc}, then for a.e.
          $t \in [0,T]$,
          \begin{displaymath}
            r(t,x) \in [0,R] \mbox{ for a.e. } x \in \Su
            \qquad
            \mbox{and}
            \qquad
            \norma{r(t)}_{\L1(\Su; \reali)} = \norma{r_o}_{\L1(\Su; \reali)} \,.
          \end{displaymath}

    \item \label{item:19} If also~\ref{item:34} and~\ref{item:28} hold,
          for all $r_o \in \L\infty(\Su; \reali)$, problem~\eqref{eq:6}
          admits a solution in the sense of Definition~\ref{def:solLoc}.

    \item \label{item:20} If also~\ref{item:34}$\cdots$\ref{item:31}
          hold, for all $r_o \in \cBV{} \left(\Su; [0,R]\right)$, the
          solution $r$ to~\eqref{eq:6} in the sense of
          Definition~\ref{def:solLoc} is such that for a.e.~$t \in [0,T]$,
          $r(t) \in \cBV{}(\Su; [0,R])$ and
          \begin{equation}
            \label{eq:50}
            \tv\left(r(t)\right)
            \leq
            \tv(r_o) \,
            e^{\kappa_o t}
            +
            R \,
            \sup_{\tau \in [0,t]}
            \norma{\partial_{xx}^2 w(\tau)}_{\L1(\Su; \L\infty(\reali; \reali))}
            \dfrac{e^{\kappa_o t} -1}{\kappa_o}
          \end{equation}
          where
          \begin{equation}
            \label{eq:1}
            \kappa_o
            \coloneqq
            \mathcal{C}
            \left(
            \norma{w}%
            _{\L\infty([0,T] \times \Su; \W1\infty(\reali; \reali))},
            \norma{\partial_x w}%
            _{\L\infty([0,T] \times \Su; \W1\infty(\reali; \reali))},
            R
            \right).
          \end{equation}

    \item \label{item:24} If $w, w_*$
          satisfy~\ref{item:26}$\cdots$\ref{item:28} and $w$ also
          satisfies~\ref{item:31}, for all $r_o \in \cBV{} (\Su; [0,R])$ the
          problems
          \begin{equation}
            \label{eq:37}
            \left\{
            \begin{array}{l}
              \partial_t r + \partial_x \left(r \; w (t,x,r)\right) =0
              \\
              r (0,x) = r_o (x)
            \end{array}
            \right.
            \qquad \mbox{ and } \qquad
            \left\{
            \begin{array}{l}
              \partial_t r_*
              +
              \partial_x \left(r_* \; w_* (t,x,r_*)\right) =0
              \\
              r_* (0,x) = r_o (x)
            \end{array}
            \right.
          \end{equation}
          admit the solutions $r$ and $r_*$ in the sense of
          \Cref{def:solLoc} that satisfy
          \begin{align}
            \label{eq:4}
             & \quad  \norma{r (t) {-} r_* (t)}_{\L1 (\Su; \reali)}
            \\
            \nonumber
             & \leq
            t\,  e^{\kappa_o t} \,
            \mathcal{C}
            \left(
            \begin{array}{c}
                \norma{w}_{\L\infty([0,T]; \W1\infty(\Su \times \reali; \reali))},
                \sup\limits_{\tau \in [0,t]}
                \norma{\partial_{xx}^2 w(\tau)}_{\L1(\Su; \L\infty(\reali; \reali ))},
                \\
                { \norma{w_*}_{\L\infty([0,T] \times \Su ; \W1\infty( \reali; \reali))},
                }
                R, \tv(r_o), T
              \end{array}
            \right)
            \\
            \nonumber
             & \quad {\times}
            \left( \!
            \norma{w {-} w_*}_{\L\infty([0,T]\times \Su; \W1\infty(\reali; \reali))}
            +
            \sup_{t \in [0,T]} \int_{\Su} \sup_{q \in \reali}
            \modulo{
              \partial_x w (t, x, q) - \partial_x w_*(t, x, q)
            } \d x
            \! \right)
          \end{align}

    \item \label{item:25} If $w$
          satisfies~\ref{item:26}$\cdots$\ref{item:31} and
          $r_o \in \cBV{}(\Su; [0,R])$, then, the solution $r$
          to~\eqref{eq:6} in the sense of Definition~\ref{def:solLoc} is
          Lipschitz continuous in time: for all $t,t_* \in [0,T]$,
          \begin{displaymath}
            \norma{r(t) - r(t_*)}_{\L1(\Su; \reali)
            }
            \leq
            \mathcal{C}
            \left(
            \begin{array}{c}
                \norma{w}_{\L\infty([0,T]\times \Su; \W1\infty(\reali; \reali ))}\,,
                \\
                \norma{\partial_x w}_{\L\infty([0,T]\times \Su; \W1\infty(\reali; \reali ))}\,,
                \\
                \sup\limits_{\tau \in [0,T]}
                \norma{\partial_{xx}^2 w(\tau)}_{\L1(\Su; \L\infty(\reali; \reali ))} \,,
                \\
                R,
                \tv(r_o), T
              \end{array}
            \right) |t-t_*| \,.
          \end{displaymath}

  \end{enumerate}
\end{theorem}

\subsection{The Nonlocal Periodic Case}
\label{sec:nonl-peri-case-1}

For a fixed $k \in \naturali \setminus \{0\}$ and
$\eta \colon \reali \to \reali^{k \times k}$,
$\eta = \left(\eta^i_j\right)$ where $i,j \in \{ 1, \ldots k\}$, we
use the notation
\begin{equation}
  \label{eq:27}
  \rho *\eta ^i
  \coloneqq
  \left(\rho^1 * \eta^i_1, \ldots, \rho^k * \eta^i_k \right)
\end{equation}
where each convolution is computed as in~\eqref{eq:14}.  Note that,
compared to Section~\ref{sec:mixed-nonlocal-local}, the kernel $\eta$
here does not have an additional column $\eta_{k+1}$ which encodes the
interaction with the \emph{local} population, here absent.

With reference to the following system of nonlocal conservation laws
on $\Su$
\begin{equation}\label{eq:25}
  \left\{
  \begin{array}{ll}
    \partial_t \rho^i
    +
    \partial_x\left(\rho^i \, V^i(t,x,\rho*\eta^i)\right) = 0
     & \mbox{in } [0,T]\times \Su \quad i \in \{1, \ldots, k\}
    \\
    \rho(t_o) =  \rho_o
     & \mbox{in } \Su \,,
  \end{array}
  \right.
\end{equation}
for $t_o \in [0,T]$, we provide the definition of solution of use in
the sequel.

\begin{definition}
  \label[definition]{def:solNonloc}
  Let $t_o \in [0,T]$,
  $V \in \Lloc\infty([0,T]\times \Su \times \reali^k; \reali^k)$ and
  $\eta \in \L1(\reali; \reali^{k \times k})$.  By \emph{solution
    to}~\eqref{eq:25} on $[0,T]$ we mean a map
  $\rho \in \C0\left([0,T]; \L1(\Su; \reali^k)\right)$ such that, for
  any $i \in \{1, \ldots, k\}$, setting
  $v^i(t,x) \coloneqq \widehat{V}^i\left(t,x,
    \left(\widehat{\rho}(t)*\eta^i\right)(x)\right)$,
  $\widehat{\rho}^i$ is a weak solution to
  \begin{equation}
    \label{eq:2}
    \left\{
    \begin{array}{ll}
      \partial_t \widehat{\rho}^i +
      \partial_x\left(\widehat{\rho}^i \, v^i(t,x)\right) = 0
       & \mbox{in } [0,T]\times \reali
      \\
      \widehat{\rho}^i(t_o) =  \widehat{\rho}_o^i \,.
       & \mbox{in }\reali \,.
    \end{array}
    \right.
  \end{equation}
\end{definition}

\noindent The following theorem extends to $\Su$ various results in
the literature, see for instance~\cite{opera-prima, MR3057143,
  keimer_pflug}.

\begin{theorem}[Well-posedness and stability of~\eqref{eq:25}]
  \label{thm:4}
  Assume $t_o \in [0,T]$ and
  \begin{enumerate}[label=\bf(N1)]
    \item \label{item:32}
          $V \in \C0\left([0,T]\times \Su \times \reali^k; \reali^k\right)
            \cap \L\infty\left([0,T]; \Wloc2\infty(\Su \times \reali^k;
            \reali^k)\right)$.
  \end{enumerate}
  \begin{enumerate}[label=\bf($\eta1$)]
    \item \label{item:41}
          $\eta \in \W2\infty(\reali;\reali^{k\times k})$,
          $\mathrm{diam}(\spt \eta)\leq H$ for a suitable
          $H \in \naturali \setminus \{0\}$.
  \end{enumerate}
  Then there exists a map
  \begin{equation}
    \label{eq:10}
    \mathcal{P} \colon [0,T]\times [0,T] \times  \cBV{}(\Su;\reali^k) \to \cBV{}(\Su;\reali^k)
  \end{equation}
  such that, for any $\mathrm{M} >0$ and setting
  $\mathrm{M}_* \coloneqq H \, \norma{\eta}_{\L\infty(\reali;
    \reali^{k\times k})} \, \mathrm{M}$, we have
  \begin{enumerate}[label=\bf(\roman*)]

    \item \label{item:4} For $\rho_o \in \cBV{}(\Su;\reali^k)$, the map
          $t \mapsto \mathcal{P}_{t_o,t}\rho_o$ is the unique solution, in
          the sense of \Cref{def:solNonloc}, on $[0,T]$ to~\eqref{eq:25}
          with initial datum $\rho_o$ assigned at time $t_o$.

    \item \label{item:10} For all $t_1,t_2,t_3 \in [0,T]$,
          $\mathcal{P}_{t_1,t_3} = \mathcal{P}_{t_2,t_3} \circ
            \mathcal{P}_{t_1,t_2}$ and $\mathcal{P}_{t_1,t_1} = \Id$.

    \item \label{item:12} Fix $i \in \{1, \ldots, k\}$.  If
          $\rho_o \in \cBV{}(\Su; \reali^k)$ and $\rho^i_o \geq 0$, then
          $(\mathcal{P}_{t_o,t}\rho_o)^i \geq 0 $ for all $t \in [0,T]$.

    \item \label{item:35} There exists
          \begin{displaymath}
            \const{3}
            \coloneqq
            \mathcal{C} \left(
            \norma{V}_{\L\infty([0,T]; \W2\infty(\Su \times B_{\mathrm{M}_*}; \reali^k))},
            \norma{\eta}_{\W2\infty(\reali; \reali^{k\times k})}, H, M \right)
          \end{displaymath}
          such that, for all $\rho_o \in \cBV{}(\Su; \reali^k)$ with
          $\norma{\rho_o}_{\L1 (\Su; \reali^k)} \leq M$, for all
          $t \in [0,T]$, for $i \in \{1,\ldots,k\}$
          \begin{eqnarray*}
            \norma{(\mathcal{P}_{t_o,t} \rho_o)^i}_{\L1(\Su; \reali)}
            & =
            & \norma{\rho_o^i}_{\L1(\Su; \reali)}
            \\
            \norma{(\mathcal{P}_{t_o,t} \rho_o)^i}_{\L\infty(\Su; \reali)}
            & \leq
            & \norma{\rho_o^i}_{\L\infty(\Su; \reali)}
            e^{\const{3} \, \modulo{t-t_o}}
            \\
            \tv\left((\mathcal{P}_{t_o,t} \rho_o)^i\right)
            & \leq
            & \left(
            \tv(\rho_o^i)
            +
            \const{3} \, \norma{\rho_o^i}_{\L1(\Su; \reali)} \, \modulo{t-t_o}
            \right)
            e^{\const{3} \, \modulo{t-t_o}}.
          \end{eqnarray*}

    \item \label{item:38} For all
          $\rho_o, \rho_{o,*} \in \cBV{}(\Su; \reali^k)$ with
          $\norma{\rho_o}_{\L1(\Su; \reali^k)}, \norma{\rho_{o,*}}_{\L1(\Su;
              \reali^k)} \leq \mathrm{M}$, then, for all $t \in [0,T]$ there
          exists
          \begin{displaymath}
            \const{4} \coloneqq
            \mathcal{C} \left(
            \begin{array}{c}
                \norma{V}_{\L\infty([0,T]; \W2\infty(\Su \times B_{\mathrm{M}_*}; \reali^k))},
                \norma{\eta}_{\W2\infty(\reali; \reali^{k \times k})},
                H,
                \\
                \norma{\rho_o}_{\cBV{}(\Su; \reali^k)},
                \norma{\rho_{o,*}}_{\L1(\Su; \reali^k)},
                T
              \end{array}
            \right)
          \end{displaymath}
          such that
          \begin{displaymath}
            \norma{\mathcal{P}_{t_o,t}\rho_o - \mathcal{P}_{t_o,t}\rho_{o,*}}_{\L1(\Su; \reali^k)}
            \leq
            \norma{\rho_o - \rho_{o,*}}_{\L1(\Su; \reali^k)}
            e^{\const{4} \, \modulo{t-t_o}}\,.
          \end{displaymath}
    \item \label{item:43} If $V, V_*$ satisfy~\ref{item:32}, for all
          $\rho_o \in \cBV{}(\Su; \reali^k)$, the problems
          \begin{equation}
            \label{eq:19}
            \left\{
            \begin{array}{@{}l}
              \partial_t \rho^i +
              \partial_x\left(\rho^i \, V^i(t,x,\rho*\eta^i)\right) = 0
              \\
              \rho(t_o) =  \rho_o
            \end{array}
            \right.
            \qquad
            \left\{
            \begin{array}{@{}l}
              \partial_t \rho_*^i +
              \partial_x\left(\rho_*^i \, V_*^i(t,x,\rho_* * \eta^i)\right) = 0
              \\
              \rho_*(t_o) =  \rho_o
            \end{array}
            \right.
          \end{equation}
          generate the processes $\mathcal{P}$ and $\mathcal{P}^*$ that
          satisfy, for all $t \in [0,T]$,
          \begin{equation}
            \label{eq:24}
            \begin{array}{rcl}
              \norma{\mathcal{P}_{t_o,t}\rho_o - \mathcal{P}^*_{t_o,t}\rho_o}_{\L1(\Su; \reali^k)}
               & \leq
               & \mathcal{C}
              \left(
              \begin{array}{@{}c@{}}
                  \norma{V}_{\L\infty([0,T]; \W2\infty(\Su \times B_{\mathrm{M}_*}; \reali^k))},
                  \\
                  \norma{V_*}_{\L\infty([0,T]; \W1\infty(\Su \times B_{\mathrm{M}_*}; \reali^k))},
                  \\
                  \norma{\eta}_{\W2\infty(\reali; \reali^{k \times k})},
                  H,
                  \norma{\rho_o}_{\cBV{}(\Su; \reali^k)}, T
                \end{array}
              \right)
              \\[24pt]
               &
               & \times
              \norma{V-V_*}_{\L\infty([0,T]; \W1\infty(\Su \times B_{\mathrm{M}_*}; \reali^{k}))} \; \modulo{t-t_o}.
            \end{array}
          \end{equation}

    \item \label{item:44} For all $\rho_o \in \cBV{} (\Su; \reali^k)$,
          the map $t \mapsto \mathcal{P}_{t_o,t}\rho_o$ is locally
          $\L1$--Lipschitz continuous in time: for all $t, t_* \in [0,T]$
          \begin{displaymath}
            \norma{\mathcal{P}_{t_o,t}\rho_o - \mathcal{P}_{t_o,t_*}\rho_o}_{\L1(\Su; \reali^k)}
            \leq
            \mathcal{C}
            \left(
            \begin{array}{c}
                \norma{V}_{\L\infty([0,T]; \W2\infty(\Su \times B_{\mathrm{M}_*}; \reali^k))},
                \\
                \norma{\eta}_{\W2\infty(\reali; \reali^{k \times k})},
                H,
                \norma{\rho_o}_{\cBV{}(\Su; \reali^k)},
                T
              \end{array}
            \right)
            \modulo{t- t_*} \,.
          \end{displaymath}
  \end{enumerate}
\end{theorem}

\begin{proposition}[Fine estimates on the solution to~\eqref{eq:25}]
  \label{prop:fine_estimates}
  Assume $t_o \in [0,T]$ and
  \begin{enumerate}[label=\bf(N2)]
    \item \label{eq:87}
          $V \in \C0\left([0,T]; \C1(\Su \times \reali^k; \reali^k)\right)
            \cap \L\infty\left([0,T]; \Wloc3\infty(\Su \times \reali^k;
            \reali^k)\right)$,
          $ V(t,\cdot, \cdot) \in \C3(\Su \times \reali^k; \reali^k)$ for
          all $t \in [0,T]$.
  \end{enumerate}
  \begin{enumerate}[label=\bf($\eta2$)]
    \item \label{item:42} $\eta \in \C3(\reali;\reali^{k\times k})$,
          $\mathrm{diam}(\spt \eta)\leq H$ for a suitable
          $H \in \naturali \setminus \{0\}$.
  \end{enumerate}
  Let $\rho_o \in \C2(\Su;\reali^k)$ and set
  $\mathrm{M}_* \coloneqq H \, \norma{\eta}_{\L\infty(\reali;
    \reali^{k\times k})}\, \norma{\rho_o}_{\L1(\Su; \reali^k)}$. Then,
  Theorem~\ref{thm:4} applies and the solution
  $\rho (t) = \mathcal{P}_{t_o,t}\rho_o$ to~\eqref{eq:25} satisfies
  \begin{enumerate}[label=\bf(\roman*)]
    \item \label{item:30} $\rho \in \C1([0,T] \times \Su; \reali^k)$,
          $\partial_{xx}^2\rho \in \C0([0,T]\times \Su; \reali^k)$,
          $\partial_{tx}^2\rho \in \L\infty([0,T]\times \Su; \reali^k)$ and
          for all $ t \in [t_o,T]$, the case $t \in [0, t_o]$ being
          analogous,
          \begin{eqnarray*}
            \norma{\rho}_{\L\infty([0,T]; \W\ell{\infty}(\Su; \reali^k))}
            & \leq
            & \mathcal{C} \left(
            \begin{array}{c}
                \norma{V}_{\L\infty([0,T]; \W{\ell+1}\infty(\Su \times B_{\mathrm{M}_*}; \reali^k))},
                \\
                \norma{\eta}_{\W{\ell+1}\infty(\reali; \reali^{k \times k})},
                H,
                \norma{\rho_o}_{\W\ell\infty(\Su; \reali^k)},
                T
              \end{array}
            \right) \quad \ell = 1,2 \,
            \\
            \norma{\rho}_{\L\infty([0,T]; \W21(\Su; \reali^k))}
            & \leq
            & \mathcal{C} \left(
            \begin{array}{c}
                \norma{V}_{\L\infty([0,T]; \W3\infty(\Su \times B_{\mathrm{M}_*}; \reali^k))},
                \\
                \norma{\eta}_{\W{3}\infty(\reali; \reali^{k \times k})},
                H,
                \norma{\rho_o}_{\W21(\Su; \reali^k)},
                T
              \end{array}
            \right).
          \end{eqnarray*}

    \item \label{item:36} $\rho$ is $\W11$ locally Lipschitz continuous
          in time, i.e. for all $t, t_* \in [0,T]$
          \begin{displaymath}
            \norma{\rho(t)- \rho(t_*)}_{\W11(\Su; \reali^k)}
            \leq
            \mathcal{C}
            \left(
            \begin{array}{c}
                \norma{V}_{\L\infty([0,T]; \W3\infty(\Su \times B_{\mathrm{M}_*}; \reali^k))},
                \\
                \norma{\eta}_{\W3\infty(\reali; \reali^{k \times k})},
                H,
                \norma{\rho_o}_{\cBV{1}(\Su; \reali^k)},
                T
              \end{array}
            \right)
            \modulo{t- t_*} \,.
          \end{displaymath}
  \end{enumerate}
\end{proposition}

\section{Proofs}
\label{sec:proofs}

For all $h \in \naturali$, $p \in [1,\infty]$ and
$m,n \in \naturali \setminus \{0\}$, we identify the spaces
\begin{equation}
  \label{eq:20}
  \small
  \!\!\!
  \begin{array}{@{\!\!\!}r@{\,}c@{\,}l@{\;}r@{\,}c@{\,}l@{\!\!}}
    \L{p} (\Su; \reali^m)
     & \coloneqq
     & \left\{
    f \colon \Su \to \reali^m \colon
    \widehat{f} \in \Lloc{p} (\reali; \reali^m)
    \right\}
     & \norma{f}_{\L{p} (\Su; \reali^m)}
     & \coloneqq
     & \norma{\smash{\widehat{f}}}_{\L{p} ([0,1];\reali^m)}
    \\
    \!
    \W{h}{p}(\Su; \reali^m)
     & \coloneqq
     & \left\{
    f \colon \Su \to \reali^m \colon
    \widehat{f} \in \Wloc{h}{p} (\reali; \reali^m)
    \right\}
     & \norma{f}_{\W{h}{p} (\Su; \reali^m)}
     & \coloneqq
     & \norma{\smash{\widehat{f}}}_{\W{h}{p} ([0,1];\reali^m)}
    \\
    \C{h} (\Su; \reali^m)
     & \coloneqq
     & \left\{
    f \colon \Su \to \reali^m \colon
    \widehat{f} \in \C{h} (\reali; \reali^m)
    \right\}
     & \norma{f}_{\C{h} (\Su; \reali^m)}
     & \coloneqq
     & \norma{\smash{\widehat{f}}}_{\C{h} ([0,1];\reali^m)}
    \\
    \cBV{} (\Su; \reali^m)
     & \coloneqq
     & \left\{
    f \colon \Su \to \reali^m \colon
    \widehat{f} \in \BVloc (\reali; \reali^m)
    \right\}
     & \norma{f}_{\cBV{} (\Su; \reali^m)}
     & \coloneqq
     & \norma{f}_{\L1 (\Su;\reali^m)} {+} \tv (f)
    \\
    \cBV{1} (\Su; \reali^m)
     & \coloneqq
     & \left\{
    f {\in} \W11 (\Su; \reali^m) \colon
    f' \in \cBV{} (\Su; \reali^m)
    \right\}
     & \norma{f}_{\cBV{1}(\Su; \reali^m)}
     & \coloneqq
     & \norma{f}_{\W11(\Su; \reali^m)} {+} \tv(f')
  \end{array}
\end{equation}
If $X,Y$ are any of the spaces above and $I$ is a real interval, we
set
\begin{equation}
  \label{eq:32}
  \norma{(f_1,f_2)}_{X\times Y}
  \coloneq
  \norma{f_1}_X + \norma{f_2}_Y
  \quad \mbox{ and } \quad
  \begin{array}{r@{\,}c@{\,}l}
    \norma{f}_{\C0 (I; X)}
     & \coloneq
     & \sup_{t \in I} \norma{f (t)}_X
    \\[6pt]
    \norma{f}_{\L\infty (I; X)}
     & \coloneq
     & \esssup_{t \in I} \norma{f (t)}_X \,.
  \end{array}
\end{equation}
For $f \in \cBV{} (\Su; \reali^m)$, define
\begin{equation}
  \label{eq:44}
  \tv(f)
  \coloneqq
  \tv (\widehat{f}; \mathopen ]0,1]) \,.
\end{equation}
In other words, to compute $\tv(f)$ for
$f \in \cBV{} (\Su; \reali^m)$, first extend $f$ to $\widehat{f}$,
then compute its distributional derivative $\mu$, which is a measure
on $\reali$, pass to its total variation $\norma{\mu}$ and obtain
$\tv(f) = \norma{\mu} (\mathopen]0,1])$.

Given $f \in \W{1}{p}(\Su; \reali^k)$, the weak derivative
$\partial_x f \in \L{p}(\Su;\reali^k)$ is such that
\begin{equation}\label{eq:41}
  \widehat{\partial_x f } (x) = \partial_x \widehat{f} (x) \
  \quad \mbox{ for } x \in \reali,
  \qquad
  \partial_x f (\xi)=
  \widecheck{\left(\partial_x \widehat{f}\right)}(\xi)
  \quad \mbox{ for } \xi \in \Su.
\end{equation}
Denote the integer part of a real number $x$ by
$\lfloor x \rfloor \coloneq \max\{k \in \interi \colon k\leq x\}$.
  {
    For a function $f \colon \Su \to \Su$, we define
    $\widetilde{f}\colon \reali \to \reali$
    a lift of $f$ to $\reali$, i.e.
    \begin{equation}
      \label{eq:16}
      \pi \circ \widetilde{f} = f \circ \pi,
    \end{equation}
    with $\pi \colon \reali \to \Su$ is the quotient map
    $\pi(x) = x - \lfloor x \rfloor$.  }

Given { either
    $ f \in \L\infty(\Su; \reali), \eta \in \L1(\reali;
      \reali)$ or
    $ f \in \L1(\Su; \reali), \eta \in \L1(\reali; \reali)$
    and $\eta$ compactly supported, } define the convolution product
$f * \eta$ in $\Su$ as
\begin{equation}
  \label{eq:14}
  f * \eta \coloneqq \widecheck{(\widehat f *
    \eta)}
  \qquad \mbox{ or, equivalently, for } x \in \reali
  \quad
  \widehat{
    (f * \eta)
  }(x)
  =
  {\int_{\reali} \widehat f (x-y) \; \eta (y) \d{y}} \,.
\end{equation}
See Lemma~\ref{lem:0} for the basic properties of this periodic
convolution.

\subsection{Proofs Related to the {Kru\v zkov} Periodic Case}
\label{sec:kruv-zkov-periodic}

\begin{proofof}{Theorem~\ref{thm:K}}
  \paragraph{Proof of~\ref{item:16}.}
  Clearly,
  $\widehat r_o, \widehat{r}_{o,*} \in \L\infty(\reali; \reali)$ and,
  as a direct consequence of definitions~\eqref{eq:20}, $\widehat w$
  satisfies~\cite[\textbf{(L1)} and~\textbf{(L2)} of
    Lemma~3.8]{opera-seconda}. Setting
  $R_o \coloneqq \max\{\norma{r_o}_{\L\infty(\Su; \reali)},
    \norma{r_{o,*}}_{\L\infty(\Su; \reali)}\}$, $\widehat{w}$ also
  satisfies~\cite[\textbf{(L3)} in~Lemma 3.8]{opera-seconda} since,
  by~\ref{item:26},
  \begin{displaymath}
    \sup
    \left\{
    \modulo{\partial_r \widehat{w}(\tau, \xi, q)}
    \colon
    \tau \in [0,T],
    \xi \in \reali,
    q \in B_{R_o}
    \right\}
    \leq
    \norma{\partial_r w}_{\L\infty([0,T] \times \Su \times B_{R_o}; \reali)}
    < + \infty \,.
  \end{displaymath}
  Setting
  \begin{displaymath}
    \sigma
    \coloneqq
    \norma{w}_{\L\infty([0,T] \times \Su \times B_{R_o}; \reali)}
    +
    R_o \,
    \norma{\partial_r w}_{\L\infty([0,T] \times \Su \times B_{R_o} ; \reali)} \,,
  \end{displaymath}
  \cite[eqn.~(3.43) in~Lemma~3.8]{opera-seconda} ensures that for all
  $t \in [0,T]$ and $a, b \in \reali$ with $a\leq b$,
  \begin{equation}
    \label{eq:17}
    \norma{\widehat{r}(t) - \widehat{r}_*(t)}%
    _{\L1([a,b]; \reali)}
    \leq
    \norma{\widehat{r}_o - \widehat{r}_{o,*}}%
    _{\L1([a-\sigma t, b+\sigma t]; \reali)}.
  \end{equation}
  Hence, for any $\mathrm{N} \in \naturali$,
  \begin{flalign*}
    2\mathrm{N} \, \norma{r(t) - r_*(t)}_{\L1(\Su; \reali)}
     & =
    \norma{\widehat r(t) - \widehat{r}_*(t)}_{\L1(B_{\mathrm{N}}; \reali)}
     & [\mbox{By~\eqref{eq:20}}]
    \\
     & \leq
    \norma{\widehat r_o - \widehat{r}_{o,*}}_{\L1(B_{\mathrm{N}+\sigma t}; \reali)}
     & [\mbox{By~\eqref{eq:17}}]
    \\
     & =
    2\mathrm{N} \, \norma{r_o - r_{o,*}}_{\L1(\Su; \reali)}
    +
    \norma{\widehat r_o - \widehat{r}_{o,*}}_{\L1(B_{\mathrm{N}+\sigma t} \setminus B_\mathrm{N}; \reali)} \,.
     & [\mbox{By~\eqref{eq:20}}]
  \end{flalign*}
  Divide by $2\mathrm{N}$ and let $\mathrm{N} \to +\infty$ to conclude
  the proof of~\ref{item:16}.

  \paragraph{Proof of~\ref{item:18}.}
  This is a direct consequence of~\cite[Theorem~3]{Kruzkov1970}, see
  also~\cite[Lemma~3.9]{opera-seconda}.

  \paragraph{Proof of~\ref{item:19}.}
  The functions $\widehat w$ and $\widehat{r}_o$ meet the requirements
  of~\cite[Lemma 3.10]{opera-seconda}.  In particular, due
  to~\ref{item:26},~\ref{item:34} and~\ref{item:28},
  $\partial_x\widehat w, \partial_r\widehat w,\partial_{xr}^2 \widehat
    w \in \L\infty([0,T] \times \reali \times \reali; \reali)$.  Hence,
  \cite[Theorem~5]{Kruzkov1970} ensures the existence of a solution
  to~\eqref{eq:38} in the sense of~\cite[Definition~1]{Kruzkov1970}.
  We are left to verify that $\widehat{r}$ is 1-periodic in space. In
  fact, the function $\mathtt{r}(t,x) \coloneqq \widehat r(t, x+1)$ is
  the {Kru\v zkov} solution to the Cauchy problem
  \begin{displaymath}
    \left\{
    \begin{array}{ll}
      \partial_t \mathtt{r} + \partial_x \left( \mathtt{r} \, \widehat w(t,x+1,\mathtt{r}) \right) = 0
       &
      \mbox{in}\,
      [0,T] \times \reali
      \\
      \mathtt{r}(0,x) = \widehat r_o (x+1)
       &
      \mbox{in}\,
      \reali
    \end{array}
    \right.
  \end{displaymath}
  which, due to the periodicity of $\widehat r_o$ and $\widehat w$,
  also solves~\eqref{eq:38}. By the classical {Kru\v zkov} uniqueness
  result~\cite[Theorem~2]{Kruzkov1970}, we conclude that
  $\widehat r(t,x) = \mathtt{r}(t,x) = \widehat r(t,x+1)$ for
  a.e.~$(t,x) \in [0,T] \times \reali$, concluding the proof
  of~\ref{item:19}.

  \paragraph{Proof of~\ref{item:20}.}
  We split the proof in two separate steps.\\
  \noindent
  \textbf{\underline{Step~1}: Construction of Approximate Solutions.}
  \\
  Introduce, for $\mathrm{N} \in \naturali \setminus\{0\}$,
  \begin{equation}
    \label{eq:26}
    \begin{array}{rcl}
      \sigma
       & \coloneqq
       & \norma{w}_{\L\infty([0,T] \times \Su \times B_R; \reali)}
      +
      R \,
      \norma{\partial_r w}_{\L\infty([0,T] \times \Su \times B_R; \reali)} \,;
      \\
      \mathcal{K}
       & \coloneqq
       & \left\{
      (t,x) \in [0,T] \times \reali \colon
      x \in
      B_{\mathrm{N} + \sigma (T-t)}
      \right\} \,;
      \\
      \mathtt{r}_{o,\mathrm{N}} (x)
       & \coloneqq
       & \left\{
      \begin{array}{lc}
        \widehat{r}_o (x)
         & x \in B_{\mathrm{N} + \sigma T},
        \\
        0
         & x \, \mbox{elsewhere}.
      \end{array}
      \right.
    \end{array}
  \end{equation}
  Fix $\delta \in \Cc\infty(\reali; \reali_+)$ so that
  $\norma{\delta}_{\L{1}\left(\reali; \reali\right)} = 1$ and
  $\spt(\delta) \subset B_1$ and consider the map
  \begin{equation}
    \label{eq:5}
    \zeta_{\mathrm{N}}
    \coloneqq
    \delta * \mathbf{1}_{B_{\mathrm{N} + \sigma T +  2}}
  \end{equation}
  where $\mathbf{1}_{B_{\mathrm{N} + \sigma T + 2}}$ is the
  characteristic function of the interval
  $B_{\mathrm{N} + \sigma T + 2}$.  Observe that
  \begin{equation}
    \label{eq:23}
    \zeta_{\mathrm{N}} \in \Cc\infty(\reali; \reali)
    \quad \mbox{ and } \quad
    \left\{
    \begin{array}{l}
      \zeta_\mathrm{N} (x) = 1
      \quad \mbox{ for all } x \in B_{\mathrm{N}+\sigma\,T}
      \\
      \zeta_\mathrm{N} (x) = 0
      \quad \mbox{ for all } x  \in \reali \setminus B_{\mathrm{N}+\sigma\,T+4}
      \\
      \norma{\zeta_\mathrm{N}}_{\C2 (\reali; \reali)} \leq \norma{\delta}_{\W21(\reali; \reali)}.
    \end{array}
    \right.
  \end{equation}
  Define
  \begin{displaymath}
    \mathtt{w}_\mathrm{N} (t,x,q)
    \coloneqq
    \zeta_\mathrm{N} (x) \;  \widehat{w} (t,x,q)
    \qquad \mbox{ for }
    (t, x, q) \in [0,T] \times \reali\times \reali\,,
  \end{displaymath}
  which by~\eqref{eq:23}, satisfies
  \begin{eqnarray}
    \label{eq:43}
    &
    & \sup\limits_{t \in [0,T], x \in \reali}
    \norma{\partial_x \mathtt{w}_\mathrm{N}(t,x,\cdot)}_{\W1\infty(\reali; \reali)}\,,
    \\
    \nonumber
    & \leq
    &
    \underbrace{
      \mathcal{C}
      \left(
      \norma{w}%
      _{\L\infty([0,T] \times \Su; \W1\infty(\reali; \reali))},\,
      \norma{\partial_x w}%
      _{\L\infty([0,T] \times \Su; \W1\infty(\reali; \reali))},\,
      \norma{\delta}_{\W21(\reali; \reali)}
      \right)
    }_{\mathrm{W}}
  \end{eqnarray}
  uniformly in $\mathrm{N}$.  Using~\ref{item:26}$\cdots$\ref{item:28} and
  the properties of $\zeta_{\mathrm{N}}$, by~\cite[Lemma~3.8 and
    Lemma~3.10]{opera-seconda}, the Cauchy problem
  \begin{equation}
    \label{eq:9}
    \left\{
    \begin{array}{l}
      \partial_t \mathtt{r}_\mathrm{N}
      +
      \partial_x \left(
      \mathtt{r}_\mathrm{N} \, \mathtt{w}_\mathrm{N}(t,x,\mathtt{r}_\mathrm{N})
      \right) = 0
      \\
      \mathtt{r}_\mathrm{N}(0,x) = \mathtt{r}_{o,\mathrm{N}}(x)
    \end{array}
    \right.
  \end{equation}
  admits a unique solution $\mathtt{r}_\mathrm{N}$ which satisfies
  $\mathtt{r}_\mathrm{N} (t,x) \in [0, R]$ for
  a.e.~$(t,x) \in [0,T]\times\reali$ by~\cite[Theorem~3]{Kruzkov1970}.
  Recall that by~\ref{item:18} in Theorem~\ref{thm:K}, for
  a.e.~$(t,x) \in [0,T]\times\Su$, it holds that $r(t,x) \in [0,R]$
  (which evidently implies that $\widehat{r}(t,x) \in [0,R]$ for
  a.e.~$(t,x) \in [0,T]\times \reali$).  We apply \Cref{prop:claudia}
  to the Cauchy problems~\eqref{eq:38} and~\eqref{eq:9}, so that
  $f(t,x,r) \coloneqq r \, \widehat{w}(t,x,r)$ and
  $f_*(t,x,r) \coloneqq r \, \mathtt{w}_\mathrm{N}(t,x,r)$. By
  construction, $f(t,x,r) = f_*(t,x,r)$ on
  $[0,T] \times B_{\mathrm{N} + \sigma T} \times \reali$.  Hence, for
  a.e.~$t \in [0,T]$,
  \begin{displaymath}
    \norma{\mathtt{r}_\mathrm{N}(t)- \widehat r(t)}_{\L1(B_{\mathrm{N} + \sigma (T-t)}; \reali )}
    \leq
    \norma{\mathtt{r}_{o,\mathrm{N}} - \widehat r_o}_{\L1(B_{\mathrm{N} + \sigma T}; \reali )} = 0
  \end{displaymath}
  by the choice~\eqref{eq:26} of the initial datum
  $\mathtt{r}_{o,\mathrm{N}}$.  In particular, for a.e.~$t \in [0,T]$,
  $\widehat{r}(t)$ coincides with $\mathtt{r}_\mathrm{N}(t)$ a.e.~on
  $B_{\mathrm{N}}$.
  \\
  \noindent
  \medskip \textbf{\underline{Step~2:} Total Variation Estimate.}
  \\
  Apply~\cite[Lemma~3.11]{opera-seconda} to~\eqref{eq:9}: for all
  $\mathrm{N} \in \naturali$ and $t \in [0,T]$,
  $\mathtt{r}_\mathrm{N}(t) \in \BV(\reali;[0,R])$ and
  \begin{eqnarray*}
    &
    & \tv\left(\mathtt{r}_\mathrm{N}(t)\right)
    \\
    & \leq
    & \tv(\mathtt{r}_{o,\mathrm{N}}) \, e^{\kappa_{o} t}
    +
    R \int_0^t e^{\kappa_{o} (t-\tau)}
    \int_\reali
    \norma{\partial_{xx}^2 \mathtt{w}_\mathrm{N}(\tau, x, \cdot)}%
    _{\L\infty(\reali; \reali)}
    \d{x} \d{\tau}
    \\
    & \leq
    &
    \tv(\mathtt{r}_{o,\mathrm{N}}) \, e^{\kappa_{o} t}
    +
    R \int_0^t e^{\kappa_{o} (t-\tau)}
    \int_{B_{\mathrm{N}+\sigma\,T}}
    \norma{\partial_{xx}^2 \widehat{w}(\tau, x, \cdot)}_{\L\infty(\reali; \reali)}
    \d{x}
    \d{\tau}
    \\
    &
    & +
    R \int_0^t e^{\kappa_{o} (t-\tau)}
    \int_{B_{\mathrm{N}+\sigma\,T+4}\setminus B_{\mathrm{N}+\sigma\,T}}
    \norma{\partial_{xx}^2 \mathtt{w}_\mathrm{N}(\tau, x, \cdot)}_{\L\infty(\reali; \reali)}
    \d{x}
    \d{\tau}
    \\
    & \leq
    &
    \tv(\mathtt{r}_{o,\mathrm{N}}) \, e^{\kappa_{o} t}
    + 2(\mathrm{N} + \sigma T ) \, R \,
    \sup\limits_{\tau \in [0,t]} \norma{\partial_{xx}^2 w(\tau)}_{\L1(\Su; \L\infty(\reali; \reali))}
    \frac{e^{\kappa_{o} t} -1 }{\kappa_{o}} + K
  \end{eqnarray*}
  where
  \begin{eqnarray*}
    \kappa_o
    & \coloneq
    & 3
    \left( 1+ R \right) \mathrm{W}
    \\
    &
    \geq
    & 3
    \left(
    \norma{\partial_x \mathtt{w}_\mathrm{N}}%
    _{\L\infty([0,T]\times \reali\times\reali;\reali)}
    +
    R \,
    \norma{\partial_{xr}^2 \mathtt{w}_\mathrm{N}}%
    _{\L\infty([0,T]\times\reali\times\reali;\reali)}
    \right)
    \,;
    \\
    K
    & \coloneq
    & \left(
    \norma{w}_{\L\infty([0,T]\times \Su \times \reali; \reali)}
    + 2\norma{\partial_{x} w}_{\L\infty([0,T]\times \Su \times \reali; \reali)}
    \right)
    8 \, R \,
    \frac{e^{\kappa_{o} t} -1 }{\kappa_{o}}
    \norma{\delta}_{\W21(\reali; \reali)}
    \\
    &
    & + \sup\limits_{\tau \in [0,t]}
    \norma{\partial_{xx}^2 w(\tau)}_{\L1(\Su; \L\infty(\reali; \reali))}
    8 \, R \,
    \frac{e^{\kappa_{o} t} -1 }{\kappa_{o}} \,.
  \end{eqnarray*}
  Exploiting the periodicity of $\widehat r(t)$ and $\widehat r_o$,
  together with the equalities $\widehat r(t) = \mathtt r_\mathrm{N} (t)$ on
  $B_\mathrm{N}$ and $\widehat r_o= \mathtt r_{o,\mathrm{N}}$ on
  $B_{\mathrm{N} + \sigma T}$, we get
  \begin{align*}
    \tv\left(\mathtt{r}_\mathrm{N} (t); B_\mathrm{N}\right)
     & =
    2 \, \mathrm{N} \, \tv\left( r (t); \Su\right) \,;
    \\
    \tv\left(
    \mathtt{r}_{o,\mathrm{N}}
    \right)
     & =
    \tv\left(\mathtt{r}_{o,\mathrm{N}}; B_{\mathrm{N} + \sigma T}\right)
    +
    \modulo{\mathtt{r}_{o,\mathrm{N}}\left((\mathrm{N}+\sigma T) -\right)}
    +
    \modulo{\mathtt{r}_{o,\mathrm{N}}\left((-\mathrm{N}-\sigma T)+\right)}
    \\
     & \leq
    2 \left(\mathrm{N} + \lfloor \sigma T \rfloor +1\right)
    \,
    \tv(r_o; \Su)
    +
    2 \, R \,.
  \end{align*}
  Hence,
  \begin{align*}
    \tv\left(r(t)\right)
     & \leq
    \frac{1}{2 \mathrm{N}}
    \Big[
      \left(
      2 (\mathrm{N} + \lfloor \sigma T \rfloor +1)
      \,
      \tv(r_o; \Su)
      +
      2 \, R
      \right)
      e^{\kappa_o t}
    \\
     & \qquad
      + 2\, \left(\mathrm{N} + \sigma T \right) \, R \;
      \sup\limits_{\tau \in [0,t]} \norma{\partial_{xx}^2 w(\tau)}_{\L1(\Su; \L\infty(\reali; \reali))}\,
      \frac{e^{\kappa_{o} t} -1 }{\kappa_{o}}
      +
      K
      \Big]
    \\
     &
    \underset{\mathrm{N}\to + \infty}{\longrightarrow}
    \tv(r_o)
    e^{\kappa_o t}
    +
    R \;
    \sup\limits_{\tau \in [0,t]} \norma{\partial_{xx}^2 w(\tau)}_{\L1(\Su; \L\infty(\reali; \reali))}
    \frac{e^{\kappa_{o} t} -1 }{\kappa_{o}},
  \end{align*}
  completing the proof of~\ref{item:20}.

  \paragraph{Proof of~\ref{item:24}.}
  {

    Assume without loss of generality
    that
    \begin{align*}
           &
      \norma{w_*}_{\L\infty([0,T] \times \Su \times B_R; \reali)}
      +
      R \,
      \norma{\partial_r w_*}_{\L\infty([0,T] \times \Su \times B_R; \reali)}
      \\
      \leq &
      \norma{w}_{\L\infty([0,T] \times \Su \times B_R; \reali)}
      +
      R \,
      \norma{\partial_r w}_{\L\infty([0,T] \times \Su \times B_R; \reali)} .
    \end{align*}

  }
  For $(t, x, q) \in [0,T] \times \reali \times \reali$, define
  $\mathtt{w}_* (t,x,q) \coloneqq \zeta_{1} (x) \;
    \widehat{w}_*(t,x,q)$ and
  $\mathtt{w} (t,x,q) \coloneqq \zeta_1 (x) \; \widehat{w} (t,x,q)$
  where $\zeta_1$ is as in~\eqref{eq:5} with
  $\mathrm{N}=1$. Similarly, recall $\mathtt{r}_{o,1}$
  from~\eqref{eq:26}.

  Consider the Cauchy problem
  \begin{equation}
    \left\{
    \begin{array}{l}
      \partial_t \mathtt{r}_*
      + \partial_x \left(
      \mathtt{r}_*
      \, \mathtt{w}_*(t,x,\mathtt{r}_*)
      \right) = 0
      \\
      \mathtt{r}_*(0,x) = \mathtt{r}_{o,1}(x) \,.
    \end{array}
    \right.
  \end{equation}
  By~\eqref{eq:43}, both $\partial_x \mathtt{w}$ and
  $\partial_x \mathtt{w}_*$ are in
  $\L\infty([0,T]\times \reali \times \reali; \reali)$ and are
  compactly supported in space by construction. Hence,
  \begin{displaymath}
    \sup_{ t \in [0,T]} \int_\reali \sup\limits_{q \in \reali }
    \modulo{\partial_x \mathtt{w}(t, x, q)} \d x
    +
    \sup_{ t \in [0,T]} \int_\reali \sup\limits_{q \in \reali }
    \modulo{\partial_x \mathtt{w}_*(t, x, q)} \d x
    < + \infty \,.
  \end{displaymath}
  Call $\mathtt{r}$ the solution to~\eqref{eq:9} with
  $\mathrm{N} = 1$. By~\cite[eqn.~(3.49) in
    Lemma~3.11]{opera-seconda}, for a.e.~$t \in [0,T]$,
  \begin{align*}
     & \norma{\mathtt{r}(t) - \mathtt{r}_*(t)}_{\L1(\reali; \reali)}
    \\
     & \leq
    t \, e^{\kappa_o t} \max \left\{
    R,\,
    \tv(\mathtt{r}_{o,1})
    +
    R \int_0^t
    \int_\reali
    \norma{\partial^2_{xx} \mathtt{w}(\tau, x, \cdot)}_{\L\infty(\reali; \reali)}
    \d x \d \tau
    \right\}
    \\
     & \quad \times
    \Bigl(
    \norma{\mathtt{w} - \mathtt{w}_*}_{\L\infty([0,T]\times \reali \times \reali; \reali )}
    +
    R \,
    \norma{\partial_r\mathtt{w} - \partial_r\mathtt{w}_*}_{\L\infty([0,T]\times \reali \times \reali; \reali )}
    \\
     & \qquad
    +
    \sup_{t \in [0,T]} \int_\reali \sup\limits_{q \in \reali }
    \modulo{\partial_x \mathtt{w}(t, x, q) - \partial_x \mathtt{w}_*(t, x, q)} \d x
    \Bigr)
    \\
     & \leq
    t \, e^{\kappa_o t} \max \left\{
    R,\,
    \tv(\mathtt{r}_{o,1})
    +
    2 \,R \,
    (5 + \sigma \, T)
    \mathcal{C}
    \left(
    \begin{array}{@{}c@{}}
        \norma{w}_{\L\infty([0,T]\times \Su \times \reali; \reali)},
        \\
        \norma{\partial_x w}_{\L\infty([0,T]\times \Su \times \reali; \reali)},
        \\
        \sup\limits_{\tau \in [0,t]} \norma{\partial_{xx}^2 w(\tau)}_{\L1(\Su; \L\infty(\reali; \reali ))}
      \end{array}
    \right)
    t
    \right\}
    \\
     & \quad \times
    \Bigl(
    \norma{w - w_*}_{\L\infty([0,T]\times \Su \times \reali; \reali )}
    +
    R \,
    \norma{\partial_r w - \partial_r w_*}_{\L\infty([0,T]\times \Su \times \reali; \reali )}
    \\
     &
    \qquad
    +
    \sup_{t \in [0,T]} \int_{B_{5 +\sigma T}} \sup_{q \in \reali}
    \modulo{
      \zeta_1'(x)
      \left(
      \widehat{w}(t, x, q) - \widehat{w}_*(t, x, q)
      \right)} \d x
    \\
     & \qquad
    +
    \sup_{t \in [0,T]} \int_{B_{5+\sigma T}} \sup_{q \in \reali}
    \modulo{
      \zeta_1(x)
      \left(
      \partial_x\widehat{w}(t, x, q) - \partial_x\widehat{w}_*(t, x, q)
      \right)} \d x
    \Bigr)
    \\
     & \leq
    t \, e^{\kappa_o t} \max \left\{
    R,\,
    \tv(\mathtt{r}_{o,1})
    +
    2\,
    R \,
    (5 + \sigma \, T)
    \mathcal{C}
    \left(
    \begin{array}{@{}c@{}}
        \norma{w}_{\L\infty([0,T]\times \Su \times \reali; \reali)},
        \\
        \norma{\partial_x w}_{\L\infty([0,T]\times \Su \times \reali; \reali)},
        \\
        \sup\limits_{\tau \in [0,t]} \norma{\partial_{xx}^2 w(\tau)}_{\L1(\Su; \L\infty(\reali; \reali ))}
      \end{array}
    \right)
    t
    \right\}
    \\
     & \quad \times
    \Bigl(
    \norma{w - w_*}_{\L\infty([0,T]\times \Su \times \reali; \reali )}
    +
    R
    \norma{\partial_r w - \partial_r w_*}_{\L\infty([0,T]\times \Su \times \reali; \reali )}
    \\
     &
    \qquad
    +
    \norma{\zeta_1'}_{\L\infty(\reali; \reali)} \,
    2\,
    (5 + \sigma \, T)
    \norma{w - w_*}_{\L\infty([0,T]\times \Su \times \reali; \reali )}
    \\
     & \qquad
    +
    \norma{\zeta_1}_{\L\infty(\reali; \reali)} \,
    2\,
    (5 + \sigma \, T)
    \sup_{t \in [0,T]} \int_{\Su} \sup_{q \in \reali}
    \modulo{
      \partial_x w (t, x, q) - \partial_x w_*(t, x, q)
    } \d x
    \Bigr)
    \\
     & \leq
    t \, e^{\kappa_o t}
    \mathcal{C}
    \left(
    \begin{array}{c}
        \norma{w}_{\L\infty([0,T]; \W1\infty(\Su \times \reali; \reali))},
        \sup\limits_{\tau \in [0,t]} \norma{\partial_{xx}^2 w(\tau)}_{\L1(\Su; \L\infty(\reali; \reali ))},
        \\
        { \norma{w_*}_{\L\infty([0,T] \times \Su ; \W1\infty( \reali; \reali))},
        } R, \,
        \tv(r_o), T
      \end{array}
    \right)
    \\
     & \quad
    \times
    \left(
    \norma{w - w_*}_{\L\infty([0,T]\times \Su; \W1\infty(\reali; \reali))}
    +
    \sup_{t \in [0,T]} \int_{\Su} \sup_{q \in \reali}
    \modulo{
      \partial_x w (t, x, q) - \partial_x w_*(t, x, q)
    } \d x
    \right).
  \end{align*}
  For a.e.~$(t,x) \in [0,T]\times[0,1]$,
  $\widehat{r} (t,x) = \mathtt{r} (t,x)$ and
  $\widehat{r}_* (t,x) = \mathtt{r}_* (t,x)$, proving~\ref{item:24}.

  \paragraph{Proof of~\ref{item:25}.}
  The proof relies on~\cite[Theorem~4.3.1]{dafermos} and on the $\tv$
  estimate~\eqref{eq:50} in~\ref{item:20} of Theorem~\ref{thm:K}. See
  also the proof of~\cite[Lemma~3.12]{opera-seconda}.
\end{proofof}

\subsection{Proofs Related to the Nonlocal Periodic Case}
\label{sec:nonl-peri-case-2}

Here we collect some results on the scalar linear periodic continuity
equation, which are not easily available in the literature.

Let $v \in \L\infty([0,T]\times \Su; \reali)$ and, for
$t_o \in [0,T]$, consider the Cauchy problem
\begin{equation}
  \label{eq:7}
  \left\{
  \begin{array}{ll}
    \partial_t \rho +
    \partial_x\left(\rho \, v(t,x)\right) = 0
     & \mbox{in } [0,T]\times \Su
    \\
    \rho(t_o) =  \rho_o
     & \mbox{in }\Su \,.
  \end{array}
  \right.
\end{equation}
Although the scalar linear equation~\eqref{eq:7} is a
particular case of~\eqref{eq:6}, the assumptions on the speed function
$v$, considered in this Section, are less restrictive than the ones
in Section~\ref{sec:kruv-zkov-periodic-1},
namely~\ref{item:26},~\ref{item:28} and~\ref{item:31}. Thus, \Cref{thm:K} can not be
used here.

Recall the following standard definition.

\begin{definition}
  \label{def:solLin}
  By \emph{weak solution to~\eqref{eq:7}} on $[0,T]$ we mean
  $\rho \in \C0\left([0,T]; \L1(\Su; \reali)\right)$ such that
  $\rho (t_o) = \rho_o$ and, for all
  $\phi \in \Cc1(]0,T[ \times \reali; \reali)$,
  \begin{displaymath}
    \int_0^T \int_{\reali} \left(\widehat \rho(t,x) \,\partial_t
    \phi(t,x) + \widehat\rho(t,x) \, \widehat{v}(t,x) \, \partial_x
    \phi(t,x) \right) \d{x} \d{t} = 0 \,.
  \end{displaymath}
\end{definition}
Under suitable assumptions, $\rho$ admits the expression
\begin{eqnarray}
  \label{eq:13}
  \begin{array}{r@{\,}c@{\,}l@{\qquad}r@{\,}c@{\,}l}
    \rho(t,x)
     & =
     & \rho_o\left(X(t_o;t,x)\right) \; \mathcal{E}(t_o,t,x)
     & (t,x)
     & \in
     & [0,T] \times \Su
    \\
    \mathcal{E}(t_1,t_2,x_o)
     & \coloneqq
     & \exp \left(
    \int_{t_2}^{t_1} \partial_x v \left(s, X(s;t_2,x_o)\right) \d{s}
    \right)
     & (t_1,t_2, x_o)
     & \in
     & [0,T]^2 \times \Su
  \end{array}
\end{eqnarray}
(see~\Cref{lem:1}) and, for $(t_o,x_o) \in [0,T]\times \Su$, the
characteristic curve $t \mapsto X(t;t_o,x_o)$ is the absolutely
continuous function solving the Cauchy problem
\begin{equation}
  \label{eq:11}
  \begin{cases}
    \frac{\d{~}}{\d{t}} X(t;t_o,x_o) = v\left(t, X(t;t_o,x_o)\right)
    \\
    X(t_o;t_o,x_o) = x_o \,.
  \end{cases}
  \quad (t_o,x_o) \in [0,T] \times \Su\,.
\end{equation}
{

Define also $\widetilde X (t; t_o, x_o)$
as the solution
to
\begin{equation}
  \label{eq:11bis}
  \begin{cases}
    \frac{\d{~}}{\d{t}} \widetilde X(t;t_o,x_o) =\widehat v\left(t, \widetilde X(t;t_o,x_o)\right)
    \\
    \widetilde X(t_o;t_o,x_o) = x_o \,.
  \end{cases}
  \quad (t_o,x_o) \in [0,T] \times \reali\,.
\end{equation}
Note that $\widetilde X(t; t_o, \cdot) \colon \reali \to \reali$ is a lift of
the map $X(t; t_o, \cdot) \colon \Su \to \Su$ to the whole real line,
according to~\eqref{eq:16}.

}

\begin{lemma}[Properties of the characteristics]
  \label{lem:characteristic}
  Let $v\in \C0([0,T]\times \Su; \reali)$ and assume there exists
  $L > 0$ such that
  \begin{equation}
    \label{eq:42}
    \modulo{v(t,x) - v(t,y)} \leq L \; \modulo{x-y}
    \quad \text{for all} \quad (t,x,y) \in [0,T]\times \Su \times \Su \,.
  \end{equation}
  Then,
  \begin{enumerate}[label=\bf(\roman*)]
    \item\label{item:13} For all $(t_o, x_o) \in [0,T] \times \Su$,
          problem~\eqref{eq:11} admits a unique solution
          $X(\cdot; t_o, x_o)\colon\! [0,T] \to \Su$.  Moreover,
          $X(t; \cdot, \cdot) \in \C0([0,T]\times \Su; \Su)$, for all
          $t \in [0,T]$.

    \item\label{item:15} If for all $t \in [0,T]$,
          $\partial_x v (t, \cdot) \in \C0 (\Su; \reali)$, then, for every
          $t_o,t \in [0,T]$, the map $X(t; t_o, \cdot) \colon \Su \to \Su$
          is a diffeomorphism. Moreover, for all $t \in [0,T]$,
          $X (t, \cdot, \cdot) \in \C1([0,T]\times \Su; \Su)$ and, using
          notation~\eqref{eq:13}, for all
          $(t, t_o, x_o) \in [0,T]^2 \times \Su$,
          \begin{equation}
            \label{eq:12}
            \partial_{x_o} X(t;t_o,x_o) =
            \mathcal{E}(t,t_o,x_o)
            \,,\qquad
            \partial_{t_o} X(t; t_o, x_o) =
            - \mathcal{E}(t, t_o, x_o) \, v(t_o,x_o)\,.
          \end{equation}
    \item \label{item:33} If
          $v \in \L\infty\left([0,T]; \W2\infty(\Su; \reali)\right)$,
          $\partial_x v \in \C0([0,T]\times \Su; \reali)$ and
          $\partial^2_{xx} v (t, \cdot) \in \C0 (\Su; \reali)$ for all
          $t \in [0,T]$, then
          $\mathcal{E}(t, \cdot, \cdot) \in \C1 ([0,T]\times\Su; \reali)$
          for all $t \in [0,T]$ and for all
          $(t, t_o, x_o) \in [0,T]^2 \times \Su$,
          \begin{align*}
            \partial_{x_o}\mathcal{E}(t,t_o,x_o)
             & =
            \mathcal{E}(t, t_o, x_o)
            \int_{t_o}^t \partial_{xx}^2 v\left(s, X(s;t_o,x_o)\right)\,
            \mathcal{E}(s, t_o, x_o)
            \dd{s} \,;
            \\
            \nonumber
            \partial_{t_o} \mathcal{E}(t, t_o, x_o)
             & = -
            \mathcal{E}(t, t_o, x_o)
            \!\!
            \left( \!\!
            \!\partial_x v(t_o, x_o)
            \!+ \!\! \int_{t_o}^t \!\!
            \partial_{xx}^2 \! v\left(s, X(s;t_o,x_o)\right)
            \mathcal{E}(s, t_o, x_o)
            \, v(t_o, x_o) \dd{s}
            \!
            \right)\!\!.
          \end{align*}

    \item\label{item:14} {  If
            $v_* \in \C0([0,T]\times \Su; \reali)$ also
            satisfies~\eqref{eq:42} then, with obvious notation, for
            all $(t, t_o, x_o) \in [0,T]^2 \times \reali$
            \begin{displaymath}
              \sup_{x_o \in \reali} \modulo{
                \widetilde{X}(t; t_o, x_o) -
                \widetilde{X}_*(t; t_o,x_o)}
              \leq
              \modulo{\int_{t_o}^t
                \norma{v(\tau) - v_*(\tau)}_{\L\infty(\Su; \reali)} \dd\tau}
              \;
              e^{L \, \modulo{t-t_o}}\,,
            \end{displaymath}
            where
            $\widetilde{X}(t; t_o, x_o)$ and, similarly,
            $\widetilde{X}_*$, solves~\eqref{eq:11bis}.
          }
  \end{enumerate}
\end{lemma}

\begin{proofof}{Lemma~\ref{lem:characteristic}}
  Consider the different items separately.

  \paragraph{Proof of~\ref{item:13}.}  The well-posedness
  of~\eqref{eq:11} directly follows from the usual Cauchy-Lipschitz
  Theorem (see, e.g., \cite[Theorem~2.1.1 and
    Theorem~2.1.3]{MR2347697}) applied to the { Cauchy
      problem~\eqref{eq:11bis}.  }

  \paragraph{Proof of~\ref{item:15}--\ref{item:33}.} By~\ref{item:13},
  $X\left(t_o; t, X(t; t_o,x_o)\right) = x_o$ for all
  $t, t_o \in [0,T]$, $x_o \in \Su$, proving that $X(t; t_o, \cdot)$
  is invertible in $\Su$.  The continuous differentiability of the
  maps $X(t,t_o,\cdot)$ and $X(t; \cdot, x_o)$ and~\eqref{eq:12}
  descend from~\cite[Theorem~2.3.2, Theorem~2.3.3]{MR2347697}.  Then,
  by~\eqref{eq:42} and $\partial_x v (t, \cdot) \in \C0(\Su; \reali)$,
  one gets $\mathcal{E}(t,\cdot, \cdot)$ is continuous, hence
  $\partial_x X(t; \cdot, \cdot)$ and $\partial_t X(t; \cdot, \cdot)$
  are, proving $X(t; \cdot, \cdot) \in \C1([0,T] \times \Su; \Su)$.
  Finally,~\ref{item:15} implies~\ref{item:33}.

  \paragraph{Proof of~\ref{item:14}.}  {  Since for all
    $t,t_o \in [0,T]$, $x_o \in \reali$, $\widetilde X(t; t_o, x_o)$
    is defined in~\eqref{eq:11bis}, the stability with respect to the
    velocity field is a direct application of Gronwall Lemma.  }
\end{proofof}

\begin{lemma}[Well-Posedness of~\eqref{eq:7} in $\Su$]
  \label{lem:1}
  Let $t_o \in [0,T]$, $v \in \C0([0,T]\times \Su; \reali)$
  satisfy~\eqref{eq:42} and fix $\rho_o \in \L\infty(\Su;
    \reali)$. Then
  \begin{enumerate}[label=\bf(\roman*)]
    \item \label{item:1} There exists a unique weak solution
          to~\eqref{eq:7} in the sense of Definition~\ref{def:solLin}.
    \item \label{item:2} The map $\rho$ defined
          in~\eqref{eq:13}--\eqref{eq:11} is the weak solution
          to~\eqref{eq:7}.
    \item \label{item:9} For all $t \in [0,T]$,
          $\norma{\rho(t)}_{\L1(\Su; \reali)} = \norma{\rho_o}_{\L1(\Su;
              \reali)}$.
    \item \label{item:6} Let $\rho_{o,*} \in \L\infty(\Su; \reali)$ and
          $\rho_*$ be the corresponding solution to~\eqref{eq:7}.  Then, for
          all $t \in [0,T]$
          \begin{equation}
            \norma{\rho(t) - \rho_*(t)}_{\L1(\Su; \reali)}
            =
            \norma{\rho_o - \rho_{o,*}}_{\L1(\Su; \reali)}.
          \end{equation}
    \item\label{item:7} For all $t \in [0,T]$,
          $\norma{\rho(t)}_{\L\infty(\Su; \reali)} \leq
            \norma{\rho_o}_{\L\infty(\Su; \reali)} \, e^{L \,
                \modulo{t-t_o}}$.
    \item \label{item:5} If
          $v \in \L\infty\left([0,T]; \W2\infty(\Su; \reali)\right)$ and
          $\rho_o \in \cBV{} (\Su; \reali)$, then for all $t \in [0,T]$
          \begin{displaymath}
            \tv\left(\rho(t)\right)
            \leq
            \left(
            \tv (\rho_o)
            +
            \norma{\partial_{xx}^2 v}_{\L\infty ([0,T]\times \Su; \reali)}
            \norma{\rho_o}_{\L1 (\Su; \reali)}
            \, \modulo{t-t_o}
            \right) e^{L\,\modulo{t-t_o}}.
          \end{displaymath}
          Moreover, if $\rho_o \in \W11(\Su; \reali)$ then
          $\rho(t) \in \W11(\Su; \reali)$ for all $t \in [0,T]$.
    \item\label{item:8} If
          $v \in \L\infty\left([0,T]; \W2\infty(\Su; \reali)\right)$ and
          $\rho_o \in \cBV{} (\Su; \reali)$, then the solution
          to~\eqref{eq:7} is $\L1$-locally Lipschitz continuous in time: for
          all $t, t_* \in [0,T]$,
          \begin{displaymath}
            \norma{\rho(t) - \rho(t_*)}_{\L1(\Su; \reali)}
            \leq
            \mathcal{C}\left(
            \begin{array}{c}
                \norma{v}_{\L\infty([0,T]; \W2\infty(\Su; \reali))},
                \norma{\rho_o}_{\cBV{}(\Su; \reali)}, T
              \end{array}
            \right) \modulo{t-t_*} \,.
          \end{displaymath}

    \item \label{item:37a} Fix $\ell =1,2$. If
          $\rho_o \in \C\ell(\Su; \reali)$,
          $v \in \L\infty\left([0,T]; \W{\ell+1}\infty(\Su; \reali)\right)$,
          $\partial_x v \in \C0([0,T]\times \Su; \reali)$ and for all
          $t \in [0,T]$,
          $\partial^2_{xx} v (t, \cdot) \in \C{\ell-1}(\Su; \reali)$, then
          $\rho \in \C1([0,T]\times \Su; \reali)$ and, if $\ell = 2$,
          $\partial_{xx}^2 \rho \in \C0([0,T]\times \Su; \reali)$ and
          $\partial_{tx}^2 \rho \in \L\infty([0,T]\times \Su; \reali)$.
          Moreover,
          \begin{align}
            \label{eq:62}
            \!\!\!\!
            \norma{\rho}_{\L\infty([0,T]; \W\ell\infty(\Su; \reali))}
            \leq
            \mathcal{C}
            \left(
            \norma{v}_{\L\infty([0,T]; \W{\ell+1}\infty(\Su; \reali))},
            \norma{\rho_o}_{\W\ell\infty(\Su; \reali)}, T
            \right)e^{(\ell+1)\, L \, T};\!\!
            \\
            \label{eq:66}
            \!\!\!\!
            \norma{\rho}_{\L\infty([0,T]; \W\ell{1}(\Su; \reali))}
            \leq
            \mathcal{C}
            \left(
            \norma{v}_{\L\infty([0,T]; \W{\ell+1}\infty(\Su; \reali))},
            \norma{\rho_o}_{\W\ell{1}(\Su; \reali)}, T
            \right)e^{(\ell+1)\, L \, T} \, .
          \end{align}
    \item \label{item:3a} Fix $\ell = 0,1$.  If
          $\rho_o \in \cBV{\ell} (\Su; \reali)$,
          $v \in \L\infty\left([0,T]; \W{\ell+2}\infty(\Su; \reali)\right)$
          and
          $v_* \in \C0([0,T]\times \Su; \reali) \cap \L\infty\left([0,T];
            \W{\ell+1}\infty(\Su; \reali)\right)$ satisfies~\eqref{eq:42},
          calling $\rho, \rho_*$ the solutions associated to $v$ and $v_*$,
          then, for all $t \in [0,T]$
          \begin{align*}
             & \quad
            \norma{\rho(t) -  \rho_*(t)}_{\W{\ell}1(\Su; \reali)}
            \\
             & \leq
            \mathcal{C}
            \left(
            \begin{array}{@{}c@{}}
                \norma{\partial_{x} v}_{\L\infty([0,T]; \W{\ell+1}\infty(\Su; \reali))}, \\
                \norma{\partial_{x} v_*}_{\L\infty([0,T]; \W{\ell}\infty(\Su; \reali))},
                \\
                \norma{\rho_o}_{\cBV{\ell}(\Su; \reali)}, T
              \end{array}
            \right)
            e^{(\ell +3 ) L\, \modulo{t-t_o}} \,
            \modulo{
            \int_{t_o}^t
            \norma{v (\tau) - v_* (\tau)}_{\W{\ell +1 }\infty(\Su; \reali)}
            \d\tau
            } \,.
          \end{align*}
  \end{enumerate}
\end{lemma}

\begin{proofof}{Lemma~\ref{lem:1}}
  Group the different items as follows.

  \paragraph{Proof of~\ref{item:1} and~\ref{item:2}.}  This is a
  direct consequence of the well-posedness of linear continuity
  equations in the non-periodic setting, see for
  instance~\cite[Lemma~2.11 and Lemma~2.12]{keimer_pflug}. We
  underline that $\rho \in \C0\left([0,T]; \L1(\Su; \reali)\right)$
  due to the translation invariance of the characteristics $X$,
  see~\eqref{eq:11}, stated in~\ref{item:13} of
  \Cref{lem:characteristic}, and of the periodicity of $\rho_o$ and
  $v$.

  \paragraph{Proof of~\ref{item:9} and~\ref{item:6}.}  Simply use the
  change of variable $\xi = X(t_o; t, x)$, see~\ref{item:15} in
  \Cref{lem:characteristic}. By linearity, also~\ref{item:6} follows.

  \paragraph{Proof of~\ref{item:7}} It directly follows
  from~\eqref{eq:13}.

  \paragraph{Proof of~\ref{item:5}.}  Assume for notational simplicity
  that $t > t_o$, the other case being analogous. Let
  $\mathrm{N} \in \naturali \setminus\{0\}$. Call
  $\mathrm{V} \coloneqq \norma{v}_{\L\infty([0,T]\times \Su; \reali)}$
  and $a = -\mathrm{N} - \lfloor\mathrm{V}\,(t-t_o)\rfloor-1$ and
  $b = \mathrm{N} + \lfloor\mathrm{V}\,(t-t_o)\rfloor+1$. Call
  $\mathtt{p}$ the solution to
  \begin{equation}
    \label{eq:54}
    \left\{
    \begin{array}{l@{\qquad}l}
      \partial_\tau \mathtt{p}
      + \partial_x \left( \mathtt{p} \, \widehat{v} (\tau,x) \right)
      = 0
       & \mbox{in } [t_o,t] \times \reali
      \\
      \mathtt{p} (t_o,x) =
      \left\{
      \begin{array}{ll}
        \widehat{ \rho}_o(x)
         & \mbox{in } [a,b]
        \\
        0
         & \mbox{in } \reali\setminus [a,b]
      \end{array}
      \right.
       & \mbox{in }\reali \, .
    \end{array}
    \right.
  \end{equation}
  Standard results on linear continuity equations, see for
  instance~\cite[eqn.~(2.5)]{stability_loc}
  or~\cite[eqn.~(4.21)]{opera-prima}, imply that
  $\mathtt{p}(t) \in \BV(\reali; \reali)$ and that
  \begin{align}
    \nonumber
     & \quad \tv\left(\mathtt{p}(t) \right)
    \\
    \nonumber
     & \leq
    \left(
    \tv\left(\mathtt{p}(t_o)\right)
    +
    \norma{\partial_{xx}^2 v}_{\L\infty([t_o,t]\times \Su; \reali)}
    \norma{\mathtt{p}(t_o)}_{\L1(\reali; \reali)}
      (t-t_o)
    \right)e^{L (t-t_o)}
    \\
    \nonumber
     & =
    \left(
    \tv(\widehat{\rho}_o; \mathopen]a,b\mathclose[)
    +
    \modulo{\widehat{\rho}_o (a+)}
    +
    \modulo{\widehat{\rho}_o (b-)}
    +
    \norma{\partial_{xx}^2 v}_{\L\infty([t_o,t]\times \Su; \reali)}
    \norma{\widehat{\rho}_o}_{\L1([a,b]; \reali)}
      (t-t_o)
    \right)e^{L (t-t_o)}
    \\
    \label{eq:56}
     & =
    \left(
    2\mathrm{N} \tv(\widehat{\rho}_o, \mathopen]0,1])
    +
    2\mathrm{N}
    \norma{\partial_{xx}^2 v}_{\L\infty([t_o,t]\times \Su; \reali)}
    \norma{\rho_o}_{\L1([0,1]; \reali)}
    (t-t_o)
    + K
    \right)
    e^{L (t-t_o)} \,,
  \end{align}
  where $K$ is the uniformly (in $N$) bounded quantity
  \begin{eqnarray*}
    K
    & \coloneqq
    & \modulo{\widehat{\rho}_o (a+)}
    +
    \modulo{\widehat{\rho}_o (b-)}
    +
    \tv(\widehat{\rho}_o; \mathopen]a, -\mathrm{N}\mathclose[  \cup \mathopen]\mathrm{N}, b\mathclose[)
    \\
    &
    & +
    \norma{\partial_{xx}^2 v}_{\L\infty([t_o,t]\times \Su; \reali)}
    \norma{\widehat{\rho}_o}_{\L1(\mathopen]a, -\mathrm{N}\mathclose[  \cup \mathopen]\mathrm{N}, b\mathclose[; \reali)}
    (t-t_o)
    \\
    & =
    & \modulo{\widehat{\rho}_o (a+)}
    +
    \modulo{\widehat{\rho}_o (b-)}
    +
    \tv\left(
    \widehat{\rho}_o;
    \mathopen]-\lfloor \mathrm{V} (t-t_o) \rfloor-1, 0\mathclose[
        \cup
        \mathopen]0, \lfloor \mathrm{V} (t-t_o)\rfloor +1 \mathclose[
    \right)
    \\
    &
    & +
    \norma{\partial_{xx}^2 v}_{\L\infty([t_o,t]\times \Su; \reali)}
    \norma{\widehat{\rho}_o}_{\L1(
    \mathopen]-\lfloor \mathrm{V} (t-t_o) \rfloor-1, 0\mathclose[
        \cup
        \mathopen]0, \lfloor \mathrm{V} (t-t_o)\rfloor +1 \mathclose[
    ; \reali)} \,
    (t-t_o) \,.
  \end{eqnarray*}
  Hence,
  \begin{flalign*}
    \tv\left( \rho (t) \right)
     & =
    \dfrac{1}{2\mathrm{N}} \,
    \tv \left( \widehat\rho (t); \mathopen]-\mathrm{N}, \mathrm{N}] \right)
    \tag*{[\mbox{By~\eqref{eq:44} and periodicity}]}
    \\
     & =
    \dfrac{1}{2\mathrm{N}} \,
    \tv \left( \mathtt{p} (t); \mathopen]-\mathrm{N}, \mathrm{N}]\right)
    \tag*{[\mbox{By~\cite[Theorem~1]{Kruzkov1970}}]}
    \\
     & \leq
    \left(
    \tv (\rho_o)
    +
    \norma{\partial_{xx}^2 v}_{\L\infty ([t_o,t]\times \Su; \reali)}
    \norma{\rho_o}_{\L1 (\Su; \reali)}
    \, (t-t_o)
    +
    \frac{K}{2\mathrm{N}}
    \right) e^{L (t-t_o)}
     & [\mbox{By~\eqref{eq:56}}]
    \\
     & \underset{\mathrm{N} \to\infty}{\longrightarrow}
    \left(
    \tv (\rho_o)
    +
    \norma{\partial_{xx}^2 v}_{\L\infty ([t_o,t]\times \Su; \reali)}
    \norma{\rho_o}_{\L1 (\Su; \reali)}
    \, (t-t_o)
    \right) e^{L (t-t_o)}
  \end{flalign*}
  completing the proof of~\ref{item:5}.

  \paragraph{Proof of~\ref{item:8}.} Assume that $t > t_*$, the other
  cases being analogous. We now follow the techniques
  in~\cite[Theorem~4.3.1]{dafermos}, for more details
  see~\cite[Lemma~3.12]{opera-seconda}.  Fix
  $\epsilon \in \mathopen]0,1\mathclose[$. We have
  \begin{align*}
     & \norma{\rho (t) - \rho (t_*)}_{\L1 (\Su; \reali)}
    \\
     & =
    \norma{\widehat\rho (t) - \widehat\rho (t_*)}_{\L1 ([0,1]; \reali)}
    \\
     & =
    \lim_{\epsilon \to 0}
    \norma{\widehat\rho (t) - \widehat\rho (t_*)}_{\L1 ([0,1+\epsilon]; \reali)}
    \\
     & =
    \lim_{\epsilon \to 0}\;
    \sup
    \left\{
    \int_0^{1+\epsilon}
    \left( \widehat\rho (t,x) - \widehat\rho (t_*,x) \right) \, \psi (x)\d{x}
    \colon
    \psi \in \Cc1([0, 1+\epsilon];\reali),
    \norma{\psi}_{\L\infty(\mathopen]0, 1+\epsilon \mathclose[;
    \reali)} \leq 1
    \right\}
    \\
     & \leq
    \lim_{\epsilon \to 0} \int_{t_*}^t
    \tv\left(  \widehat\rho(\tau,x) \, \widehat{v}(\tau,x);
    \mathopen]0, 1+\epsilon \mathclose[\right)
    \dd \tau
    \tag*{[\mbox{As in~\cite[Theorem~4.3.1]{dafermos}}]}
    \\
     & \leq
    \lim_{\epsilon \to 0} \int_{t_*}^t
    \left(
    \tv\left( \widehat\rho(\tau);
    \mathopen]0, 1+\epsilon \mathclose[ \right)
    \norma{v}_{\L\infty([t_*,t]\times \Su; \reali)}
    +
    \, \norma{\rho(\tau)}_{\L\infty(\Su; \reali) }
    \tv\left(\widehat{v}(\tau); \mathopen]0, 1+\epsilon \mathclose[\right)
    \right)
    \dd \tau
    \\
     & =
    \int_{t_*}^t
    \left(
    \tv\left( \widehat\rho(\tau);
    \mathopen]0, 1 ] \right)
    \norma{v}_{\L\infty([t_*,t]\times \Su; \reali)}
    +
    \, \norma{\rho(\tau)}_{\L\infty(\Su; \reali) }
    \tv\left(\widehat{v}(\tau); \mathopen]0, 1]\right)
    \right)
    \dd \tau
    \\
     & \leq
    \left(
    \tv (\rho_o)
    +
    \norma{\partial_{xx}^2 v}_{\L\infty ([0,T]\times \Su; \reali)}
    \norma{\rho_o}_{\L1 (\Su; \reali)}
    \, T
    \right) e^{L\,T}
    \norma{v}_{\L\infty([0,T]\times \Su; \reali)}
    (t-t_*)
    \tag*{[\mbox{By~\ref{item:5}}]}
    \\
     & +
    L\, \norma{\rho_o}_{\L\infty(\Su; \reali)} \, e^{L T} (t-t_*) \,.
    \tag*{[\mbox{By~\ref{item:7}}]}
  \end{align*}
  completing the proof of~\ref{item:8}.

  \paragraph{Proof of~\ref{item:37a}.}
  The fact that $\rho \in \C1([0,T]\times \Su; \reali)$ follows from
  expression~\eqref{eq:13}--\eqref{eq:11}
  and~\ref{item:15}--\ref{item:33} in \Cref{lem:characteristic}.  The
  formul\ae~for the computation of
  $\partial_x \rho, \partial_{xx}^2\rho$ and $\partial_{tx}^2\rho$
  in~\cite[Proof of~\textbf{(N1)}-\textbf{(N5)-\textbf{(N6)}} in
    Lemma~3.5]{opera-seconda}, show that
  $\partial_{xx}^2 \rho\in \C0([0,T]\times \Su; \reali)$,
  $\partial_{tx}^2\rho\in \L\infty([0,T]\times \Su; \reali)$
  and~\eqref{eq:62}--\eqref{eq:66}.

  \paragraph{Proof of~\ref{item:3a}.}
  Assume $t > t_o$. The proof is an adaptation of standard results
  about the stability on velocity fields,
  see~\cite[Lemma~4.6]{opera-prima} and~\cite[Proof of
    Lemma~3.7]{opera-seconda}:
  \begin{align*}
     & \norma{\rho(t) - \rho_*(t)}_{\L1(\Su; \reali)}
    \\
     & =
    \int_{\Su} \modulo{
      \rho_o\left(X(t_o;t,x)\right)\mathcal{E}(t_o,t,x)
      -
      \rho_{o}\left(X_*(t_o;t,x)\right)\mathcal{E}_*(t_o,t,x)
    }
    \d{x}
    \tag*{[\mbox{By~\eqref{eq:13}}]}
    \\
     & \leq
    \int_{\Su} \modulo{
      \rho_o\left(X(t_o;t,x)\right)
      -
      \rho_{o}\left(X_*(t_o;t,x)\right)
    }
    \mathcal{E}(t_o,t,x)
    \d{x}
    \\
     & \quad
    +
    \int_{\Su} \modulo{
      \rho_o\left(X_*(t_o;t,x)\right)
    }
    \modulo{
      \mathcal{E}(t_o,t,x)
      -
      \mathcal{E}_*(t_o,t,x)
    }
    \d{x}
    \\
     & \leq
    \int_{\Su} \modulo{
      \rho_o(\xi)
      -
      \rho_{o}\left(X_*\left(t_o;t,X(t;t_o,\xi)\right)\right)
    }
    \d{\xi}
    \tag*{[\mbox{By~\cite[eqn.~(3.35)]{opera-seconda}} and~\ref{item:14} in \Cref{lem:characteristic}]}
    \\
     & \;\,
    {+}
    \!\!\int_{\Su} \modulo{
      \rho_o\left(X_*(t_o;t,x)\right)
    }
    \d{x}
    \left(
    \norma{\partial_{xx}^2 v}_{\L\infty([0,T]\times \Su; \reali)} (t{-}t_o)
      {+} 1
    \right) e^{2 L (t{-}t_o)}
    \norma{v {-} v_*}_{\L1([t_o,t]; \W1\infty(\Su; \reali))}
    \\
     & \leq
    2 \esssup_{x \in \reali}{\modulo{
        \widetilde{X}_*(t_o;t, x)- \widetilde{X}(t_o;t,x)}
    } \;
    \tv(\rho_o)
    \tag*{[\mbox{By \Cref{lem:BVper}}]}
    \\
     & \quad
    +
    \norma{\rho_o}_{\L1(\Su; \reali)}
    \,
    \left(
    \norma{\partial_{xx}^2 v}_{\L\infty([0,T]\times \Su; \reali)}\, T
    + 1
    \right) e^{3 L (t-t_o)}
    \norma{v - v_*}_{\L1([t_o,t]; \W1\infty(\Su; \reali))}
    \\
     & \leq
    2 \,
    \tv(\rho_o) \,
    e^{L (t-t_o)} \norma{v - v_*}_{\L1([t_o,t]; \L\infty(\Su; \reali))}
    \tag*{[\mbox{By~\ref{item:14} in \Cref{lem:characteristic}}]}
    \\
     & \quad
    +
    \norma{\rho_o}_{\L1(\Su; \reali)}
    \,
    \left(
    \norma{\partial_{xx}^2 v}_{\L\infty([0,T]\times \Su; \reali)}\, T
    + 1
    \right) e^{3 L (t-t_o)} \,
    \norma{v - v_*}_{\L1([t_o,t]; \W1\infty(\Su; \reali))}\,,
  \end{align*}
  completing the estimate of the $\L1$ distance. Using similar
  computations (see for instance~\cite[Proof
    of~Lemma~3.7]{opera-seconda}), we get also the estimate on the
  $\L1(\Su; \reali)$ distance of the derivatives
  \begin{eqnarray*}
    &
    & \norma{\partial_x \rho(t) - \partial_x \rho_*(t)}_{\L1(\Su; \reali)}
    \\
    & \leq
    & \mathcal{C}
    \left(\!\!
    \begin{array}{c}
        \norma{\partial_{xx}^2 v}_{\L\infty([0,T]; \W1\infty(\Su; \reali))},
        \\
        \norma{\partial_{xx}^2 v_*}_{\L\infty([0,T]; \L\infty(\Su; \reali))},
        \\
        \norma{\rho_o}_{\cBV{1}(\Su; \reali)}, T
      \end{array}\!\!
    \right)
    e^{4 L (t-t_o)}
    \norma{v - v_*}_{\L1([t_o,t]; \W2\infty(\Su; \reali))}\, ,
  \end{eqnarray*}
  completing the proof.
\end{proofof}

\begin{proofof}{Theorem~\ref{thm:4}}
  We assume below that $t_o = 0$, the general case being entirely
  similar.  It is not restrictive to assume that
  $\norma{\rho_o}_{\L1 (\Su; \reali^k)} \leq \mathrm{M}$.  The proof
  of~\ref{item:4} relies on a fixed point argument in
  \begin{displaymath}
    \mathcal{D}_{T^*} \coloneqq
    \left\{ \rho \in \C0\left([0,T^*]; \L1(\Su; \reali^k)\right) \colon
    \norma{\rho}_{\C0([0,T^*]; \L1(\Su; \reali^k))} \leq
    \mathrm{M}
    \right\} \,,
  \end{displaymath}
  for a suitable $T^* \in \mathopen]0, T]$, equipped with the
  norm~\eqref{eq:32}. Introduce the set
  \begin{displaymath}
    \mathcal{V}_T
    \coloneqq
    \C0\left([0,T]\times \Su; \reali^k\right) \cap
    \L\infty\left([0,T]; \W2\infty(\Su; \reali^k)\right)
  \end{displaymath}
  equipped with the norm
  \begin{displaymath}
    \norma{v}_{\mathcal{V}_T}
    \coloneqq
    \esssup_{t \in [0,T]} \norma{v (t)}_{\W1\infty (\Su; \reali^k)}
  \end{displaymath}
  and define the maps
  \begin{equation}
    \label{eq:29}
    \Pi_V \colon
    \C0\left([0,T]; \L1(\Su; \reali^k)\right)
    \to
    \mathcal{V}_T
    \quad \mbox{ and } \quad
    \Sigma_{\rho_o}:
    \mathcal{V}_T
    \to
    \C0\left([0,T]; \L1(\Su; \reali^k)\right)
  \end{equation}
  where
  \begin{displaymath}
    (\Pi_V \rho) (t,x)
    \coloneqq
    \left[
      \begin{array}{c}
        V^1\left(t,x,(\rho(t)*\eta^1)(x)\right)
        \\
        \vdots
        \\
        V^k\left(t,x,(\rho(t)*\eta^k)(x)\right)
      \end{array}\right]
  \end{displaymath}
  and for all $i \in \{1, \ldots, k\}$, $(\Sigma_{\rho_o} v)^i$ is the
  solution to~\eqref{eq:2} with initial datum $\rho_o^i$ and speed
  $v^i$.

  \paragraph{\underline{Step~1}: $\Pi_V$ is well-defined and bounded.}
  For all $\rho \in \C0\left([0,T]; \L1(\Su; \reali^k)\right)$,
  $\Pi_V\rho$ is well-defined on $[0,T] \times \Su$ due to the
  definition of convolution product~\eqref{eq:14}.

  For all $\rho \in \mathcal{D}_T$, \Cref{lem:0} implies that for all
  $t \in [0,T]$
  \begin{displaymath}
    \max_{i = 1, \ldots, k}
    \norma{\rho(t)*\eta^i}_{\L\infty(\Su; \reali^k)}
    \leq
    \mathrm{M}_*
    \quad \mbox{ where } \quad
    \mathrm{M}_*
    =
    \norma{\eta}_{\L\infty(\reali; \reali^{k \times k})}\, H \, \mathrm{M} \,.
  \end{displaymath}
  The regularities of $V$ and $\eta$, through long but elementary
  computations analogous to the ones for instance in~\cite[Proof of
    Lemma~4.3]{opera-prima}, show that
  $\Pi_V\rho \in \L\infty\left([0,T]; \W2\infty(\Su; \reali^k)\right)$
  and
  \begin{equation}\label{eq:30}
    \!\!\!\!
    {
    \norma{\Pi_V \rho}_{\L\infty([0,T]; \W2\infty(\Su; \reali^k))}}
    \leq
    \underbrace{
      \mathcal{C}
      \left(
      \norma{V}_{\L\infty([0,T]; \W2\infty(\Su \times B_{\mathrm{M}_*}; \reali^k))},
      \norma{\eta}_{\W2\infty(\reali; \reali^{k \times k})},
      H,
      \mathrm{M}
      \right)
    }_{\mathcal{B}(\mathrm{M})}
    ,
  \end{equation}
  proving boundedness.

  Moreover, $\Pi_V \rho \in \C0([0,T]\times \Su; \reali^k)$. In fact,
  the assumption
  $V \in \C0([0,T]\times \Su \times \reali^k; \reali^k)$ ensures by
  Heine-Cantor Theorem (see~\cite[Lemma~8.4]{Folland}) that $V$ is
  uniformly continuous on $[0,T] \times \Su \times B_{\mathrm{M}_*}$,
  namely
  \begin{align}\label{eq:55}
    \forall \, \epsilon >0 \quad
    \exists \, \delta= \delta (\epsilon) >0
    \, \text{ such that } \,
    \forall (t,x,q), (t',x',q') \in [0,T] \times \Su \times B_{\mathrm{M}_*},
    \\
    \nonumber
    \norma{(t',x',q')-(t,x,q)} \leq \delta \implies
    \norma{V(t',x',q')-V(t,x,q)} \leq \epsilon.
  \end{align}
  As $\eta \in \C0(\reali; \reali^{k \times k})$, for all
  $t \in [0,T]$ the map $\rho(t)*\eta \in \C0(\Su; \reali^k)$ and
  there exists $\delta_1 = \delta_1(t, \delta) >0$ such that
  \begin{displaymath}
    \forall x, x' \in \Su \quad \modulo{x'-x} \leq \delta_1
    \implies
    \norma{\left( \rho(t)*\eta \right)(x) - \left( \rho(t)*\eta \right)(x')}
    \leq \delta/2 \,.
  \end{displaymath}
  On the other hand,
  $\rho \in \C0\left([0,T]; \L1(\Su; \reali^k)\right)$, hence there
  exists $\delta_2 = \delta_2(\delta) >0$ such that by
  Lemma~\ref{lem:0}
  \begin{align*}
    \forall t, t' \in [0,T] \quad
    \modulo{t'-t} \leq \delta_2
    \implies
     & \norma{\rho(t) - \rho(t')}_{\L1(\Su; \reali^k)}
    \leq
    \frac{\delta}{2 \, H \, \norma{\eta}_{\L\infty(\reali; \reali^{k\times k})}}
    \\
    \implies
     & \norma{\rho(t)*\eta - \rho(t')*\eta}_{\L\infty(\Su; \reali^k)}
    \leq
    \delta/2 \,.
  \end{align*}
  Hence, for all $(t,x), (t',x') \in [0,T] \times \Su$ such that
  $\modulo{t-t'} \leq \min\left\{\delta, \delta_2\right\}$ and
  $\modulo{x-x'} \leq \min\left\{\delta, \delta_1\right\}$
  \begin{align*}
     & \norma{\left( \rho(t)*\eta \right)(x)
      - \left( \rho(t')*\eta \right)(x')}
    \\
     & \leq
    \norma{\left( \rho(t)*\eta \right)(x) - \left( \rho(t)*\eta \right)(x')}
    +
    \norma{\left( \rho(t)*\eta \right)(x') - \left( \rho(t')*\eta \right)(x')}
    \\
     & \leq
    \delta
  \end{align*}
  and, by~\eqref{eq:55},
  $\norma{V\left(t',x', \rho(t')*\eta(x')\right)-V\left(t,x,
      \rho(t)*\eta(x)\right)} \leq \epsilon$, proving
  \textbf{\underline{Step~1}}.

  \paragraph{\underline{Step~2}: $\Pi_V$ is locally Lipschitz
    continuous in $\rho$.}

  For all the functions $\rho, \rho_* \in \mathcal{D}_T$, the same
  computations as in~\cite[Lemma~4.3]{opera-prima} yield
  \begin{align*}
     & \quad
    \norma{\Pi_V \rho - \Pi_V \rho_*}_{\L\infty([0,T]; \W1\infty(\Su; \reali^k))}
    \\
     & \leq
    \mathcal{C}
    \left(
    \begin{array}{c}
        \norma{\nabla_\rho V}_{\L\infty([0,T]; \W1\infty(\Su \times B_{\mathrm{M}_*}; \reali^{k\times k}))},
        \\
        \norma{\eta}_{\W1\infty(\reali; \reali^{k \times k})},
        H,
        \mathrm{M}
      \end{array}
    \right)
    \norma{\rho - \rho_*}_{\C0([0,T]; \L1(\Su; \reali^k))}\,.
  \end{align*}

  \paragraph{\underline{Step~3}: $\Pi_V$ is locally Lipschitz
    continuous in $V$.}

  For all $\rho \in \mathcal{D}_T$ and for all
  $V, V_* \in \C0([0,T]\times \Su \times \reali^k; \reali^k) \cap
    \L\infty\left([0,T]; \Wloc2\infty(\Su \times \reali^k;
    \reali^k)\right)$, by direct computations,
  \begin{displaymath}
    \norma{\Pi_V \rho - \Pi_{V_*} \rho}_{\L\infty([0,T]; \W1\infty(\Su; \reali^k))}
    \leq
    \mathcal{C}
    \left(\!
    \begin{array}{@{}c@{}}
        \norma{\eta'}_{\L\infty(\reali; \reali^{k \times k})},
        \\
        H , \mathrm{M}
      \end{array}\!
    \right)
    \norma{V - V_*}_{\L\infty([0,T]; \W1\infty(\Su \times B_{\mathrm{M}_*}; \reali^{k}))}.
  \end{displaymath}

  \paragraph{\underline{Step~4}: $\Sigma_{\rho_o}$ is well-defined,
    bounded and locally Lipschitz continuous.}

  For $v \in \mathcal{V}_T$, set $\rho \coloneqq \Sigma_{\rho_o}
    v$. Hence, \ref{item:1}, \ref{item:2}, \ref{item:9}
  and~\ref{item:3a} with $\ell=0$ in \Cref{lem:1} apply, proving local
  Lipschitz continuity. More precisely, for later use, we have:
  \begin{align}
    \nonumber
     &
    \norma{(\Sigma_{\rho_o} v) (t) - (\Sigma_{\rho_o} v_*) (t)}_{\L1 (\Su; \reali^k)}
    \\
    \label{eq:21}
    \leq
     &
    \underbrace{
      \mathcal{C}
      \left(
      \begin{array}{@{}c@{}}
          \norma{\partial_{x} v}_{\L\infty([0,T]; \W{1}\infty(\Su; \reali^k))}, \\
          \norma{\partial_{x} v_*}_{\L\infty([0,T]; \L\infty(\Su; \reali^k))},
          \\
          \norma{\rho_o}_{\cBV{}(\Su; \reali^k)}, T
        \end{array}
      \right)
    }_{\Lip (\Sigma_{\rho_o})}
    \, t \,
    \norma{v - v_*}_{\L\infty([0,t]; \W{1}{\infty}(\Su; \reali^k))} .
  \end{align}

  \paragraph{\underline{Step~5}: $\mathcal{T}$ is well-defined and it
    is a contraction for a sufficiently small time.}

  Define $\mathcal{T}: \mathcal{D}_T \to \mathcal{D}_T$ by
  $\mathcal{T} = \Sigma_{\rho_o} \circ \Pi_V$, with the
  notation~\eqref{eq:29}. The well-posedness and Lipschitz continuity
  of $\mathcal{T}$ follow by \textbf{\underline{Step~1}},
  \textbf{\underline{Step~2}} and \textbf{\underline{Step~4}}.
  Consider $\rho, \rho_* \in \mathcal{D}_T$ and call $v = \Pi_V \rho$
  and $v_* = \Pi_V \rho_*$.  By~\ref{item:3a} in \Cref{lem:1}, the
  bound~\eqref{eq:30} and \textbf{\underline{Step~2}}, for all
  $t \in [0,T]$ we have
  \begin{align*}
    \norma{\mathcal{T} \rho - \mathcal{T} \rho_*}_{\C0([0,t]; \L1( \Su; \reali^k))}
     & \leq
    \mathcal{C}\left(
    \begin{array}{c}
        \norma{v}_{\L\infty([0,T]; \W2\infty(\Su; \reali^k))},
        \\
        \norma{v_*}_{\L\infty([0,T]; \W1\infty(\Su; \reali^k))},
        \\
        \norma{\rho_o}_{\cBV{}(\Su; \reali^k)},
        T
      \end{array}
    \right)
    t \,
    \norma{v -  v_*}_{\L\infty([0,t]; \W1\infty(\Su; \reali^k))}
    \\
     &
    \leq
    \underbrace{\mathcal{C}\!\left(\!
      \begin{array}{c}
          \norma{V}_{\L\infty([0,T]; \W2\infty(\Su \times B_{\mathrm{M}_*}; \reali^k))},
          \\
          \norma{\eta}_{\W2\infty(\reali; \reali^{k \times k})},
          H,
          \\
          \norma{\rho_o}_{\cBV{}(\Su; \reali^k)},
          T
        \end{array}\!
      \right)}_{\mathcal{L} \left(\norma{\rho_o}_{\cBV{}(\Su; \reali^k)}\right)}
    \; t \,
    \norma{\rho - \rho_*}_{\C0([0,t]; \L1(\Su; \reali^k))}.
  \end{align*}
  Hence, choosing $T^* \in \mathopen]0,T]$ such that
  $\mathcal{L} \left(\norma{\rho_o}_{\cBV{}(\Su; \reali^k)}\right) \,
    T^* <1$, $\mathcal{T}$ is a contraction on $\mathcal{D}_{T^*}$.

  \paragraph{\underline{Step~6}: Proof of~\ref{item:35}.}

  By Definition~\ref{def:solNonloc}, any solution $\rho$
  to~\eqref{eq:25} is a fixed point of $\mathcal{T}$, hence,
  by~\ref{item:9} in \Cref{lem:1},
  \begin{displaymath}
    \norma{\rho^i(t)}_{\L1(\Su; \reali)} = \norma{\rho^i_o}_{\L1(\Su; \reali)}
    \qquad \mbox{for all } t \in [0,T].
  \end{displaymath}
  By~\eqref{eq:30}, concerning the constant $\const{3}$
  in~\ref{item:35}, define
  $\const{3} \coloneq \mathcal{B}\left( \norma{\rho_o}_{\L1(\Su;
        \reali^k)} \right)$ and observe that
  $\norma{\Pi_V \rho}_{\L\infty([0,t]; \W2\infty(\Su; \reali^k))} \leq
    \mathcal{B}\left( \norma{\rho(t)}_{\L1(\Su; \reali^k)} \right) =
    \mathcal{B}\left( \norma{\rho_o}_{\L1(\Su; \reali^k)} \right) =
    \const{3}$.

  At last, \ref{item:7} and~\ref{item:5} in \Cref{lem:1} yield that
  $\mathcal{T}\rho(t) =\rho(t) \in \BV(\Su; \reali^k)$ and
  \begin{align*}
    \nonumber
    \norma{\rho(t)}_{\L\infty(\Su; \reali^k)}
     & \leq
    \norma{\rho_o}_{\L\infty(\Su; \reali^k)}
    e^{
    \norma{\partial_x(\Pi_V \rho)}_{\L\infty([0,t]\times \Su; \reali^k)}\, t}
    \\
     & \leq
    \norma{\rho_o}_{\L\infty(\Su; \reali^k)}
    \, e^{\const{3}\, t}
    \\
    \nonumber
    \tv\left(\rho(t)\right)
     & \leq \!
    \left(\!
    \tv(\rho_o)
    \!+\!
    \norma{\partial_{xx}^2 (\Pi_V \rho)}_{\L\infty([0,t]\times \Su; \reali^k)}
    \norma{\rho_o}_{\L1(\Su; \reali^k)}
    \,
    t\!
    \right)
    e^{
    \norma{\partial_x(\Pi_V \rho)}_{\L\infty([0,t]\times \Su; \reali^k)}
    \, t
    }
    \\
     & \leq
    \left(
    \tv(\rho_o)
    +
    \const{3} \,
    \norma{\rho_o}_{\L1(\Su; \reali^k)}
    \,
    t
    \right)
    e^{
        \const{3}
        \, t
      },
  \end{align*}
  proving~\ref{item:35}.

  \paragraph{\underline{Step~7}: Proof of~\ref{item:4}$\cdots$\ref{item:12}.}
  By means of Banach Fixed Point Theorem, \textbf{\underline{Step~5}}
  ensures that the map $\mathcal{T}$ admits a unique fixed point
  $\rho \in \C0\left([0,T^*]; \L1(\Su; \reali^k)\right)$ which is the
  unique solution to~\eqref{eq:25} on $[0,T^*]$.  Moreover, the a
  priori bounds in \textbf{\underline{Step~6}} guarantee that there
  exists a constant
  \begin{displaymath}
    K
    \coloneqq \mathcal{C} \left(
    \norma{V}_{\L\infty([0,T]; \W2\infty(\Su \times B_{\mathrm{M}_*}; \reali^k))},
    \norma{\eta}_{\W2\infty(\reali;
      \reali^{k\times k})}, H, \norma{\rho_o}_{\cBV{}(\Su;
      \reali^k)} , T \right)
  \end{displaymath}
  such that for all $t \in [0,T]$,
  $ \norma{\rho(t)}_{\cBV{}(\Su; \reali^k)} \leq K$. Hence, selecting
  the time $T^*$ so that $\mathcal{L}(K) \, T^* <1$, ensures that
  $T^*$ depends only on the parameters of the problem and on the
  initial datum. We can iterate the argument above obtaining the
  solution to~\eqref{eq:25} on $[T^*, 2T^*]$, $[2T^*, 3T^*]$,
  $\ldots$, covering the whole time interval $[0,T]$ and completing
  the proof of~\ref{item:4}.

  The proof of~\ref{item:10} follows the same lines as
  in~\cite[\textbf{(1)} of Theorem~2.2]{MR4894089}. The proof
  of~\ref{item:12} follows from~\eqref{eq:13}.

  \paragraph{\underline{Step~8}: Proof of~\ref{item:38}.}
  Define the functions
  $v(t,x) \coloneqq (\Pi_V \mathcal{P}_{0,t}\rho_o) (x)$ and
  $v_*(t,x) \coloneqq (\Pi_{V} \mathcal{P}_{0,t} \rho_{o,*}) (x)$.
  By~\ref{item:6} and~\ref{item:3a} in \Cref{lem:1}, for all
  $t \in [0,T]$,
  \begin{align*}
     & \norma{\mathcal{P}_{0,t} \rho_o - \mathcal{P}_{0,t} \rho_{o,*}}_{\L1(\Su; \reali^k)}
    \\
     &
    \leq
    \norma{\rho_o - \rho_{o,*}}_{\L1(\Su; \reali^k)}
    +
    \mathcal{C}
    \left(
    \begin{array}{c}
        \norma{\partial_x v}_{\L\infty([0,T]; \W1\infty(\Su; \reali^k))},
        \\
        \norma{\partial_x v_*}_{\L\infty([0,T]; \L\infty(\Su; \reali^k))},
        \\
        \norma{\rho_o}_{\cBV{}(\Su; \reali^k)},
        t
      \end{array}
    \right)
    \norma{v - v_*}_{\L1([0,t]; \W1\infty(\Su; \reali^k))}
    \\
     &
    \leq
    \norma{\rho_o - \rho_{o,*}}_{\L1(\Su; \reali^k)}
    \\
     & \quad +
    \mathcal{C}
    \left(
    \begin{array}{c}
        \norma{V}_{\L\infty([0,T]; \W2\infty(\Su \times B_{\mathrm{M}_*}; \reali^k))},
        \\
        \norma{\eta}_{\W2\infty(\reali; \reali^{k \times k})},
        H,
        \\
        \norma{\rho_o}_{\cBV{}(\Su; \reali^k)},
        \norma{\rho_{o,*}}_{\L1(\Su; \reali^k)},
        t
      \end{array}
    \right)
    \int_0^t
    \norma{
      \mathcal{P}_{0,\tau} \rho_o -
      \mathcal{P}_{0,\tau} \rho_{o,*}}_{\L1(\Su; \reali^k)} \d\tau
    \tag*{[\mbox{By~\eqref{eq:30}}]}\, .
  \end{align*}
  An application of the Gronwall Lemma concludes the proof
  of~\ref{item:38}.

  \paragraph{\underline{Step~9}: Proof of~\ref{item:43}.}
  We follow the same strategy as in the proof
  of~\cite[$(\mathcal{P}5)$ in Theorem~2.3]{opera-prima}, based
  on~\cite[Theorem~2.9]{MR1816648}. Then,
  \begin{equation}
    \label{eq:22}
    \!\!\!\!\!
    \norma{\mathcal{P}_{0,t}\rho_o -
      \mathcal{P}^*_{0,t}\rho_o}_{\L1 (\Su;\reali^k)}
    \leq
    e^{\const{4} T}
    \int_0^t \liminf_{h\to 0^+} \!
    \dfrac{\norma{\mathcal{P}_{\tau,\tau+h}
        \mathcal{P}^*_{0,\tau} \rho_o -
        \mathcal{P}^*_{\tau,\tau+h} \mathcal{P}^*_{0,\tau}
        \rho_o}_{\L1 (\Su;\reali^k)}}{h} \d\tau \,
  \end{equation}
  where $\const{4}$ is as in~\ref{item:38}.

  For $\tau \in [0,t]$, call $\rho = \mathcal{P}^*_{0,\tau}
    \rho_o$. For $h$ sufficiently small, by \underline{\textbf{Step~5}}
  the map $\mathcal{T} \colon \mathcal{D}_h \to \mathcal{D}_h$, given
  by $\mathcal{T} = \Sigma_\rho \circ \Pi_V$, is a contraction and is
  a Lipschitz continuous function also of the speed $V$
  by~\underline{\textbf{Step~3}}. Hence, setting
  $v (t,x) \coloneq \Pi_V \rho(t,x)$ and
  $v^*(t,x) \coloneq \Pi_{V_*} \rho(t,x)$
  \begin{eqnarray*}
    &
    & \norma{\mathcal{P}_{\tau, \tau+h}\rho -
      \mathcal{P}_{\tau, \tau+h}^* \rho}_{\L1(\Su; \reali^k)}
    \\
    & \leq
    & \dfrac{\Lip (\Sigma_\rho) \;  \mathcal{C}
      \left(\norma{\eta'}_{\L\infty(\reali; \reali^{k \times k})},
      H , \mathrm{M}
      \right)
    }{1-\mathcal{L} \left(\norma{\rho_o}_{\cBV{}(\Su; \reali^k)}\right) \,  T^*}
    \, \norma{V-V_*}_{\L\infty([0,T]; \W1\infty(\Su \times B_{\mathrm{M}_*}; \reali^{k}))}
    \\
    & \leq
    &     \mathcal{C}
    \left(
    \begin{array}{@{}c@{}}
        \norma{\partial_{x} v}_{\L\infty([0,T]; \W{1}\infty(\Su; \reali^k))}, \\
        \norma{\partial_{x} v^*}_{\L\infty([0,T]; \L\infty(\Su; \reali^k))},
        \\
        \norma{\rho_o}_{\cBV{}(\Su; \reali^k)}, T
      \end{array}
    \right)
    \dfrac{\mathcal{C}
      \left(\norma{\eta'}_{\L\infty(\reali; \reali^{k \times k})},
      H , \mathrm{M}
      \right)
    }{1-\mathcal{L} \left(\norma{\rho_o}_{\cBV{}(\Su; \reali^k)}\right) \,  T^*}
    \, h
    \\
    && \times
    \norma{V-V_*}_{\L\infty([0,T]; \W1\infty(\Su \times B_{\mathrm{M}_*}; \reali^{k}))}
    \\
    & \leq
    & \mathcal{C}
    \left(
    \begin{array}{@{}c@{}}
        \norma{V}_{\L\infty([0,T]; \W2\infty(\Su \times B_{\mathrm{M}_*}; \reali^k))},
        \\
        \norma{V_*}_{\L\infty([0,T]; \W1\infty(\Su \times B_{\mathrm{M}_*}; \reali^k))},
        \\
        \norma{\eta}_{\W2\infty(\reali; \reali^{k \times k})},
        H,
        \norma{\rho_o}_{\cBV{}(\Su; \reali^k)}, T
      \end{array}
    \right) \, h \,
    \norma{V-V_*}_{\L\infty([0,T]; \W1\infty(\Su \times B_{\mathrm{M}_*}; \reali^{k}))}
  \end{eqnarray*}
  where we used~\eqref{eq:21} and~\eqref{eq:30}.  Inserting the latter
  expression above into~\eqref{eq:22}, we obtain~\eqref{eq:24},
  completing the proof to \Cref{item:43}.

  \paragraph{\underline{Step~10}: Proof of~\ref{item:44}.}
  Using the notation in the proof
  of~\ref{item:4}-\ref{item:12}-\ref{item:35}, call
  $v(t,x) \coloneqq (\Pi_V \mathcal{P}_{0,t}\rho_o) (x)$.  An
  application of~\ref{item:8} in \Cref{lem:1} with
  $L = \norma{v}_{\L\infty([0,T]; \W1\infty(\Su; \reali^k))}$ and
  estimate~\eqref{eq:30} yield
  \begin{align*}
    \norma{\mathcal{P}_{0,t}\rho_o - \mathcal{P}_{0,t_*}\rho_o}_{\L1(\Su; \reali)}
     & \leq
    \mathcal{C}\left(
    \begin{array}{c}
        \norma{v}_{\L\infty([0,T]; \W2\infty(\Su; \reali))},
        \norma{\rho_o}_{\cBV{}(\Su; \reali^k)}, T
      \end{array}
    \right)  \modulo{t-t_*}
    \\
     & \leq
    \mathcal{C}\left(
    \begin{array}{c}
        \norma{V}_{\L\infty([0,T]; \W2\infty(\Su \times B_{\mathrm{M}_*}; \reali^k))},
        \\
        \norma{\eta}_{\W2\infty(\reali; \reali^{k \times k})},
        H,
        \norma{\rho_o}_{\cBV{}(\Su; \reali^k)}, T
      \end{array}
    \right)  \modulo{t-t_*},
  \end{align*}
  proving~\ref{item:44}.
\end{proofof}

\begin{proofof}{Proposition~\ref{prop:fine_estimates}}
  For all $i \in \{1, \ldots, k\}$, call
  $v^i(t,x) \coloneqq V^i\left(t,x, \left(\rho(t)*\eta^i\right)(x)
    \right)$ and recall that for all $t \in [0,T]$,
  $\norma{\rho(t)}_{\L1(\Su; \reali^k)} = \norma{\rho_o}_{\L1(\Su;
      \reali^k)}$ by~\ref{item:35} in \Cref{thm:4}. The regularity
  assumptions on $V$ and on $\eta$ imply that
  \begin{displaymath}
    \begin{array}{c}
      v \in \L\infty\left([0,T]; \W3\infty(\Su;\reali^k)\right),
      \quad
      \partial_x v \in \C0([0,T]\times \Su; \reali^k),
      \\
      \partial^2_{xx} v (t, \cdot) \in \C1 (\Su; \reali^k)
      \quad \mbox{for all}\, t \in [0,T] \, .
    \end{array}
  \end{displaymath}
  In particular, computations analogous to those in
  \textbf{\underline{Step~1}} in the proof of \Cref{thm:4} yield
  \begin{align}
    \label{eq:61}
    \norma{v}_{\L\infty([0,t]; \W{\ell+1}\infty(\Su; \reali^k))}
    \leq
    \mathcal{C} \left(
    \begin{array}{c}
        \norma{V}_{\L\infty([0,T]; \W{\ell+1}\infty(\Su \times B_{\mathrm{M}_*}; \reali^k))},
        \\
        \norma{\eta}_{\W{\ell+1}\infty(\reali; \reali^{k \times k})},
        H,
        \norma{\rho_o}_{\L1(\Su; \reali^k)}
      \end{array}
    \right)
    \qquad \ell=1,2 .
  \end{align}
  Hence,~\ref{item:37a} in \Cref{lem:1} apply, proving
  $\rho \in \C1([0,T] \times \Su; \reali^k)$,
  $\partial_{xx}^2\rho \in \C0([0,T]\times \Su; \reali^k)$ and
  $\partial_{tx}^2 \rho \in \L\infty([0,T]\times \Su; \reali^k)$.
  Inserting~\eqref{eq:61} in~\eqref{eq:62} and~\eqref{eq:66} concludes
  the proof of~\ref{item:30}. The proof of~\ref{item:36} is a
  straightforward adaptation to the periodic setting
  of~\cite[$(\mathbf{N}3)$ in Lemma~3.5]{opera-seconda}, using
  \Cref{lem:BVper}.
\end{proofof}

\subsection{Proofs Related to the Mixed Nonlocal-Local Case}
\label{sec:proof-WP}

The following result provides a stability property for
Definition~\ref{def:sol} under perturbations of both the initial data
and the coupling terms.  We omit the proof since it is a
straightforward adaptation of~\cite[Proposition~2.2]{opera-seconda}.

\begin{proposition}[Stability of Definition~\ref{def:sol}]
  \label{lem:GoodSol}
  Assume~\ref{item:17}, \ref{item:21}, and~\ref{item:hyp-eta}
  hold. Fix $T>0$ and
  $(\rho_{o,n}, r_{o,n}) \in \W11 (\Su;\reali^k) \times \L1
    (\Su;\reali)$,
  $(\epsilon_n,e_n) \in \C0 \left([0,T];\W11 (\Su; \reali^k) \!\times\!
    \L1 (\Su;\reali)\right)$, and a function $(\rho_\infty, r_\infty)$
  belonging to
  $ \C0\left([0,T];\W11 (\Su;\reali^k) \times \L1 (\Su;\reali)\right)$
  such that for all $n \in \naturali$ the problem
  \begin{displaymath}
    \begin{cases}
      \partial_t\rho_n^i
      +
      \partial_x\left(
      \rho_n^i \,
      v^i_{NL}\left(
      \Sigma \rho_n * \eta^i
      +
      (r_n + e_n)  * \eta^i_{k+1}
      \right)\right)=0 \qquad i \in \{1, \ldots, k\}
      \\
      \partial_t r_n
      +
      \partial_x\left(
      r_n \, v_L\left( \Sigma (\rho_n+\epsilon_n) + r_n \right)
      \right)=0
      \\
      \rho_n(0,x) = \rho_{o,n}(x)
      \\
      r_n(0,x)=r_{o,n}(x)
    \end{cases}
  \end{displaymath}
  admits as solution $(\rho_n,r_n)$ in the sense of
  Definition~\ref{def:sol}. Moreover, assume that
  \begin{equation}
    \label{eq:72}
    \begin{array}{@{}rcll@{}}
      \displaystyle
      \lim_{n\to+\infty} (\rho_n,r_n)
       & =
       & (\rho_\infty,r_\infty)
       & \mbox{ in }
      \C0 \left([0,T]; \W11(\Su;\reali^k) \times
      \L1(\Su;\reali)\right)\,;
      \\
      \displaystyle
      \lim_{n\to+\infty} (\epsilon_n,e_n)
       & =
       & (0,0)
       & \mbox{ in }
      \C0 \left([0,T]; \W11(\Su;\reali^k) \times
      \L1(\Su;\reali)\right)\,.
    \end{array}
  \end{equation}
  Then, $(\rho_\infty,r_\infty)$ solves, in the sense of
  Definition~\ref{def:sol}, the problem
  \begin{equation}
    \begin{cases}
      \partial_t\rho^i
      +
      \partial_x\left(
      \rho^i \,
      v^i_{NL}\left(
      \Sigma \rho * \eta^i
      +
      r * \eta^i_{k+1}
      \right)\right)=0 \qquad i \in \{1, \ldots, k\}
      \\
      \partial_t r
      +
      \partial_x\left(
      r \, v_L\left( \Sigma \rho + r \right)
      \right)=0
      \\
      \rho(0,x) = \rho_\infty(0,x)
      \\
      r(0,x)=r_\infty(0,x) \,.
    \end{cases}
  \end{equation}
\end{proposition}

We present some preliminary results on the interplay between solutions
to systems of nonlocal conservation laws and solutions to \emph{local}
conservation laws.  {
    In particular, we verify that the regularity
    assumptions required by the well-posedness results
    of~\Cref{sec:kruv-zkov-periodic} and~\Cref{sec:nonl-peri-case-2}  are
    indeed satisfied by the velocity functions arising in the mixed system
    \eqref{eq:3}.}

\begin{lemma}[Regularity of $V$ and $w$]
  \label{lem:buono}
  Assume~\ref{item:hyp-eta},\ref{item:17} and~\ref{item:21}.
  \begin{enumerate}[label=\bf(\roman*)]
    \item \label{item:22} For all
          $r \in \C0\left([0,T]; \L1(\Su; \reali)\right)$, the map
          \begin{displaymath}
            V^i(t,x,q)\coloneqq v^i_{NL}\left(\Sigma q +
            \left(\, r(t)*\eta^i_{k+1}\right)\,(x)\right) \quad (t,x,q) \in [0,T]\times \Su \times \reali^k,
            \,\, i \in \{1, \ldots, k\}\,.
          \end{displaymath}
          satisfies~\ref{eq:87} in \Cref{prop:fine_estimates} and
          \begin{equation}
            \label{eq:40}
            \norma{V}_{\L\infty([0,t]; \W3\infty(\Su \times \reali^k; \reali^k))}
            \leq
            \mathcal{C}
            \left(
            \begin{array}{c}
                \norma{v_{NL}}_{\W3\infty(\reali; \reali^k)},
                \norma{\eta_{k+1}}_{\W3\infty(\reali; \reali^k)},
                \\
                H,
                \norma{r}_{\C0([0,t]; \L1(\Su; \reali))}
              \end{array}
            \right)\!\, .
          \end{equation}
    \item \label{item:11} For all
          $\rho \in \C1([0,T]\times \Su; \reali_+^k)$ satisfying
          $\partial_{xx}^2 \rho \in \C0([0,T] \times \Su; \reali^k)$ and
          $\partial_{tx}^2\rho \in \L\infty([0,T]\times \Su; \reali^k)$, the
          map
          \begin{displaymath}
            w(t,x,r) \coloneqq v_L\left(\Sigma \rho(t,x) + r\right)
            \qquad
            (t,x,r) \in [0,T] \times \Su \times \reali
          \end{displaymath}
          satisfies~\ref{item:26}$\cdots$\ref{item:31} as soon as
          $R \geq \max\left\{ R_L, \norma{\rho}_{\L\infty ([0,T] \times\Su;
            \reali^k)} \right\}$ and
          \begin{align}
            \label{eq:51}
            \norma{w}_{\L\infty([0,T]\times \Su; \W1\infty(\reali; \reali ))}
             & \leq
            \norma{v_L}_{\W1\infty(\reali; \reali)};
            \\
            \label{eq:52}
            \norma{\partial_x w}_{\L\infty([0,T] \times \Su; \W1\infty( \reali; \reali ))}
             & \leq
            \mathcal{C}\left(
            \norma{v_L'}_{\W1\infty(\reali; \reali)},
            \norma{\partial_x\rho}_{\L\infty([0,T]\times \Su; \reali^k)}
            \right);
            \\
            \label{eq:63}
            \!\!\!\!\!\!\!\!\!
            \esssup_{\tau \in [0,T]}
            \norma{\partial_{xx}^2 w(\tau)}_{\L1(\Su; \L\infty(\reali; \reali ))}
             & \leq
            \mathcal{C}\left(
            \norma{v_L'}_{\W1\infty(\reali; \reali)},
            \norma{\partial_x\rho}_{\L\infty([0,T]; \W11(\Su; \reali^k))}
            \right) \,.
          \end{align}
  \end{enumerate}
\end{lemma}

\begin{proofof}{Lemma~\ref{lem:buono}}
  Consider~\ref{item:22}.  Direct computations as in~\cite[Proof to
    Lemma 3.6]{opera-seconda} allow concluding that
  $V \in \C0\left([0,T]; \C1(\Su \times \reali^k; \reali^k)\right)
    \cap \L\infty\left([0,T]; \W3\infty(\Su \times \reali^k;
    \reali^k)\right)$, for all $t \in [0,T]$ the map
  $ V(t,\cdot, \cdot) \in \C3(\Su \times \reali^k; \reali^k)$
  and~\eqref{eq:40} holds.

  To prove~\ref{item:11}, observe that by the regularity assumptions
  on $v_L$ and $\rho$, direct computations show
  that~\ref{item:26},~\ref{item:28},~\ref{item:31} and
  \eqref{eq:51}--\eqref{eq:63} hold.

  Note that~\ref{item:27} is directly implied by~\ref{item:28}
  and~\ref{item:31}. At last, assumption~\ref{item:21} ensures that
  any choice
  $R \geq \max\left\{ R_L, \norma{\rho}_{\L\infty ([0,T] \times \Su;
    \reali^k)} \right\}$ implies that~\ref{item:34} holds.
\end{proofof}

We now establish stability estimates on the dependence of the nonlocal
(respectively, local) component of the solution to the mixed system on
an assigned local (respectively, nonlocal) map.

\begin{lemma}[Stability of the nonlocal solution]
  \label{lem:stabilita_loc_nonloc}
  Assume~\ref{item:hyp-eta},\ref{item:17} and fix $\mathrm{R} >0$.
  For all $r, r_* \in \C0\left([0,T]; \L1(\Su; \reali)\right)$ with
  \begin{displaymath}
    \norma{r}_{\C0([0,T]; \L1(\Su; \reali))} +
    \norma{r_*}_{\C0([0,T]; \L1(\Su; \reali))}
    \leq \mathrm{R},
  \end{displaymath}
  $\rho_o\in \cBV{1}(\Su; \reali^k)$ and
  $ \rho_{o,*} \in \W11(\Su; \reali^k)$, call $\rho,\rho_*$ the
  solutions, in the sense of \Cref{def:solNonloc}, to
  \begin{displaymath}
    \left\{
    \begin{array}{@{}l}
      \partial_t \rho^i {+} \partial_x
      \left(
      \rho^i v^i_{NL} \left(
      \Sigma \rho *\eta^i {+} r*\eta^{i}_{k+1}
      \right)
      \right) = 0
      \\
      \rho(0,x) = \rho_o(x)
    \end{array}
    \right.
    \,
    \left\{
    \begin{array}{@{}lr@{\,}lc@{}}
      \partial_t \rho^i_* {+} \partial_x \!
      \left(
      \rho^i_* v^i_{NL} \left(
      \Sigma \rho_* *\eta^i {+} r_* * \eta^{i}_{k+1}
      \right)
      \right) = 0
      \\
      \rho_*(0,x) = \rho_{o,*}(x) \, .
    \end{array}
    \right.
  \end{displaymath}
  Then, there exists a constant
  \begin{displaymath}
    \mathcal{C}_{5} \coloneqq
    \mathcal{C}\left(
    \norma{v_{NL}}_{\W3\infty(\reali; \reali^k)},
    \norma{\eta}_{\W3\infty(\reali; \reali^{k\times(k+1)})},
    H,
    \norma{\rho_o}_{\cBV{1}(\Su; \reali^k)},
    \norma{\rho_{o,*}}_{\cBV{}(\Su; \reali^k)},
    \mathrm{R},
    t
    \right)
  \end{displaymath}
  such that for all $t \in [0,T]$,
  \begin{displaymath}
    \norma{\rho(t) - \rho_* (t)}_{\W11(\Su; \reali^k)}
    \leq
    \left(1 + \mathcal{C}_{5} t\right)
    \norma{\rho_o -  \rho_{o,*}}_{\W11(\Su; \reali^k)}
    +
    \mathcal{C}_{5} \int_0^t \norma{r(\tau)- r_* (\tau)}_{\L1(\Su; \reali)}
    \dd \tau.
  \end{displaymath}
\end{lemma}

The present result extends~\cite[Lemma~3.7]{opera-seconda} to the
periodic setting. To this aim, we exploit the stability of solutions
to the linear continuity equation on the speed.

\begin{proofof}{Lemma~\ref{lem:stabilita_loc_nonloc}}
  This proof consists in applying~\ref{item:3a} in \Cref{lem:1} to the
  case
  \begin{displaymath}
    \begin{array}{l@{}l}
      v^i(t,x)
      \coloneqq
       &
      v^i_{NL} \left(
      \left(\Sigma \rho(t)*\eta^i\right)(x)
      +
      \left(r(t)*\eta^i_{k+1}\right)(x)
      \right)
      \\
      v^i_*(t,x)
      \coloneqq
       & v^i_{NL} \left(
      \left(\Sigma \rho_*(t)*\eta^i\right)(x)
      +
      \left(r_*(t)*\eta^i_{k+1}\right)(x)
      \right)
    \end{array}
    \quad (t,x) \in [0,T]\times \Su
    \,\, i \in \{1, \ldots, k\}\, .
  \end{displaymath}
\end{proofof}

\begin{lemma}[Stability of the local solution]
  \label{lem:stabilita_loc_nonloc_2}
  Let~\ref{item:21} hold. Fix
  $r_o, r_{o,*} \in \cBV{}(\Su; [0,R_L])$ and
  $\rho, \rho_* \in \C1([0,T] \times \Su; \reali_+^k)$ with
  $\partial_{xx}^2 \rho, \partial_{xx}^2 \rho_* \in \C0([0,T] \times
    \Su; \reali^k)$ and
  $\partial_{tx}^2 \rho, \partial_{tx}^2\rho_* \in
    \L\infty([0,T]\times \Su; \reali^k)$.  Call $r, r_*$ the solutions
  to
  \begin{equation}\label{eq:53}
    \left\{
    \begin{array}{l}
      \partial_t r + \partial_x \left(r\, v_L \left(
        \Sigma \rho + r
        \right)\right) = 0
      \\
      r(0,x) = r_o(x)
    \end{array}
    \right.
    \qquad
    \left\{
    \begin{array}{l}
      \partial_t r_* + \partial_x \left(r_* \, v_L \left(
        \Sigma \rho_* + r_*
        \right)\right) = 0
      \\
      r_*(0,x) = r_{o,*}(x) \, .
    \end{array}
    \right.
  \end{equation}
  Then, there exists a constant
  \begin{equation}
    \label{eq:35}
    \mathcal{C}_{6} \coloneqq
    \mathcal{C}
    \left(
    R_L,
    \norma{v_L}_{\W2\infty(\reali; \reali)},
    \norma{ \rho}_{\L\infty([0,T]; \W21(\Su; \reali^k))},
    \norma{\rho_*}_{\L\infty([0,T]\times \Su; \reali^k)},
    \tv(r_{o,*}), T
    \right)
  \end{equation}
  such that for all $t \in [0,T]$,
  \begin{displaymath}
    \norma{r(t) - r_*(t)}_{\L1(\Su; \reali)}
    \leq
    \norma{r_o - r_{o,*}}_{\L1(\Su; \reali)}
    +
    \mathcal{C}_{6} \,
    \norma{\rho - \rho_*}_{\C0([0,t]; \W11(\Su; \reali^k ))}
    \, t.
  \end{displaymath}
\end{lemma}

\begin{proofof}{Lemma~\ref{lem:stabilita_loc_nonloc_2}}
  Call $w(t,x,q) \coloneqq v_L\left(\Sigma \rho (t,x) + q\right)$,
  $w_*(t,x,q) \coloneqq v_L\left(\Sigma \rho_* (t,x) + q\right)$ and
  let
  \begin{displaymath}
    \mathrm{P} \coloneqq \max \{ \norma{\rho}_{\L\infty([0,T] \times \Su; \reali^k)},
    \norma{\rho_*}_{\L\infty([0,T] \times \Su; \reali^k)} \}.
  \end{displaymath}
  By~\ref{item:11} in \Cref{lem:buono}, both $w, w_*$
  satisfy~\ref{item:26}$\cdots$\ref{item:31} (with the choice
  $R \coloneqq \max\{R_L , \mathrm{P}\}$). Hence, Theorem~\ref{thm:K}
  ensures the well-posedness of problems~\eqref{eq:53}.  In
  particular, \ref{item:16} and~\ref{item:24} in Theorem~\ref{thm:K}
  yield that for all $t \in [0,T]$,
  \begin{align*}
    \norma{r(t) - r_*(t)}_{\L1(\Su; \reali)}
     & \leq
    \norma{r_o - r_{o,*}}_{\L1(\Su; \reali)}
    +
    t \,
    \mathrm{K}
    \,
    \Bigl(
    \norma{w - w_*}_{\L\infty([0,T] \times \Su; \W1\infty(\reali; \reali))}
    \\
     &
    \qquad \qquad \qquad
    +
    \sup_{t \in [0,T]} \int_{\Su} \esssup_{q \in \reali}
    \modulo{
      \partial_x w (t, x, q) - \partial_x w_*(t, x, q)
    } \d x
    \Bigr)
  \end{align*}
  where
  \[ \mathrm{K} \coloneqq \mathcal{C} \left(
    \begin{array}{c}
        \norma{w}_{\L\infty([0,T]\times \Su; \W1\infty(\reali; \reali))},
        \norma{\partial_x w}_{\L\infty([0,T]\times \Su; \W1\infty(\reali; \reali))},
        \\
        \esssup_{\tau \in [0,t]}
        \norma{\partial_{xx}^2 w(\tau)}_{\L1(\Su; \L\infty(\reali; \reali ))},
        { \norma{w_*}_{\L\infty([0,T] \times \Su ; \W1\infty( \reali; \reali))},
        }
        R,
        \tv(r_{o,*}),
        T
      \end{array}
    \right).\] Hence, by making the expressions of $w$, $w_*$, $R$
  explicit and applying~\eqref{eq:51}-\eqref{eq:52}-\eqref{eq:63}, we
  complete the proof.
\end{proofof}

The proof of the well-posedness of the mixed nonlocal-local
system~\eqref{eq:3} in the periodic setting relies on the interplay
between the well-posedness and stability results for the local and
nonlocal cases developed in the previous sections.

\begin{proofof}{\Cref{thm:WP}}
  As a first step, we assume that the initial datum meets higher
  regularity assumptions, which will be later relaxed by means of a
  bootstrap argument.

  Fix $\mathrm{M}, T > 0$ and suppose
  $\rho_o \in \C2(\Su; \reali_+^k)$ and $r_o \in \cBV{}(\Su; [0,R_L])$
  satisfy
  \begin{equation}\label{eq:67}
    \norma{\rho_o}_{\W21(\Su; \reali^k)}
    +
    \norma{r_o}_{\cBV{}(\Su; \reali)}
    \leq
    \mathrm{M} \,.
  \end{equation}
  We iteratively define the sequence $(\rho_n, r_n)$. In the case
  $n=1$, $(\rho_1,r_1)$ is the solution to
  \begin{displaymath}
    \begin{array}{c}
      \left\{
      \begin{array}{lc@{\,\,}lc}
        \partial_t \rho_1^i +
        \partial_x \left( \rho_1^i v_{NL}^i \left( \Sigma \rho_1 * \eta^i + (r_o * \eta^{i}_{k+1})(x) \right)\right) = 0
         &
        \mbox{in}
         &
        [0,T] \times \Su
         &
        \qquad i \in \{1, \ldots, k\}
        \\
        \rho_1(0,x) = \rho_o(x)
         &
        \mbox{in}
         &
        \Su
         &
      \end{array}
      \right.
      \\
      \left\{
      \begin{array}{lc@{\,\,}l}
        \partial_t r_1 +
        \partial_x \left( r_1\, v_{L} \left( \Sigma \rho_1(t,x) + r_1 \right)\right) = 0
         &
        \mbox{in}
         &
        [0,T] \times \Su
        \\
        r_1(0,x) = r_o(x)
         &
        \mbox{in}
         &
        \Su
      \end{array}
      \right.
    \end{array}
  \end{displaymath}
  and, for $n \geq 2$,
  \begin{displaymath}
    \left\{
    \begin{array}{lc@{\,\,}l@{\,}c}
      \partial_t \rho_n^i +
      \partial_x \left( \rho_n^i v_{NL}^i \left( \Sigma \rho_n * \eta^i + (r_{n-1} * \eta^i_{k+1})(t,x) \right)\right) = 0
       &
      \mbox{in}
       &
      [0,T] \times \Su
       &
      \quad i \in \{1, \ldots, k\}
      \\
      \rho_n(0,x) = \rho_o(x)
       &
      \mbox{in}
       &
      \Su
       &
    \end{array}
    \right.
  \end{displaymath}
  \begin{displaymath}
    \left\{
    \begin{array}{lc@{\,\,}l}
      \partial_t r_n +
      \partial_x \left( r_n\, v_{L} \left( \Sigma \rho_n(t,x) + r_n \right)\right) = 0
       &
      \mbox{in}
       &
      [0,T] \times \Su
      \\
      r_n(0,x) = r_o(x)
       &
      \mbox{in}
       &
      \Su.
    \end{array}
    \right.
  \end{displaymath}

  \paragraph{\underline{Step~0}: If $(\rho_n,r_n)$ are defined, then
    for all $(t,x) \in [0,T]\times\reali$, $\rho_n (t,x) \geq 0$ and
    $r_n (t,x) \in [0,R_L]$.}

  Use the Maximum Principle~\cite[Theorem~3]{Kruzkov1970}. First,
  compare $\rho_n$ with the zero solution to the nonlocal
  equation. Then, compare $r_n$ with the constant solutions attaining
  the values $0$ and $R_L$, thanks to the positivity of $\rho_n$.

  \paragraph{\underline{Step~1}: For all $n$, the couple
    $\left(\rho_n, r_n\right)$ is well-defined,
    $\rho_n \in \C0\left([0,T]; \W11(\Su; \reali^k)\right)$ and
    $r_n \in \C0\left([0,T]; \L1(\Su; \reali)\right)$.}

  We refer for the details of this proof to~\cite[Step~1 in the Proof
    of Theorem~2.3]{opera-seconda}. Here, we recall its structure
  through the diagram:
  \begin{center}
    \begin{tabular}{cc@{\hspace{-.025cm}}c@{\hspace{-.025cm}}c}
      $n=1$
       & \ref{item:17}, \ref{item:hyp-eta} and $r_o \in \L1(\Su; [0,R_L])$
       & $\implies$
       & \begin{tabular}{c}
           $V$ in \Cref{lem:buono} satisfies
           \\
           \ref{eq:87} in \Cref{prop:fine_estimates};
         \end{tabular}
      \\[10pt]
       & \begin{tabular}{c}
           $V$ satisfies~\ref{eq:87} in \Cref{prop:fine_estimates}
           \\
           and $\rho_o \in \C2(\Su; \reali_+^k)$
         \end{tabular}
       & $\implies$
       & \begin{tabular}{c}
           $\rho_1$ satisfies~\ref{item:30}
           in \Cref{prop:fine_estimates}
           \\
           and~\ref{item:12} in Theorem~\ref{thm:4}.
         \end{tabular}
      \\[10pt]
      $n\geq1$
       & \begin{tabular}{c}
           \ref{item:21}, $\rho_n^i \geq 0$ and $\rho_n$ satisfies~\ref{item:30}
           \\
           in \Cref{prop:fine_estimates}
         \end{tabular}
       & $\implies$
       & \begin{tabular}{c}
           $w$ in \Cref{lem:buono} satisfies
           \\
           \ref{item:26}$\cdots$\ref{item:31};
         \end{tabular}
      \\[10pt]
       & \begin{tabular}{c}
           $w$ satisfies
           \ref{item:26}$\cdots$\ref{item:31}
           \\
           and $r_o \in \cBV{}(\Su; [0,R_L])$
         \end{tabular}
       & $\implies$
       & \begin{tabular}{c}
           $r_n \in \C0\left([0,T]; \L1(\Su; [0, R_L])\right)$
           \\
           by~\ref{item:16} and~\ref{item:19} in Theorem~\ref{thm:K}
           \\
           and~\ref{item:25} in Theorem~\ref{thm:K};
         \end{tabular}
      \\[14pt]
       & \begin{tabular}{c}
           \ref{item:17}, \ref{item:hyp-eta},
           $\rho_o \in \C2(\Su; \reali_+^k)$
           \\
           and
           $r_n \in \C0 \left([0,T];\L1 (\Su;[0,R_L])\right)$
         \end{tabular}
       & $\implies$
       & \begin{tabular}{c}
           $\rho_{n+1} \in \C0\left([0,T]; \W11(\Su; \reali_+^k)\right)$
           \\
           satisfies~\ref{item:30} in \Cref{prop:fine_estimates}.
         \end{tabular}
    \end{tabular}
  \end{center}

  \paragraph{\underline{Step~2:}
    There exists a constant
    \begin{displaymath}
      \mathcal{C}_{7} \coloneqq \mathcal{C} (
      \norma{v_{NL}}_{\W3\infty(\reali,\reali^k)},
      \norma{\eta}_{\W3\infty(\reali,\reali^{k \times (k+1)})},
      H,
      \norma{\rho_o}_{\cBV{1}(\Su; \reali^k)},
      \norma{r_o}_{\L1(\Su; \reali)},
      T)
    \end{displaymath}
    such that
    $\norma{\rho_n(t)- \rho_{n-1}(t)}_{\W11(\Su;\reali^k)} \leq
      \mathcal{C}_{7} \, \norma{r_{n-1} - r_{n-2}} _{\C0([0,t];
        \L1(\Su;\reali))} \, t$ for all $t \in [0,T]$.  }

  By~\ref{item:16} in Theorem~\ref{thm:K}, for all $t \in [0,T]$ it
  holds that
  \begin{equation}\label{eq:64}
    \max \left\{
    \norma{r_{n-1}(t)}_{\L1(\Su; \reali)},
    \norma{r_{n-2}(t)}_{\L1(\Su; \reali)}
    \right\}
    =
    \norma{r_o}_{\L1(\Su; \reali)}\, .
  \end{equation}
  Hence, an application of \Cref{lem:stabilita_loc_nonloc} with
  $\mathrm{R} \coloneqq \norma{r_o}_{\L1(\Su; \reali)}$ proves
  \textbf{\underline{Step~2}}.

  \paragraph{\underline{Step~3:} There exists a constant
    \begin{displaymath}
      \mathcal{C}_{8} \coloneqq \mathcal{C}
      \left(
      \begin{array}{c}
          \norma{v_{NL}}_{\W3\infty(\reali,\reali^k)},
          R_L,
          \norma{v_L}_{\W2\infty(\reali; \reali)},
          \norma{\eta}_{\W3\infty(\reali;\reali^{k \times (k+1)})},
          H,
          \\
          \norma{\rho_o}_{\W21(\Su; \reali^k)},
          \norma{r_o}_{\cBV{}(\Su; \reali)},
          T
        \end{array}
      \right)
    \end{displaymath}
    such that
    $\norma{r_n(t)- r_{n-1}(t)}_{\L1(\Su;\reali)} \le \mathcal{C}_{8}
      \, \norma{r_{n-1} - r_{n-2}} _{\C0([0,t]; \L1(\Su;\reali))} \, t^2
    $ for all $t \in [0,T]$.}

  An application of \Cref{lem:stabilita_loc_nonloc_2} yields
  \begin{align*}
     & \norma{r_n(t) - r_{n-1}(t)}_{\L1(\Su; \reali)}
    \\
     & \leq
    \mathcal{C}
    \left(
    \begin{array}{@{}c@{}}
        R_L,
        \norma{v_L}_{\W2\infty(\reali; \reali)},
        \\
        \norma{\rho_{n}}_{\L\infty([0,T]; \W21(\Su; \reali^k))},
        \\
        \norma{\rho_{n-1}}_{\L\infty([0,T]\times \Su; \reali^k)},
        \\
        \tv(r_o), T
      \end{array}
    \right)
    \norma{\rho_{n} - \rho_{n-1}}_{\C0([0,t]; \W11(\Su; \reali^k ))}
    \, t
    \\
     &
    \leq
    \mathcal{C}
    \left(
    \begin{array}{@{}c@{}}
        \norma{v_{NL}}_{\W3\infty(\reali; \reali^k)},
        R_L,
        \norma{v_L}_{\W2\infty(\reali; \reali)},
        \\
        \norma{\eta}_{\W3\infty(\reali; \reali^{k\times (k+1)})},
        H,
        \\
        \norma{\rho_o}_{\W21(\Su; \reali^k)},
        \tv(r_o),
        \\
        \norma{r_{n-1}}_{\C0([0,t]; \L1(\Su; \reali))},
        \norma{r_{n-2}}_{\C0([0,t]; \L1(\Su; \reali))},
        T
      \end{array}
    \right)
    \norma{\rho_{n} - \rho_{n-1}}_{\C0([0,t]; \W11(\Su; \reali^k ))}
    \, t
  \end{align*}
  using~\ref{item:35} in \Cref{thm:4},~\ref{item:30} in
  \Cref{prop:fine_estimates} and~\eqref{eq:40} in \Cref{lem:buono}.
  Using~\eqref{eq:64} and \textbf{\underline{Step~2}} completes the
  proof to \textbf{\underline{Step~3}}.

  \paragraph{\underline{Step~4:} A priori estimates on the solutions
    to~\eqref{eq:3} in the sense of \Cref{def:sol}.  }

  Call
  $V^i(t,x,q) \coloneqq v_{NL}^i\left(\Sigma q + \left(
    r(t)*\eta^i_{k+1} \right)(x)\right)$ and
  $w(t,x,q) \coloneqq v_L\left( \Sigma \rho(t,x) + q \right)$ and
  observe that they respectively satisfy~\ref{eq:87} in
  \Cref{prop:fine_estimates} and~\ref{item:26}$\cdots$\ref{item:31}.
  In fact,
  \begin{center}
    \begin{tabular}{ccl}
      $r \in \C0\left([0,T]; \L1(\Su; \reali)\right)$
       & $\implies$
       & \ref{item:22} in \Cref{lem:buono} holds
      \\[8pt]
       & $\implies$
       &
      \begin{tabular}{l}
        $\rho$ satisfies~\ref{item:30} in \Cref{prop:fine_estimates}
        \\
        and $\rho \geq 0$ by~\ref{item:12} in Theorem~\ref{thm:4}
      \end{tabular}
      \\[12pt]
       & $\implies$
       & \ref{item:11} in \Cref{lem:buono} holds.
    \end{tabular}
  \end{center}
  A further use of the Maximum
  Principle~\cite[Theorem~3]{Kruzkov1970}, thanks to $\rho \geq 0$,
  ensures that for all $(t,x) \in [0,T] \times\reali$,
  $r (t,x) \in [0,R_L]$.

  Any solution $(\rho,r)$ to~\eqref{eq:3}, satisfies the following
  estimates as long as it exists:
  \begin{align*}
    \norma{r(t)}_{\L1(\Su; \reali)}
     & =
    \norma{r_o}_{\L1(\Su; \reali)},
    \tag*{
      [\mbox{By~\ref{item:18} in Theorem~\ref{thm:K}}]
    }
    \\
    \norma{\rho(t)}_{\L\infty(\Su; \reali^k)}
     &
    \leq
    \mathcal{C} \left(
    \begin{array}{c}
        \norma{V}_{\L\infty([0,t]; \W2\infty(\Su \times \reali^k; \reali^k))},
        \\
        \norma{[\eta_1, \ldots, \eta_k]}_{\W{2}\infty(\reali; \reali^{k \times k})},
        \\
        H,
        \norma{\rho_o}_{\L\infty(\Su; \reali^k)}, t
      \end{array}
    \right)
    \tag*{
      [\mbox{By~\ref{item:35} in \Cref{thm:4}}]
    }
    \\
     &
    \leq
    \mathcal{C} \left(
    \begin{array}{c}
        \norma{v_{NL}}_{\W3\infty(\reali; \reali^k)},
        \\
        \norma{\eta}_{\W{3}\infty(\reali; \reali^{k \times (k+1)})},
        H,
        \\
        \norma{\rho_o}_{\L\infty(\Su; \reali^k)},
        \norma{r_o}_{ \L1(\Su; \reali)}, t
      \end{array}
    \right)
    \tag*{
      [\mbox{By~\eqref{eq:40} in \Cref{lem:buono}}]
    }
    \\
    \norma{\rho(t)}_{\W21(\Su; \reali^k)}
     &
    \leq
    \mathcal{C} \left(
    \begin{array}{c}
        \norma{V}_{\L\infty([0,t]; \W3\infty(\Su \times B_{\mathrm{M}_*}; \reali^k))},
        \\
        \norma{[\eta_1, \ldots, \eta_k]}_{\W{3}\infty(\reali; \reali^{k \times k})},
        H,
        \\
        \norma{\rho_o}_{\W21(\Su; \reali^k)},
        t
      \end{array}
    \right)
    \tag*{
      [\mbox{By~\ref{item:30} in \Cref{prop:fine_estimates}}]
    }
    \\
     &
    \leq
    \underbrace{
      \mathcal{C} \left(
      \begin{array}{c}
          \norma{v_{NL}}_{\W3\infty(\reali; \reali^k)},
          \\
          \norma{\eta}_{\W{3}\infty(\reali; \reali^{k \times (k+1)})},
          H,
          \\
          \norma{\rho_o}_{\W21(\Su; \reali^k)},
          \norma{r_o}_{\L1(\Su; \reali)},
          t
        \end{array}
      \right)
    }_{\mathcal{C}_{9}}
    \tag*{
      [\mbox{By~\ref{eq:40} in \Cref{lem:buono}}]
    }
    \\
    \tv\left(r(t)\right)
     &
    \leq
    \mathcal{C} \left(
    \begin{array}{@{}c@{}}
        \norma{w}_{\L\infty([0,T]\times \Su; \W1\infty(\reali; \reali))},
        \\
        \norma{\partial_x w}_{\L\infty([0,T]\times \Su; \W1\infty(\reali; \reali))},
        \\
        \esssup_{\tau \in [0,t]}
        \norma{\partial_{xx}^2 w(\tau)}_{\L1(\Su; \L\infty(\reali; \reali ))},
        \\
        \tv(r_o), R, t
      \end{array}
    \right)
    \tag*{
      [\mbox{By~\ref{item:20} in Theorem~\ref{thm:K}}]
    }
    \\
     &
    \leq
    \mathcal{C} \left(
    \begin{array}{c}
        R_L,
        \norma{v_L}_{\W2\infty(\reali; \reali)},
        \norma{\rho}_{\L\infty([0,T]\times \Su; \reali^k)},
        \\
        \norma{\partial_x \rho}_{\L\infty([0,T]; \W11(\Su; \reali^k))},
        \tv(r_o),
        t
      \end{array}
    \right)
    \tag*{
      [By~\eqref{eq:51}--\eqref{eq:63} in \Cref{lem:buono}]}
    \\
     &
    \leq
    \underbrace{
      \mathcal{C} \left(
      \begin{array}{c}
          \norma{v_{NL}}_{\W3\infty(\reali; \reali^k)},
          R_L,
          \norma{v_L}_{\W2\infty(\reali; \reali)},
          \\
          \norma{\eta}_{\W{3}\infty(\reali; \reali^{k \times (k+1)})},
          H,
          \\
          \norma{\rho_o}_{\W21(\Su; \reali^k)},
          \norma{r_o}_{\cBV{}(\Su; \reali)},
          t
        \end{array}
      \right)
    }_{\mathcal{C}_{10}} \, .
  \end{align*}

  \paragraph{\underline{Step~5:} Existence of a solution
    $(\rho, r) \in \C0\left(\reali_+; \W11(\Su; \reali^k) \times
      \L1(\Su; \reali)\right)$ to~\eqref{eq:3}.}

  Fix an arbitrary $T > 0$. By the \emph{a priori} estimates in
  \textbf{\underline{Step~4}}, for all $t \in [0,T]$,
  \begin{equation}\label{eq:65}
    \norma{\rho(t)}_{\W21(\Su; \reali^k)} +
    \norma{r(t)}_{\cBV{}(\Su; \reali)}
    \leq
    \underbrace{\mathcal{C}_{9} + \mathcal{C}_{10} + \norma{r_o}_{\L1(\Su; \reali)}}_K \,,
  \end{equation}

  Fix $t \in [0,T]$ and estimate for $n \geq 2$
  \begin{align*}
     & \norma{r_{n+1} (t) - r_n (t)}_{\L1 (\Su; \reali)}
    \\
     & \leq
    \mathcal{C}_8  \,
    \norma{r_{n} - r_{n-1}}_{\C0([0,t];\L1 (\Su; \reali))} \, t^2
    \tag*{[\mbox{By \textbf{\underline{Step~3}}}]}
    \\
     & \leq
    \mathcal{C}\left(
    \begin{array}{@{}c@{}}
        \norma{v_{NL}}_{\W3\infty(\reali; \reali^k)},
        R_L,
        \norma{v_L}_{\W2\infty (\reali; \reali)},
        \\
        \norma{\eta}_{\W3\infty(\reali; \reali^{k\times(k+1)})},
        H, K, T
      \end{array}
    \right)
    \norma{r_n - r_{n-1}}_{\C0 ([0,t]; \L1 (\Su; \reali))}
    \, t^2
    \tag*{[\mbox{By~\eqref{eq:65} }]}
  \end{align*}
  For a suitable choice of $T^*$, $r_n$ is a Cauchy sequence in
  $\C0 \left([0,t]; \L1 (\Su; \reali^k)\right)$.  Then,
  \textbf{\underline{Step~2}} ensures that also $\rho_n$ is a Cauchy
  sequence in $\C0 \left([0,t]; \W11 (\Su; \reali)\right)$.

  Hence, $\left(\rho_n, r_n\right)$ converges to a limit point
  $(\rho, r) \in \C0\left([0,T^*]; \W11(\Su; \reali^k) \times \L1(\Su;
    \reali)\right)$, which is a solution of~\eqref{eq:3} in the sense
  of \Cref{def:sol} by \Cref{lem:GoodSol}. Moreover, by
  \textbf{\underline{Step~0}}, for all $(t,x) \in [0,T^*] \times \Su$,
  $(\rho,r) (t,x) \in \reali_+^k \times [0,R_L]$.

  By \Cref{prop:fine_estimates} and~\ref{item:20} in
  Theorem~\ref{thm:K},
  $\left( \rho(T^*), r(T^*) \right) \in \C2(\Su; \reali^k) \times
    \cBV{}(\Su; \reali)$ and~\eqref{eq:65} holds.  Hence, we can iterate
  the same argument in order to extend the solution to the time
  intervals $[T^*, 2T^*], [2T^*, 3T^*], \ldots $, covering the whole
  time interval $[0,T]$.  By the arbitrariness of $T$, this completes
  the proof of \textbf{\underline{Step~5}}.

  \paragraph{\underline{Step~6:} Locally Lipschitz dependence on the
    initial datum. }

  Fix $T>0$ and
  $(\rho_{o,*}, r_{o,*}) \in \left(\C2(\Su; \reali^k_+) \times
    \cBV{}(\Su; [0, R_L])\right)$.  Consider the problems
  \begin{displaymath}
    \begin{cases}
      \partial_t \rho^i
      {+}
      \partial_x \!\left(\!\rho^i \,
      v_{NL}^i\left( \Sigma \rho*\eta^i {+} r*\eta^i_{k+1}
      \right)\right)
      {=}0
      \\
      \partial_t r
        {+} \partial_x \left(r \, v_L\left(\Sigma \rho{+}r \right)\right)
      {=} 0
      \\
      \rho(0) = \rho_o
      \\
      r(0) = r_o
    \end{cases}
    \begin{cases}
      \partial_t \rho^i_*
      {+}
      \partial_x \!\left(\!\rho^i_* \,
      v_{NL}^i\left( \Sigma \rho_**\eta^i {+} r_**\eta^i_{k+1}
      \right)\right)
      {=}0
      \\
      \partial_t r_*
      + \partial_x \left(r_* \;
      v_L\left(\Sigma \rho_*+ r_* \right)\right)
      {=} 0
      \\
      \rho_*(0) = \rho_{o,*}
      \\
      r_*(0) = r_{o,*} \,.
    \end{cases}
  \end{displaymath}
  Call
  $\mathrm{R} \coloneqq \max \{ \norma{r_o}_{\L1(\Su; \reali)},
    \norma{r_{o,*}}_{\L1(\Su; \reali)} \}$.  Then, for all
  $t \in [0,T]$, using \Cref{lem:stabilita_loc_nonloc_2}, the \emph{a
    priori} estimates in \textbf{\underline{Step~4}} and
  \Cref{lem:stabilita_loc_nonloc}, we have
  \begin{align*}
     & \norma{r(t) - r_*(t)}_{\L1(\Su; \reali)}
    \\
     & \leq
    \norma{r_o - r_{o,*}}_{\L1(\Su; \reali)}
    +
    \mathcal{C}
    \left(
    \begin{array}{c}
        R_L,
        \norma{v_L}_{\W2\infty(\reali; \reali)},
        \\
        \norma{\rho}_{\L\infty([0,T]; \W21(\Su; \reali^k))},
        \\
        \norma{\rho_*}_{\L\infty([0,T]\times \Su; \reali^k)},
        \\
        \tv(r_{o,*}), T
      \end{array}
    \right)
    \norma{\rho - \rho_*}_{\C0([0,t]; \W11(\Su; \reali^k ))}
    \, t
    \\
     & \leq
    \norma{r_o - r_{o,*}}_{\L1(\Su; \reali)}
    \\
     & \quad
    +
    \underbrace{
      \mathcal{C}
      \left(
      \begin{array}{c}
          \norma{v_{NL}}_{\W3\infty(\reali; \reali^k)},
          R_L,
          \norma{v_L}_{\W2\infty(\reali; \reali)},
          \\
          \norma{\eta}_{\W{3}\infty(\reali; \reali^{k \times (k+1)})},
          H,
          \\
          \norma{\rho_o}_{\W21(\Su; \reali^k)},
          \norma{r_o}_{\L1(\Su; \reali)},
          \\
          \norma{\rho_{o,*}}_{\L\infty(\Su; \reali^k)},
          \norma{r_{o,*}}_{\cBV{}(\Su; \reali)},
          T
        \end{array}
      \right)
    }_{\mathcal{C}_{11}}
    \norma{\rho - \rho_*}_{\C0([0,t]; \W11(\Su; \reali^k ))}
    \, t
    \\
     & \leq
    \norma{r_o - r_{o,*}}_{\L1(\Su; \reali)}
    \\
     & \quad
    +
    \mathcal{C}_{11} \, t
    \left(
    \left(1 + \mathcal{C}_{5} t\right)
    \norma{\rho_o -  \rho_{o,*}}_{\W11(\Su; \reali^k)}
    +
    \mathcal{C}_{5} \int_0^t \norma{r(\tau)- r_* (\tau)}_{\L1(\Su; \reali)}
    \dd \tau
    \right)
  \end{align*}

  So, one obtains that there exists a constant
  \begin{displaymath}
    \mathcal{C}_{12} \coloneqq
    \mathcal{C} \left(
    \begin{array}{c}
        \norma{v_{NL}}_{\W3\infty(\reali; \reali^k)},
        R_L,
        \norma{v_L}_{\W2\infty(\reali; \reali)},
        \norma{\eta}_{\W{3}\infty(\reali; \reali^{k \times (k+1)})},
        H,
        \\
        \norma{\rho_o}_{\W21(\Su; \reali^k)},
        \norma{r_o}_{\L1(\Su; \reali)},
        \norma{\rho_{o,*}}_{\cBV{}(\Su; \reali^k)},
        \norma{r_{o,*}}_{\cBV{}(\Su; \reali)},
        T
      \end{array}
    \right)
  \end{displaymath}
  such that
  \begin{eqnarray*}
    &
    & \norma{r(t) - r_*(t)}_{\L1(\Su; \reali)}
    \\
    & \leq
    & \norma{r_o - r_{o,*}}_{\L1(\Su; \reali)}
    +
    \mathcal{C}_{12} \, t \left(
    \norma{\rho_o - \rho_{o,*}}_{\W11(\Su; \reali^k)}
    +
    \int_0^t \norma{r(\tau)- r_*(\tau)}_{\L1(\Su; \reali)} \d\tau
    \right).
  \end{eqnarray*}
  An application of Gronwall Lemma concludes that for all
  $t \in [0,T]$
  \begin{displaymath}
    \norma{r(t)-r_*(t)}_{\L1(\Su; \reali)}
    \leq
    \left(
    \norma{r_o - r_{o,*}}_{\L1(\Su; \reali)}
    +
    \mathcal{C}_{12} \, t \, \norma{\rho_o - \rho_{o,*}}_{\W11(\Su; \reali^k)}
    \right) e^{ \mathcal{C}_{12} t^2}\, .
  \end{displaymath}
  By \Cref{lem:stabilita_loc_nonloc}, we get
  \begin{eqnarray*}
    \norma{\rho(t) -\rho_*(t)}_{\W11(\Su; \reali^k)}
    & \leq
    & (1 + \mathcal{C}_{5} t)
    \norma{\rho_o - \rho_{o,*}}_{\W11(\Su; \reali^k)}
    \\
    &
    & +
    \mathcal{C}_{5}
    \left(
    \norma{r_o - r_{o,*}}_{\L1(\Su; \reali)} t
    +
    \mathcal{C}_{12} t^2
    \norma{\rho_o - \rho_{o,*}}_{\W11(\Su; \reali^k)}
    \right)
    e^{ \mathcal{C}_{12} t^2} \,.
  \end{eqnarray*}
  In conclusion, for all $t \in [0,T]$ there exists a constant
  \begin{equation}
    \label{eq:78}
    \mathcal{C}_{13}
    \coloneqq
    \mathcal{C} \left(
    \begin{array}{c}
        \norma{v_{NL}}_{\W3\infty(\reali; \reali^k)},
        R_L,
        \norma{v_L}_{\W2\infty(\reali; \reali)},
        \norma{\eta}_{\W{3}\infty(\reali; \reali^{k \times (k+1)})},
        H,
        \\
        \norma{\rho_o}_{\W21(\Su; \reali^k)},
        \norma{r_o}_{\L1(\Su; \reali)},
        \norma{\rho_{o,*}}_{\cBV{}(\Su; \reali^k)},
        \norma{r_{o,*}}_{\cBV{}(\Su; \reali)},
        T
      \end{array}
    \right)
  \end{equation}
  such that
  \begin{eqnarray}
    &
    & \norma{\rho(t)-  \rho_*(t)}_{\W11(\Su; \reali^k)}
    +
    \norma{r(t)- r_*(t)}_{\L1(\Su; \reali)}
    \\
    \nonumber
    & \leq
    &
    (1+\mathcal{C}_{13} \, t)
    \left(
    \norma{\rho_o - \rho_{o,*}}_{\W11(\Su; \reali^k)}
    +
    \norma{r_o - r_{o,*}}_{\L1(\Su; \reali)}
    \right).
  \end{eqnarray}

  \paragraph{\underline{Step~7:} The solution is Locally Lipschitz
    continuous in time. }

  Fix $T > 0$. Call $V$, $w$ and $R$ as in \Cref{lem:buono} and
  apply~\ref{item:36} in \Cref{prop:fine_estimates} and~\ref{item:25}
  in Theorem~\ref{thm:K} to obtain that for all $t, t_* \in [0,T]$,
  \begin{flalign}
    \nonumber
     & \norma{\rho (t) - \rho (t_*)}_{\W11 (\Su; \reali^k)}
    \\
    \nonumber
     & \leq
    \mathcal{C}\left(
    \begin{array}{c}
        \norma{V}_{\L\infty([0,T]; \W{3}\infty(\Su \times \reali^k; \reali^k))},
        \\
        \norma{[\eta_1, \ldots, \eta_k]}_{\W{3}\infty(\reali; \reali^{k \times k})},
        H,
        \norma{\rho_o}_{\cBV{1}(\Su; \reali^k)},
        T
      \end{array}
    \right)
    \modulo{t-t_*}
     & [\mbox{Use~\eqref{eq:40}}]
    \\
    \label{eq:76}
     & \leq
    \mathcal{C}
    \left(
    \begin{array}{c}
        \norma{v_{NL}}_{\W3\infty(\reali; \reali^k)},
        \norma{\eta}_{\W{3}\infty(\reali; \reali^{k \times (k+1)})},
        H,
        \\
        \norma{\rho_o}_{\cBV{1}(\Su; \reali^k)},
        \norma{r_o}_{\L1(\Su; \reali)},
        T
      \end{array}
    \right)
    \modulo{t-t_*} \,.
    \\
    \nonumber
     & \norma{r (t) - r (t_*)}_{\L1 (\Su;\reali)}
    \\
    \nonumber
     & \leq
    \mathcal{C}
    \left(
    \begin{array}{c}
        \norma{w}_{\L\infty([0,T]\times \Su; \W1\infty(\reali; \reali ))},
        \\
        \norma{\partial_x w}_{\L\infty([0,T]\times \Su; \W1\infty(\reali; \reali ))},
        \\
        \esssup_{\tau \in [0,T]}
        \norma{\partial_{xx}^2 w(\tau)}_{\L1(\Su; \L\infty(\reali; \reali ))},
        \\
        R, \tv(r_o), T
      \end{array}
    \right)
    \modulo{t-t_*}
     &
    \\
    \nonumber
     & \leq
    \mathcal{C}
    \left(
    \begin{array}{c}
        R_L,
        \norma{v_L}_{\W2\infty(\reali; \reali)},
        \\
        \norma{ \rho}_{\L\infty([0,T]\times \Su; \reali^k)},
        \norma{\partial_x \rho}_{\L\infty([0,T]; \W11(\Su; \reali^k))}
        \\
        \tv(r_o), T
      \end{array}
    \right)
    \modulo{t-t_*}
     & [\mbox{By~\ref{item:11} in \Cref{lem:buono}}]
    \\
    \nonumber
     & \leq
    \mathcal{C}
    \left(
    \begin{array}{c}
        \norma{v_{NL}}_{\W3\infty(\reali; \reali^k)},
        R_L,
        \norma{v_L}_{\W2\infty(\reali; \reali)},
        \\
        \norma{\eta}_{\W{3}\infty(\reali; \reali^{k \times (k+1)})},
        H,
        \\
        \norma{\rho_o}_{\W21(\Su; \reali^k)},
        \norma{r_o}_{\cBV{}(\Su; \reali)},
        T
      \end{array}
    \right)
    \modulo{t-t_*}
     & [\mbox{By \textbf{\underline{Step~4}}}],
  \end{flalign}
  completing the proof of \textbf{\underline{Step~7}}.

  \paragraph{\underline{Step~8:} Uniqueness of the solution
    to~\eqref{eq:3} in the sense of \Cref{def:sol}. } Fix
  $(\rho_o, r_o) \in \C2(\Su;\reali^k_+) \times \cBV{} (\Su;[0,
      R_L])$. Assume there exist $(\rho, r), (\rho_*, r_*)$ solutions
  to~\eqref{eq:3} in the sense of \Cref{def:sol}, so that
  $(\rho, r), (\rho_*, r_*) \in \C0\left([0,T]; \W11(\Su;
    \reali^k)\times \L1(\Su; \reali)\right)$.
  Definition~\ref{def:sol} allows using
  \Cref{lem:stabilita_loc_nonloc} with
  $\mathrm{R} = \norma{r_o}_{\L1(\Su; \reali)}$ to obtain:
  \begin{equation}
    \label{eq:83}
    \norma{\rho (t) - \rho_* (t)}_{\W11 (\Su; \reali^k)}
    \leq
    \mathcal{C}_{5} \int_0^t
    \norma{r (\tau) - r_* (\tau)}_{\L1 (\Su;\reali)} \d\tau \,.
  \end{equation}
  By~\ref{item:30} in Proposition~\ref{prop:fine_estimates}
  and~\ref{item:22} in Lemma~\ref{lem:buono}, $\rho$ and $\rho_*$ are
  sufficiently regular to apply \Cref{lem:stabilita_loc_nonloc_2}
  which yields
  \begin{align*}
     & \norma{r(t) - r_*(t)}_{\L1(\Su; \reali)}
    \\
     & \leq
    \mathcal{C} \left(
    \begin{array}{c}
        R_L,
        \norma{v_L}_{\W2\infty(\reali; \reali)},
        \norma{\rho_*}_{\L\infty([0,T]; \W21(\Su; \reali^k))},
        \\
        \norma{\rho}_{\L\infty([0,T]\times \Su; \reali^k)},
        \tv(r_{o}), T
      \end{array}
    \right)
    t \, \norma{\rho- \rho_*}_{\C0([0,t]; \W11(\Su; \reali^k))}
    \\
     & \leq
    \mathcal{C} \left(
    \begin{array}{c}
        \norma{v_{NL}}_{\W3\infty(\reali; \reali^k)},
        R_L,
        \norma{v_L}_{\W2\infty(\reali; \reali)},
        \\
        \norma{\eta}_{\W{3}\infty(\reali; \reali^{k \times (k+1)})},
        H,
        \\
        \norma{\rho_o}_{\W{2}{1}(\Su; \reali^k)},
        \norma{r_o}_{\cBV{}(\Su; \reali)},
        T
      \end{array}
    \right)
    t \, \norma{\rho - \rho_*}_{\C0([0,t]; \W11(\Su; \reali^k))}
    \tag*{[\mbox{By \textbf{\underline{Step~4}}}]}
  \end{align*}
  Complete the previous expression with~\eqref{eq:83} and apply the
  Gronwall Lemma to obtain
  $\norma{r(t) - r_*(t)}_{\L1(\Su; \reali)} = 0$ and, returning
  to~\eqref{eq:83}, also
  $\norma{\rho (t) - \rho_* (t)}_{\W11 (\Su; \reali^k)} = 0$.

  \paragraph{\underline{Step~9:} Definition and properties of
    $\mathcal{S}$. }

  The above steps allow defining the map
  \begin{displaymath}
    \mathcal{S}
    \colon
    \reali_+
    \times \C2(\Su;\reali^k_+)
    \times \cBV{} (\Su; [0, R_L])
    \longrightarrow
    \C2(\Su;\reali^k_+)
    \times \cBV{} (\Su; [0, R_L])
  \end{displaymath}
  by
  $\mathcal{S}_t (\rho_o,r_o) \coloneqq \left(\rho (t), r (t)\right)$,
  where $(\rho,r)$ is the solution constructed above.  By
  Definition~\ref{def:sol} and by the uniqueness of solutions, it
  immediately follows that $\mathcal{S}$ is a semigroup.

  Call $B_\mathrm{M}$ the set of pairs $(\rho_o, r_o)$
  satisfying~\eqref{eq:67}.  The restriction of $\mathcal{S}$ to
  $[0,T] \times B_\mathrm{M}$ is Lipschitz continuous with respect to
  the modulus in $t$ and the norm
  $\W11(\Su; \reali^k) \times \L1(\Su; \reali)$ in $(\rho_o,r_o)$ with
  a Lipschitz constant that, by~\textbf{\underline{Step~6}}
  and~\textbf{\underline{Step~7}}, depends only on the kernel, the
  speed laws, the final time and $\mathrm{M}$.  Hence, by uniform
  continuity, $\mathcal{S}$ can be extended to
  $[0,T] \times \overline{B}_{\mathrm{M}}$, where
  $\overline{B}_{\mathrm{M}}$ is the closure of $B_M$ with respect to
  the norm $\W11(\Su; \reali^k) \times \L1(\Su; \reali)$.  By an
  adaptation of~\cite[Lemma~3.2]{opera-seconda}, we get that
  \begin{equation}
    \label{eq:80}
    \overline{B}_{\mathrm{M}} =
    \left\{
    (\rho, r) \in \cBV1 (\Su;\reali_+^k) \times \cBV{} (\Su; [0,R_L])
    \colon
    \norma{\rho}_{\cBV1(\Su; \reali^k)}
    +
    \norma{r}_{\cBV{} (\Su;\reali)} \leq \mathrm{M}
    \right\}\,.
  \end{equation}
  The arbitrariness of $T$ and $\mathrm{M}$ allow to further extend
  $\mathcal{S}$ as in~\eqref{eq:79}.

  To prove the properties of $\mathcal{S}$, fix $(\rho_o,r_o)$ in the
  domain of $\mathcal{S}$ and choose a sequence
  $(\rho_{o,n}, r_{o,n})$ in
  $\C2(\Su;\reali^k_+) \times \cBV{} (\Su;[0, R_L])$ converging to
  $(\rho_o,r_o)$ in $\W11(\Su; \reali^k) \times \L1(\Su; \reali)$. For
  all $t_1,t_2 \in \reali_+$ we have
  $(\mathcal{S}_{t_1} \circ \mathcal{S}_{t_2}) (\rho_{o,n},r_{o,n}) =
    \mathcal{S}_{t_1+t_2} (\rho_{o,n},r_{o,n})$ and passing to the limit
  $n \to +\infty$ we prove~\ref{item:29} using \Cref{lem:GoodSol}.

  To prove the existence result stated in~\ref{item:23} apply
  \Cref{lem:GoodSol} with $\epsilon_n = e_n = 0$.
  Properties~\ref{item:40}, \ref{item:39} and~\ref{item:3} follow from
  \textbf{\underline{Step~4}}, \textbf{\underline{Step~6}},
  \textbf{\underline{Step~7}}, \textbf{\underline{Step~8}} and the
  very construction of $\mathcal{S}$ via the bootstrap argument above.
  The proof is completed.
\end{proofof}

\appendix

\section{Complements to the Classical {Kru\v zkov}
  Result} \label{sec:compl-class-kruv}

\begin{proposition}[Localized coinciding fluxes]
  \label{prop:claudia}
  Fix $T > 0$. Let
  $f, f_* \in \C1([0,T]\times \reali \times \reali ; \reali)$ with
  $\partial_{tr}^2 f, \partial_{xr}^2 f, \partial_{tr}^2 f_*,
    \partial_{xr}^2 f_* \in \Lloc\infty([0,T]\times \reali \times
    \reali; \reali)$. Choose
  $r_o , r_{o,*} \in \L\infty(\reali; \reali)$. Assume that the
  problems
  \begin{displaymath}
    \left\{
    \begin{array}{l}
      \partial_t r + \partial_x f (t,x,r) =0
      \\
      r (0,x) = r_o (x)
    \end{array}
    \right.
    \qquad \mbox{ and } \qquad
    \left\{
    \begin{array}{l}
      \partial_t r_*
      +
      \partial_x f_* (t,x,r_*) =0
      \\
      r_* (0,x) = r_{o,*} (x)
    \end{array}
    \right.
  \end{displaymath}
  respectively admit the Kru\v zkov solutions
  $r, r_* \in \L\infty([0,T] \times \reali; \reali)$. Let
  $\mathrm{N}>0$ and define
  \begin{eqnarray}
    \nonumber
    R
    & \coloneqq
    & \max\{\norma{r}_{\L\infty([0,T] \times \reali; \reali)},
    \norma{r_*}_{\L\infty([0,T] \times \reali; \reali)}\}
    \\
    \label{eq:48}
    \sigma
    & \coloneqq
    & \max \left\{
    \norma{\partial_r f}_{\L\infty ([0,T]\times B_{\mathrm{N}}\times B_{R}; \reali)}
    \,,\;
    \norma{\partial_r f_*}_{\L\infty ([0,T]\times B_{\mathrm{N}}\times B_{R}; \reali)}
    \right\}
    \\
    \nonumber
    \mathcal{K}
    & \coloneqq
    & \left\{
    (t,x) \in [0,T] \times \reali \colon
    \modulo{x} \leq \mathrm{N} - \sigma \, t
    \mbox{ and } \sigma \,t \leq \mathrm{N}
    \right\}.
  \end{eqnarray}
  Suppose that
  \begin{equation}
    \label{eq:46}
    \forall\, (t,x,r) \in \mathcal{K} \times \reali
    \qquad
    f (t,x,r) = f_* (t,x,r) \,.
  \end{equation}
  then for almost all
  $t \in [0, \min\{T, \mathrm{N}/\sigma \}]$,
  \begin{displaymath}
    \norma{r(t) - r_*(t)}_{\L1(B_{\mathrm{N} - \sigma t}; \reali)}
    \leq
    \norma{r_o - r_{o,*}}_{\L1(B_{\mathrm{N}}; \reali)} \,.
  \end{displaymath}
\end{proposition}

\begin{proofof}{Proposition~\ref{prop:claudia}}
  By~\cite[Definition~1]{Kruzkov1970}, we have that for all
  $\kappa \in \reali$ and for all
  $\phi \in \Cc1 (\mathopen]0,T\mathclose[ \times \reali; \reali_+)$
  with $\spt \phi \subseteq \mathcal{K}$,
  \begin{eqnarray*}
    \iint_{\mathcal{K}}
    \modulo{r(\tau,x) - \kappa} \, \partial_\tau \phi(\tau,x)
    \d{x}\d{\tau}
    \\
    +
    \iint_{\mathcal{K}}
    \sgn(r(\tau,x) - \kappa) \,
    \left[
      f \left(\tau,x,r (\tau,x)\right)
      - f (\tau,x,\kappa)
      \right] \partial_x \phi(\tau,x) \d{x}\d{\tau}
    \\
    -
    \iint_{\mathcal{K}}
    \sgn(r(\tau,x) - \kappa) \,
    \partial_x  f (\tau,x,\kappa) \, \phi(\tau,x) \d{x}\d{\tau}
    & \geq
    & 0 \,;
    \\
    \iint_{\mathcal{K}}
    \modulo{r_*(s,y) - \kappa} \, \partial_s \phi(s,y)
    \d{y} \d{s}
    \\
    +
    \iint_{\mathcal{K}}
    \sgn(r_*(s,y) - \kappa) \,
    \left[
      f \left(s,y,r_*(s,y)\right)
      - f (s,y,\kappa)
      \right] \partial_y \phi(s,y) \d{y} \d{s}
    \\
    -
    \iint_{\mathcal{K}}
    \sgn(r_*(s,y) - \kappa) \,
    \partial_x  f (s,y,\kappa) \, \phi(s,y) \d{y} \d{s}
    & \geq
    & 0 \,.
  \end{eqnarray*}
  by~\eqref{eq:46}. Consider now
  $\psi \in \Cc1 (\mathopen]0,T\mathclose[ \times \reali \times
      \mathopen]0,T\mathclose[ \times \reali ; \reali)$ with
  $\spt \psi \subset \mathcal{K}^2$, say $\psi = \psi(\tau, x, s, y)$.
  An application of the doubling of variables
  method~\cite[Theorem~2]{Kruzkov1970} yields that
  \begin{align*}
    \iint_{\mathcal{K}} \iint_{\mathcal{K}}
     & \Big\{
    \modulo{
      r(\tau,x) - r_*(s,y)} \, \partial_\tau \psi(\tau,x,s,y)
     &
    \\
     & +
    \sgn\left(r(\tau,x) - r_*(s,y)\right) \,
    \left[
      f \left(\tau,x,r (\tau,x)\right)
      - f \left(\tau,x,r_* (s,y)\right)
      \right]
    \partial_x \psi(\tau,x,s,y)
     &
    \\
     & -
    \sgn\left(r(\tau,x) - r_*(s,y)\right) \,
    \, \partial_x f \left(\tau,x,r_*(s,y)\right)
    \, \psi(\tau,x,s,y)
    \Big\}
    \d{x}\d{\tau}\d{y}\d{s}
     & \geq
    0 \,;
    \\
    \iint_{\mathcal{K}} \iint_{\mathcal{K}}
     & \Big\{
    \modulo{r_*(s,y) - r(\tau,x)} \, \partial_s \psi(\tau, x, s,y)
    \\
     & +
    \sgn\left(r_*(s,y) - r(\tau,x)\right) \,
    \left[
      f \left(s,y,r_* (s,y)\right)
      -  f \left(s,y,r(\tau,x)\right)
      \right] \partial_y \psi(\tau,x,s,y)
    \\
     & -
    \sgn\left(r_*(s,y) - r(\tau,x)\right) \,
    \partial_y f \left(s,y,r(\tau,x)\right) \, \psi(\tau,x,s,y)
    \Big\}
    \d{y}\d{s}\d{x}\d{\tau}
     & \geq
    0 \,.
  \end{align*}
  Adding the two expressions above, we obtain an inequality analogous
  to~\cite[eqn.~(3.4)]{Kruzkov1970},
  \begin{align}
    \label{eq:18}
    \iint_{\mathcal{K}}
     & \iint_{\mathcal{K}}
    \Big\{
    \modulo{r (\tau,x) - r_* (s,y)}
    \left(
    \partial_\tau \psi (\tau,x,s,y)+ \partial_s \psi (\tau,x,s,y)
    \right)
    \\
    \nonumber
     & +
    \sgn\left( r(\tau, x) - r_*(s,y) \right)
    \left[
      f\left(\tau,x,r(\tau,x) \right)
      -
      f\left(s,y,r_*(s,y) \right)
      \right]
    \\
     &
    \quad
    \times
    \left(
    \partial_x \psi(\tau,x,s,y) + \partial_y \psi(\tau,x,s,y)
    \right)
    \\
    \nonumber
     & +
    \sgn\left( r(\tau, x) - r_*(s,y) \right)
    \Bigl(
    \left[
      f\left(s,y,r_*(s,y) \right)
      -
      f\left(\tau,x,r_*(s,y) \right)
      \right]
    \partial_x \psi(\tau,x,s,y)
    \\
    \nonumber
     & \quad -
    \partial_x f\left( \tau,x, r_*(s,y) \right)
    \psi (\tau,x,s,y)
    +
    \left[
      f\left( s,y, r(\tau,x) \right)
      -
      f\left( \tau,x, r(\tau,x) \right)
      \right]
    \partial_y \psi(\tau,x,s,y)
    \\
    \nonumber
     & \quad +
    \partial_y f\left(s,y, r(\tau,x) \right)
    \psi(\tau,x,s,y)
    \Bigr)
    \Big\}
    \d{x}\d{\tau}\d{y}\d{s}
    \qquad
    \qquad
    \qquad
    \qquad
    \qquad
    \qquad
    \qquad
    \quad
    \geq 0 \,.
  \end{align}
  Let $\delta \in \Cc\infty (\reali; \reali_+)$ be such that
  $\spt \delta \subseteq [-1,1]$ and
  $\int_\reali \delta(\xi) \d{\xi} = 1$.  For any $h >0$ we define
  $\delta_h(x) \coloneqq \frac{1}{h} \delta \left(\frac{x}{h}\right)$
  and
  $\alpha_h(x) \coloneqq \int_{-\infty}^{x} \delta_h(\xi) \d{\xi}$.

  Choose an arbitrary
  $\Phi \in \Cc1 (\mathopen]0,T\mathclose[ \times \reali; \reali_+)$
  such that $\spt \Phi \subset \pint{\mathcal{K}}$. Call
  \begin{displaymath}
    h_*
    \coloneqq
    \min\left\{
    \norma{(\tau,x) - (s,y)}
    \colon
    (\tau,x) \in \spt\Phi \mbox{ and }
    (s,y) \in \partial\mathcal{K}
    \right\} / 2
  \end{displaymath}
  For $h \in \mathopen]0,h_*\mathclose[$ define
  \begin{displaymath}
    \psi_h (\tau, x,s,y)
    \coloneqq
    \Phi \left(
    \frac{\tau+s}{2}, \frac{x+y}{2}
    \right)
    \delta_h\left(
    \frac{\tau-s}{2}
    \right)
    \delta_h\left(
    \frac{x-y}{2}
    \right) \,.
  \end{displaymath}
  Observe that
  \begin{eqnarray*}
    \psi_h (\tau,x,s,y) \not= 0
    & \implies
    & \modulo{\tau-s}/2 \leq h
    \mbox{ and }
    \modulo{x-y}/2 \leq h
    \mbox{ and }
    \left( (\tau+s)/ 2, (x+y)/2  \right)\in \spt\Phi
    \\
    & \implies
    & \modulo{\tau-s} \leq 2h
    \mbox{ and }
    \modulo{x-y} \leq 2h
    \mbox{ and }
    \left( (\tau+s)/ 2, (x+y)/2  \right)\in \spt\Phi
    \\
    & \implies
    & d \left( (\tau,x) , \spt\Phi\right) \leq \sqrt2 h
    \mbox{ and }
    d \left( (s,y) , \spt\Phi\right) \leq \sqrt2 h
    \\
    & \implies
    & (\tau,x) \in \pint{\mathcal{K}}
    \mbox{ and }
    (s,y) \in \pint{\mathcal{K}} \,,
  \end{eqnarray*}
  proving that
  $\spt\psi_h \subseteq \pint{\mathcal{K}} \times
    \pint{\mathcal{K}}$. Compute the limit as $h \to 0$ of~\eqref{eq:18}
  following~\cite[eqn.~(3.12)]{Kruzkov1970} to
  obtain
  \begin{equation}
    \label{eq:28}
    \iint_{\mathcal{K}}
    \Big\{
    \modulo{r - r_*} \, \partial_t\Phi
    + \sgn(r - r_*)
    \left[f(\tau,x,r) - f(\tau,x,r_*)\right]
    \partial_x \Phi
    \Big\}
    \d{x}\d{\tau}
    \geq 0 \, .
  \end{equation}
  We now closely follow~\cite[Proof of Theorem~1]{Kruzkov1970}.
  Recall that by the definition of
  solutions~\cite[Definition~1]{Kruzkov1970}, there exists a zero
  measure set $\mathcal E_r \subset [0,T]$ (respectively,
  $\mathcal E_{r_*}$) such that $r(\tau,x)$ (respectively,
  $r_*(\tau, x)$) is defined for a.e.~$x$ in $\reali$ for all
  $\tau \in [0,T] \setminus \mathcal E_r$ (respectively,
  $\tau \in [0,T] \setminus  \mathcal E_{r_*}$).  Let
  $\mathcal E_\mu \subset [0, \min\{T, \mathrm{N}/\sigma\}]$ be the
  zero measure set of times which are not Lebesgue points for the map
  \begin{displaymath}
    \mu (\tau) \coloneqq \int_{B_{\mathrm{N} - \sigma \tau}} \modulo{r(\tau,x) - r_*(\tau,x)} \d{x}.
  \end{displaymath}
  Then, the set
  $\mathcal E \coloneqq \mathcal E_r \cup \mathcal E_{r_*} \cup
    \mathcal E_\mu$ is also of zero measure.  Set
  $T_* \coloneqq \min\{T, \mathrm{N}/\sigma \}$, fix
  $\tau_1, \tau_2 \in [0,T] \setminus \mathcal E $ such that
  $0 < \tau_1 < \tau_2 < T_*$,
  $0 < h \leq \min\{ \tau_1, T_* - \tau_2\}$ and
  $0 <\varepsilon < \frac{1}{2}(\mathrm{N} - \sigma \tau_2)$, consider
  as test function the map
  $\Phi_{h,\varepsilon} \in \Cc1
    (\mathopen]0,T\mathclose[\times\reali; \reali_+)$ given by
  \begin{align}
    \label{eq:49}
    \Phi_{h,\varepsilon}(\tau, x)
    \coloneqq
     &
    \left[
      \alpha_h(\tau- \tau_1) - \alpha_h(\tau - \tau_2)
      \right] \chi_\varepsilon(\tau,x) \,;
    \\
    \nonumber
    \chi_\varepsilon (\tau,x)
    \coloneqq
     & \,
    1 - \alpha_\varepsilon (\modulo{x}
    - \mathrm{N}
    + \sigma \tau + \varepsilon) \, .
  \end{align}
  Observe that, due to the choices of $h$ and $\varepsilon$, the
  support of $\Phi_{h,\varepsilon}$ is contained in
  $\pint{\mathcal{K}}$.  Moreover, for all $(\tau,x) \in \mathcal{K}$,
  it holds that
  \begin{displaymath}
    0 = \partial_{\tau} \chi_\varepsilon (\tau,x)
    + \sigma \,
    \modulo{\partial_x \chi_\varepsilon (\tau,x)}
    \geq
    \partial_{\tau} \chi_\varepsilon (\tau,x) +
    \frac{
      f(\tau,x,r)
      -
      f_*(\tau,x,r_*)
    }{r -r_* } \, \partial_x \chi_\varepsilon (\tau,x),
  \end{displaymath}
  by~\eqref{eq:48}.  Hence, it holds that
  \begin{align*}
    0
    \geq
    \left[
      \partial_{\tau} \chi_\varepsilon(\tau,x)  +
      \frac{
        f(\tau,x,r)
        -
        f_*(\tau,x,r_*)
      }{r(\tau,x) -r_*(\tau,x) } \, \partial_x \chi_\varepsilon(\tau,x)
      \right]
    \,
    \left[
      \alpha_h(\tau- \tau_1) - \alpha_h(\tau - \tau_2)
      \right] \,;
    \\
    \left[
      \delta_h(\tau- \tau_1) - \delta_h(\tau - \tau_2)
      \right] \chi_\varepsilon(\tau,x)
    \geq
    \partial_t \Phi_{h,\varepsilon}(\tau,x) +
    \frac{
      f(\tau,x,r)
      -
      f_*(\tau,x,r_*)
    }{r -r_* } \,
    \partial_x \Phi_{h,\varepsilon}(\tau,x)
  \end{align*}
  where we used~\eqref{eq:49}. Then,
  \begin{align}
    \label{eq:45}
     & \modulo{r(\tau,x) - r_*(\tau,x)}
    \left[
      \delta_h(\tau- \tau_1) - \delta_h(\tau - \tau_2)
      \right] \chi_\varepsilon (\tau,x)
    \\
    \nonumber
     & \geq
    \modulo{r(\tau,x) - r_*(\tau,x)} \, \partial_t \Phi_{h,\varepsilon}(\tau,x)
    \\
    \nonumber
     & \qquad +
    \sgn
    \left(
    r (\tau,x)- r_*(\tau,x)
    \right)
    \left(
    f(\tau,x,r)
    -
    f_* (\tau,x,r_*)
    \right)
    \,
    \partial_x \Phi_{h,\varepsilon}(\tau,x) \,.
  \end{align}
  Plugging $\Phi_{h,\varepsilon}$ into~\eqref{eq:28}
  using~\eqref{eq:45} yields
  \begin{align*}
    \iint_{\mathcal{K}}
    \Big\{
    \left[
      \delta_h(\tau - \tau_1) - \delta_h(\tau - \tau_2)
      \right]
    \chi_\varepsilon (\tau,x)
    \modulo{r(\tau,x) - r_*(\tau,x)}
    \Big\}
    \d{x}\d{\tau}
     & \geq 0
     & \mbox{and}
    \\
    \int_0^T
    \left[
      \delta_h(\tau - \tau_1) - \delta_h(\tau - \tau_2)
      \right]
    \Big\{
    \int_{B_{\mathrm{N}  - \sigma \tau}}
    \modulo{r(\tau,x) - r_*(\tau,x)}
    \d{x}
    \Big\}
    \d{\tau}
     & \geq 0
     & \mbox{as } \varepsilon \to 0 \,.
  \end{align*}
  In the limit $h \to 0$, we obtain
  \begin{align*}
    \int_{B_{\mathrm{N}  - \sigma \tau_2}} \modulo{r(\tau_2, x) - r_*(\tau_2,x)} \d{x}
     & \leq
    \int_{B_{\mathrm{N}  - \sigma \tau_1}}\modulo{r(\tau_1, x) - r_*(\tau_1, x)} \d{x} \,;
    \\
    \int_{B_{\mathrm{N} -\sigma t}} \modulo{r(t, x) - r_*(t,x)} \d{x}
     & \leq
    \int_{B_{\mathrm{N}}}\modulo{r_o( x) - r_{o,*}(x)} \d{x}\qquad
    \mbox{as } \tau_1 \to 0 , \, t=\tau_2 \,,
  \end{align*}
  proving the claim.
\end{proofof}

\section{Functional Analysis Lemmas}
\label{sec:technique}

\begin{lemma}[Properties of periodic convolutions]
  Let $a \in \reali$ and
  $H \in \naturali \setminus\{0\}$.\label{lem:0}
  \begin{displaymath}
    \begin{array}{cccl@{\!\!}}
      \rho  \in \L\infty (\Su; \reali)
       & \eta \in \L1 (\reali; \reali)
       & \implies
       & \norma{\rho*\eta}_{\L\infty(\Su; \reali)}
      \leq
      \norma{\eta}_{\L1(\reali; \reali)}\,
      \norma{\rho}_{\L\infty(\Su; \reali)} \,.
      \\
      \rho  \in \L1 (\Su; \reali)
       &
      \begin{array}[c]{@{}c@{}}
        \eta \in \L\infty (\reali; \reali)
        \\
        \spt \eta \subseteq [a, a+H]
      \end{array}
       & \implies
       & \norma{\rho*\eta}_{\L\infty(\Su; \reali)}
      \leq
      H \,
      \norma{\eta}_{\L\infty(\reali; \reali)}\,
      \norma{\rho}_{\L1(\Su; \reali)} \,.
    \end{array}
  \end{displaymath}
\end{lemma}

\begin{proofof}{Lemma~\ref{lem:0}}
  The first inequality is a direct consequence of the
  standard properties of the convolution. To prove the second one,
  evaluate
  \begin{eqnarray*}
    \norma{\rho*\eta}_{\L\infty(\Su; \reali)}
    & =
    & \esssup_{x \in [0, 1]} \modulo{\int_\reali
      \widehat\rho(x-y) \, \eta(y) \d{y}}
    \\
    & \leq
    & \norma{\eta}_{\L\infty(\reali; \reali)} \,
    \esssup_{x \in [0, 1]} \int_a^{a+H} \modulo{\widehat\rho(x-y)} \dd y
    \\
    & \leq
    & H \, \norma{\eta}_{\L\infty(\reali; \reali)} \,
    \norma{\rho}_{\L1(\Su; \reali)} \,,
  \end{eqnarray*}
  completing the proof.
\end{proofof}

\begin{lemma}[Nonlinear translation estimate]
  \label{lem:BVper}
  For any $f \in \cBV{}(\Su; \reali)$ and any measurable
  $h \colon \Su \to \Su$, denoting by $\widetilde{h}$
  {
      a lift on $\reali$ of $h$ as in~\eqref{eq:16}
    },
  \begin{displaymath}
    \int_{\Su} \modulo{f\left(h(x)\right) - f(x)} \dd x
    \leq
    2 \, \esssup_{\reali} \modulo{\widetilde{h} (x) - x} \, \tv(f) \,.
  \end{displaymath}
\end{lemma}

\begin{proofof}{Lemma~\ref{lem:BVper}}
  Since $f \in \cBV{}(\Su; \reali)$, then $\widehat{f} \in \BVloc\left(
    \reali; \reali\right)$ and so it can be written as a sum of two monotone
  functions. Hence, $\widehat{f}$ is Borel measurable and thus
  $f$ is also Borel measurable. Consequently, the composition map $f \circ h$ is measurable.

  Fix $a \leq 0$.  Fix $\mathrm{N} \in \naturali \setminus \{0\}$ and
  call ${ F(x) \coloneqq \tv(\widehat f; \mathopen]a, x])}$ for
  all $x \geq a$. The function $F$ is non-decreasing and
  \begin{equation}\label{eq:15}
    \mbox{for all } x,y \in [a, +\infty\mathclose[ \mbox{ with } y \leq x
    \qquad
    \modulo{\widehat{f}(x) - \widehat{f}(y)}
    \leq
    F(x) - F(y) \, .
  \end{equation}
  Following~\cite[Lemma~2.3]{MR2347697}, observe that for all
  $c \geq 0 $
  \begin{align*}
    \int_a^{\mathrm{N}+a} \!\!\!\left( F(x + c) - F(x) \right) \dd x
     & =
    \text{meas}
    \left\{
    (x, y) \in \reali^2 \colon
    x \in [a, \mathrm{N}+a] \text{ and }
    y \in [F(x), F(x+c)]
    \right\}
    \\
     & =
    \int_0^{F(\mathrm{N} + a + c)}
    \text{meas}
    \left\{
    x \in [0,\mathrm{N}] \colon
    y \in \left[F(x) ,\, F(x + c) \right]
    \right\} \dd y
    \\
     & \leq
    \int_0^{
      {
          \tv\left(\widehat{f}; \mathopen]a, \mathrm{N} +a + c ]\right)
        }
    }
    c \dd y
    \\
     & =
    c \, {
        \tv\left(\widehat{f}; \mathopen]a, \mathrm{N} + a + c \mathclose]\right)
      }
    \\
     & =
    c \, \left(
    \mathrm{N} \tv(f) +
    {
      \tv\left(\widehat{f}; \mathopen]a, a+c\mathclose]
      \right)
    }
    \right)\!\!
    \tag*{
      \mbox{[By periodicity of $\widehat{f}$]}
    }
  \end{align*}
  Then, defining
  $a \coloneq \min\left\{\! \essinf_{\xi\in\reali}\left( \widetilde{h}
      (\xi)-\xi \right), 0 \right\}$ and
  $b \coloneq \max\left\{ \!\esssup_{\xi \in \reali}\left( \widetilde{h}
      (\xi)-\xi \right), 0 \right\}$,
  \begin{align*}
     & \int_{\Su} \modulo{f\left(h(x)\right) - f(x)} \d{x}
    \\
     & =
    \int_0^1
    \modulo{ \widehat{f}\left(\widetilde{h}(x)\right) - \widehat{f}(x)} \dd x
    \\
     & =
    \dfrac{1}{\mathrm{N}}
    \int_0^{\mathrm{N}}
    \modulo{ \widehat{f}\left(\widetilde{h}(x)\right) - \widehat{f}(x)} \dd x
    \\
     & \leq
    \dfrac{1}{\mathrm{N}}
    \int_0^{\mathrm{N}}
    \left[
      F\left(\max\left\{\widetilde{h} (x),x \right\} \right)
      -
      F\left(\min\left\{\widetilde{h} (x),x \right\} \right)
      \right]\dd x
     & \mbox{[By~\eqref{eq:15}]}
    \\
     & =
    \dfrac{1}{\mathrm{N}}
    \int_0^{\mathrm{N}}
    \left[
      F\left(x+\max\left\{\widetilde{h} (x)-x,0 \right\} \right)
      -
      F\left(x+\min\left\{\widetilde{h} (x)-x,0 \right\} \right)
      \right]\dd x
    \\
     & \leq
    \dfrac{1}{\mathrm{N}}
    \int_0^{\mathrm{N}}
    \left[
      F(x + b) - F (x+a)
      \right]\dd x
    \\
     & =
    \dfrac{1}{\mathrm{N}}
    \int_a^{a+\mathrm{N}}
    \left[ F(x + b - a) - F(x)\right]\dd x
    \\
     & \leq
    (b-a) \, \left(
    \tv(f) +
    \dfrac{1}{\mathrm{N}} \,
    {
      \tv\left(\widehat{f}; \mathopen]a, b\mathclose]
      \right)}
    \right)
  \end{align*}
  Then, letting $\mathrm{N} \to +\infty$ completes the proof.
\end{proofof}

\begin{lemma}[Characterization of $\cBV{1}$]
  \label{lem:last}
  The set $\cBV{1}(\Su; \reali^k)$ is the set of $\W11$ limits of
  $\W21$ bounded sequences in $\C2(\Su; \reali^k)$. Moreover, if
  $\rho \in \cBV{1}(\Su; \reali^k)$ there exists a sequence
  $\rho_n \in \C2(\Su; \reali^k)$ such that
  \begin{displaymath}
    \norma{\rho}_{\cBV{1}(\Su; \reali^k)}
    =
    \lim_{n \to \infty} \norma{\rho_n}_{\W21(\Su; \reali^k)}\,.
  \end{displaymath}
\end{lemma}

\begin{proofof}{Lemma~\ref{lem:last}}
  Consider any bounded sequence
  $\rho_n \in \left(\C2 \cap \W21\right)(\Su; \reali^k)$ converging in
  $\W11$ to a limit function $\rho$. Evidently,
  $\rho \in \W11(\Su; \reali^k)$.

  By assumption, for all $n \in \naturali$
  $\widehat \rho_n'' \in \Lloc1(\Su; \reali^k)$ and there exists
  $K >0$ such that
  $\sup_{n \in \naturali } \norma{\rho_n''}_{\L1(\Su; \reali)} \leq K
  $.  Call $\mu_n$ the associated Radon measure on $\reali$ and
  observe that for all compact set $\Omega \subset \reali $
  \begin{displaymath}
    \sup_{n \in \naturali }\norma{\mu_n} \left(\Omega\right)
    =
    \sup_{n \in \naturali}\int_\Omega \modulo{\widehat\rho_n''(x) }\d{x}
    \leq
    \left(\mbox{diam}(\Omega) +1 \right)
    \sup_{n \in \naturali}\norma{\rho_n''}_{\L1(\Su; \reali)}
    \leq
    \left(\mbox{diam}(\Omega) +1 \right) K \,.
  \end{displaymath}
  Due to the weak* compactness property of Radon measures
  (see~\cite[Corollary~1.60]{MR1857292}) there exists a Radon measure
  $\mu$ on $\reali$ such that, up to a subsequence, $\mu_n$ locally
  weakly* converges to $\mu$, i.e.
  \begin{displaymath}
    \int_\reali f(x) \d \mu  =
    \lim_{n \to + \infty} \int_\reali f (x) \d \mu_n
    =
    \lim_{n \to + \infty} \int_\reali f(x) \, \widehat \rho_n''(x) \d x
    \quad \mbox{for all } \, f \in \Cc0(\reali; \reali)\,.
  \end{displaymath}
  In particular, for all $f \in \Cc1(\reali; \reali^k)$, since
  $\rho_n' \to \rho'$ in $\L1(\Su; \reali)$,
  \begin{displaymath}
    \lim_{n \to + \infty} \int_\reali f(x) \, \widehat \rho_n''(x) \d x
    =
    -
    \lim_{n \to + \infty} \int_\reali f'(x) \, \widehat \rho_n'(x) \d x
    =
    -
    \int_\reali f'(x) \, \widehat \rho'(x) \d x \, .
  \end{displaymath}
  Hence, combining the previous expressions, one gets that
  \begin{displaymath}
    \int_\reali f(x) \d \mu
    =
    - \int_\reali f'(x) \, \widehat \rho'(x) \d x
    \quad
    \mbox{for all} \,
    f \in \Cc1(\reali; \reali^k),
  \end{displaymath}
  proving that $\widehat \rho' \in \BVloc(\reali; \reali^k)$ and,
  consequently, $\rho \in \cBV{1}(\Su; \reali^k)$.

  Conversely, fix $\rho \in \cBV{1}(\Su; \reali^k)$ and
  $\delta \in \Cc\infty (\reali; \reali_+)$ be such that
  $\spt \delta \subseteq [-1,1]$ and
  $\int_\reali \delta(\xi) \d{\xi} = 1$.  For any $n \in \naturali$ we
  define
  $
    \delta_n(x) \coloneqq n \delta \left(n x\right)$
  and $\rho_n \coloneqq \widecheck{\widehat{\rho}*\delta_n}$. Note
  that $\widehat{\rho}*\delta_n$ is periodic due to \Cref{lem:0}.  By
  standard convolution properties,
  $\widehat{\rho}_n \to \widehat{\rho}$ in
  $\Wloc11(\reali; \reali^k)$, hence $\rho_n \to \rho$ in
  $\W11(\Su; \reali^k)$.  Moreover, recalling the definition of
  $\tv(\rho)$ in~\eqref{eq:44}, by~\cite[Proposition 8.49]{Folland},
  one gets
  \begin{equation}
    \label{eq:68}
    \norma{\rho_n}_{\W11(\Su; \reali^k)}
    \leq
    \norma{\rho}_{\W11(\Su; \reali^k)}
    \quad \mbox{ and } \quad
    \norma{\partial_{xx}^2\rho_n}_{\L1(\Su; \reali^k)} \leq \tv(\rho')\, .
  \end{equation}
  The lower semicontinuity of the total variation with respect to the
  $\Lloc1$ topology and~\eqref{eq:68}, right, yield that
  $\tv(\rho') = \lim_{n \to \infty}
    \norma{\partial_{xx}^2\rho_n}_{\L1(\Su; \reali^k)}$, completing the
  proof.
\end{proofof}

\section*{Acknowledgments}
The authors were partly supported by the GNAMPA~2025 project
\emph{Modellli di Traffico, di Biologia e di Dinamica dei Gas Basati
  su Sistemi di Equazioni Iperboliche}. They also acknowledge the
PRIN~2022 project \emph{Modeling, Control and Games through Partial
  Differential Equations} (CUP~D53D23005620006), funded by the
European Union - Next Generation EU.  MG was partly supported by the
European Union-NextGeneration EU (National Sustainable Mobility Center
CN00000023, Italian Ministry of University and Research Decree
n.~1033-17/06/2022, Spoke 8) and by the Project funded under the
National Recovery and Resilience Plan (NRRP) of Italian Ministry of
University and Research funded by the European Union-NextGeneration
EU. Award Number: ECS\_00000037, CUP: H43C22000510001,
MUSA-Multilayered Urban Sustainability Action. CN was supported by
the PNRR project \emph{Ricerca DM 118/2023}, CUP~F13C23000360007.

\section*{Conflict  of Interest}
The authors declare no conflicts of interest in this paper.

  {\small

    \bibliographystyle{abbrv}

    \bibliography{wp_copy}

  }

\end{document}